\documentclass[12pt]{article}

\usepackage{verbatim,amsmath,amssymb,mathrsfs,amsthm,amsfonts,bm}
\usepackage{varwidth}
\usepackage{extarrows}
\usepackage{enumitem}
\usepackage{fancyhdr}
\usepackage{booktabs}
\usepackage[dvipsnames]{xcolor}
\usepackage{graphicx,wrapfig}
\usepackage{subfigure}
\usepackage{caption}
\usepackage{amsmath}
\usepackage{makecell}
\usepackage{multirow}
\usepackage{array}
\usepackage{setspace}
\usepackage{float}
\usepackage{cite}
\usepackage{color,mathtools}
\usepackage{epstopdf}
\usepackage{lipsum}
\usepackage{appendix}
\usepackage{algorithm}
\usepackage{algorithmic}
\usepackage{blindtext}
\usepackage{ulem}
\usepackage[left]{lineno}

\numberwithin{equation}{section}
\numberwithin{figure}{section}
\numberwithin{table}{section}
\allowdisplaybreaks

\newtheorem{theorem}{Theorem}[section]

\newtheorem{remark}{Remark}[section]
\newtheorem{example}{Example}[section]


\usepackage{hyperref}
\hypersetup{
   colorlinks=true, linkcolor=black,
   urlcolor=cyan, citecolor=black,
}

\begin{document}
\title{Well-balanced nodal finite volume WENO scheme with flux localization for hyperbolic balance laws}
\date{}

\author{Jiahui Zhang\thanks{School of Mathematical Sciences, University of Science and Technology of China, Hefei, Anhui 230026, P.R. China.  E-mail: zjh55@mail.ustc.edu.cn.}
	\and Yinhua Xia\thanks{School of Mathematical Sciences, University of Science and Technology of China, Hefei, Anhui 230026, P.R. China.  E-mail: yhxia@ustc.edu.cn. The research of Yinhua Xia was partially supported by Anhui Provincial Natural Science Foundation grant No. 2408085J004, Shanghai Institute for Mathematics and Interdisciplinary Sciences grant No. SIMIS-ID-2024-XD and NSFC grant No. 12271498.
	}
	\and Yan Xu\thanks{Corresponding author. School of Mathematical Sciences, University of Science and Technology of China, Hefei, Anhui 230026, P.R. China and Laoshan Laboratory, Qingdao 266237, P.R. China.  E-mail: yxu@ustc.edu.cn. The work of Yan Xu was partially supported by NSFC grant No. 12471347,  and Laoshan Laboratory (No. LSKJ202300305).
	}
}

\maketitle

\begin{abstract}
In this paper, we develop a general framework for constructing arbitrarily high-order well-balanced ``nodal" finite volume weighted essentially non-oscillatory (WENO) schemes, specifically tailored for hyperbolic balance laws. To handle complex systems involving nonconservative products and steady-state solutions with global integral quantities, we introduce a flux localization technique. This technique defines the numerical fluxes, incorporating the entire source term and nonconservative products, locally within each computational cell. By exactly preserving the nodal equilibrium state and conservatively updating the cell averages within a finite volume framework, the method avoids both reference-state reconstruction and the use of ODE solvers. The well-balanced property is achieved by performing high-order WENO reconstruction and employing a linear segment path in terms of equilibrium variables, rather than conservative ones. A rigorous theoretical analysis confirms that the scheme exactly preserves all the targeted stationary solutions. The effectiveness and versatility of the proposed method are demonstrated through a series of numerical experiments on various systems of hyperbolic balance laws.
\end{abstract}


\section{Introduction}\label{se:in}
We focus on a class of hyperbolic balance laws that have been extensively applied across various disciplines, including chemistry, biology, fluid dynamics, astrophysics, and meteorology. In the one-dimensional case, such systems can be written as:
\begin{equation}\label{ncs}
  \boldsymbol {u}_t + \boldsymbol{f}(\boldsymbol {u})_x +  \mathcal{G}(\boldsymbol {u})\boldsymbol {u}_x = \boldsymbol{r}(\boldsymbol {u})w_x + \boldsymbol{s}(\boldsymbol{u}),
\end{equation}
where $\boldsymbol {u}$ denotes the conservative variable, and $\boldsymbol{f}(\boldsymbol {u})$ is the physical flux. The terms $\boldsymbol{r}(\boldsymbol {u})$ and $\boldsymbol{s}(\boldsymbol{u})$ represent source terms arising from geometrical, biological, reactive, or other physical effects, with $w$ being a known function. 
The term $\mathcal{G}(\boldsymbol {u}) \boldsymbol {u}_x$ represents a nonconservative product, which cannot be expressed in divergence form. As a result, the system belongs to the broader class of nonconservative hyperbolic partial differential equations (PDEs).

Such PDE systems arise in a wide range of fluid models, including shallow water equations, multiphase flow, and gas dynamics. A key feature of these models is the existence of non-trivial steady-state solutions, where the non-zero flux gradient is exactly balanced by the source terms. For the system above, steady-state solutions $\boldsymbol{u}$ satisfy
\begin{equation}\label{steady}
  \boldsymbol{f}(\boldsymbol {u})_x +  \mathcal{G}(\boldsymbol {u})\boldsymbol {u}_x = \boldsymbol{r}(\boldsymbol {u})w_x + \boldsymbol{s}(\boldsymbol{u}).
\end{equation}
Due to the presence of source terms and nonconservative products, many standard numerical methods fail to capture such steady states accurately, even on highly refined meshes. This motivates the development of well-balanced numerical methods, which can preserve steady-state solutions at the discrete level and accurately resolve small perturbations around them using relatively coarse meshes. Such methods are referred to as ``well-balanced", a term first introduced by Bermudez and Vazquez \cite{bermudez1994upwind}. 

Over the past two decades, the development of numerical methods with desirable properties for systems of the form \eqref{ncs}, including the special case where $\mathcal{G}(\boldsymbol {u})=0$ , has become a highly active area of research. Significant mathematical effort has been devoted to constructing well-balanced schemes for hyperbolic systems featuring nonconservative products and/or source terms.
Par{\'e}s and collaborators proposed well-balanced finite volume methods \cite{castro2020well} and finite difference methods \cite{pares2021well} with arbitrary-order accuracy for balance laws, utilizing specially designed reconstruction operators. A key challenge in these methods lies in computing the stationary solution in each cell or node at every time step, which involves solving an ODE system subject to a prescribed average or point value over the computational domain.
Further generalizations were developed in \cite{gomez2021high,gomez2021collocation}, where the authors introduced a flexible framework capable of preserving a wide range of stationary solutions, not only hydrostatic equilibrium but also moving equilibrium, thus making the approach applicable to general systems of balance laws. Berberich et al. \cite{berberich2021high} designed high-order finite volume schemes for multi-dimensional balance laws under the assumption that steady-state solutions are provided in either analytical or discrete form.
Another influential technique is flux globalization, which addresses steady states involving spatial integrals. This approach, when combined with the path-conservative central-upwind (PCCU) method free of Riemann solvers, has been successfully applied to various systems, such as shallow water equations \cite{cao2022flux,diaz2019path,kurganov2023well,cao2023flux} and the more complex rotating shallow water magnetohydrodynamics \cite{chertock2024locally}. Ciallella et al. \cite{ciallella2023arbitrary} further extended this approach to achieve arbitrarily high-order accuracy using a finite volume WENO framework.
In the discontinuous Galerkin (DG) setting, Guerrero et al. \cite{guerrero2021well} proposed a well-balanced DG scheme by incorporating a local projection step with Runge-Kutta and ADER time integration. 
Mantri et al. \cite{mantri2021well} developed a well-balanced DG scheme for 2$\times$2 hyperbolic balance laws by solving a linear system for the equilibrium variables derived from global fluxes. This framework was further extended in \cite{mantri2024fully} via a DG spectral element method. Additionally, a novel well-balanced DG method for the rotating shallow water equations was introduced in \cite{zhang2025well}, also based on the flux globalization strategy. The concept of equilibrium-preserving spaces in DG methods has also been established for various hyperbolic balance laws, including the single-layer shallow water equations \cite{zhang2023moving}, the Ripa model and Euler equations with gravity \cite{2023arXivequilibrium}, and the two-layer shallow water equations \cite{zhang2024arxivwell}. For further developments on well-balanced methods across other hyperbolic systems, see \cite{dumbser2023wellbalanced,xu2025well,zhang2023structure,gaburro2021well,xing2014exactly}.

The finite volume method is a widely used approach for solving hyperbolic conservation laws in numerical simulations. One of its key advantages is the ability to handle unstructured meshes while preserving both high-order accuracy and conservation properties.
A major challenge in solving hyperbolic problems is the appearance of strong shocks, which can induce nonphysical oscillations known as Gibbs phenomena. To address this, special care must be taken in constructing non-oscillatory subcell reconstructions. Several such methods have been developed, including total variation diminishing (TVD) schemes \cite{leveque1992numerical} and the monotonic upwind scheme for conservation laws (MUSCL) \cite{van1979towards}. However, TVD schemes are inherently limited to second-order accuracy and tend to degenerate to first-order near extrema.
To overcome these limitations, significant efforts have been devoted to the development of high-order, high-resolution schemes. In 1987, Harten et al. \cite{harten1987uniformly} introduced the essentially non-oscillatory (ENO) scheme within the finite volume framework. This was followed by Shu and Osher's \cite{shu1988efficient} finite difference ENO formulation in 1988. The ENO approach chooses the smoothest stencil among multiple candidates, based on local difference quotients.
In 1994, Liu et al. \cite{liu1994weighted} proposed the weighted essentially non-oscillatory (WENO) scheme as a natural extension of ENO. Rather than relying solely on the smoothest stencil, WENO forms a convex combination of all candidate stencils, blending them using nonlinear weights. Jiang and Shu \cite{jiang1996efficient} further advanced this idea by introducing a new smoothness indicator, resulting in the classical WENO-JS scheme with higher-order accuracy.  Since then, numerous variants have been proposed to improve the accuracy and robustness of WENO reconstructions, especially near critical points.
For instance, WENO-M \cite{henrick2005mapped} and WENO-Z \cite{castro2011high} improve accuracy at smooth extrema.
Central WENO schemes \cite{levy1999central} offer simplicity and robustness by avoiding characteristic decomposition.
Hybrid WENO schemes \cite{li2010hybrid} combine multiple indicators to identify troubled cells, significantly reducing computational cost while improving shock resolution. 
For a comprehensive overview of the latest developments in WENO methodology, we refer readers to the recent review by Shu \cite{shu2020essentially}.

The main objective of this work is to design a well-balanced finite volume WENO (FV-WENO) scheme for a class of hyperbolic balance laws, aiming to achieve high-order accuracy while precisely preserving all the targeted stationary solutions in one dimension, including those not only described by simple algebraic relation, but also those involving spatially global integrals. In two dimensions, we restrict ourselves to quasi one-dimensional steady states. The key highlights are the following:
\begin{itemize}
\item {Nodal evolution with high-order well-balancing:}  The framework achieves an arbitrarily high-order well-balanced property by evolving nodal point values directly, instead of taking traditional cell averages as the primary degrees of freedom. The cell averages are instead computed dynamically as a linear combination of these adjacent nodal values. The method preserves the conservative finite volume structure and facilitates the preservation of target equilibria at the nodes.

\item {Flux localization technique for global extremes:}  A novel flux localization technique is introduced to handle complex systems involving nonconservative products and steady-state solutions with global integral quantities. This method incorporates the entire source term and nonconservative products locally within each computational cell, eliminating the need for reference-state reconstructions or ODE solvers.

\item {Equilibrium variable reconstruction and linear paths:} To precisely preserve all the targeted stationary solutions, the high-order WENO reconstruction is applied directly to equilibrium variables rather than conservative ones. Furthermore, scheme stability is guaranteed by utilizing linear segment paths to connect the left and right states of these equilibrium variables.
\end{itemize}


The outline of this paper is as follows. Section \ref{se:fvnv} introduces the basic concept of the nodal finite volume WENO method for a scalar hyperbolic conservation law. In Section \ref{se:fv}, the well-balanced nodal finite volume WENO scheme with flux localization is proposed for a class of one-dimensional hyperbolic balance laws. This method is then extended to the two-dimensional case in Section \ref{se:e2d}. Section \ref{se:nu} presents numerical results under various scenarios to demonstrate that the proposed nodal FV-WENO scheme is well-balanced, high-order accurate, and effective in shock-capturing. Finally, Section \ref{se:co} provides some concluding remarks.

\section{Nodal FV-WENO scheme}\label{se:fvnv}
To provide a clearer and more comprehensive introduction to the novel well-balanced finite volume scheme for hyperbolic balance laws in the next section, we begin by briefly outlining the core idea of the proposed nodal finite volume scheme for one-dimensional scalar hyperbolic conservation laws
\begin{equation}\label{scalar}
  U_t + K(U)_x = 0.
\end{equation} 

\subsection{Scheme formulation}
We assume that the computational domain $\mathcal{D}$ is partitioned into $N$ uniform cells denoted by $I_j$, such that
$$\mathcal {D} = \bigcup_{j=1}^{N} I_j = \bigcup  _{j=1}^{N}[x_{j-\frac{1}{2}}, x_{j+\frac{1}{2}}].$$
The $x_{j\pm \frac{1}{2}}$ represent the left and right interfaces of cell $I_j$.  The cell width is defined as $\Delta x = x_{j+\frac{1}{2}}-x_{j-\frac{1}{2}}$, and the center of each cell is given by $x_j = \frac{1}{2}(x_{j-\frac{1}{2}} + x_{j+\frac{1}{2}})$. 

In traditional finite volume schemes, the primary computational variables are the cell averages $\overline{{U}}_j$, which approximate the average of the exact solution over cells $I_j$, defined as $$\overline{{U}} (x_j,t) = \frac{1}{\Delta x}\int_{x_{j-{1}/{2}}}^{x_{j+{1}/{2}}} {U}(x,t) \ dx.$$ 
{Traditional high-order well-balanced finite volume schemes usually require an ODE solver with prescribed cell averages in order to recover the generalized moving steady state solution. In this work, in} contrast to most existing finite volume methods, we adopt an arbitrary nodal value within each cell $I_j$ as the degree of freedom, rather than computing the cell average directly. For simplicity and ease of implementation, we choose the cell midpoint as the representative node. In this setting, the cell averages are then reconstructed as linear combinations of the nodal values from neighboring cells. {Our proposed scheme remains entirely within the finite volume framework.}
Specifically, the nodal values ${U}_j$ are used to approximate the pointwise solution ${U}(x_j,t)$.
The resulting finite volume scheme takes the form
\begin{equation}\label{scheme_scalar}
  \dfrac{d}{dt}\int_{I_j} {U} \ dx= - \left(\widehat{{K}}_{j+\frac{1}{2}}   - \widehat{{K}}_{j-\frac{1}{2}}\right):=\text{RHS}_j({U}).
\end{equation}
The integral $\int_{I_j} {U} \ dx$ on the left-hand side of \eqref{scheme_scalar} denotes the numerical cell integral, which is computed as a fixed linear combination of the numerical nodal values.
Here, $\widehat{{K}}_{j\pm\frac{1}{2}} = \widehat{K}(U_{j+\frac{1}{2}}^-,U_{j+\frac{1}{2}}^+)$ denotes the numerical flux, where $U_{j+\frac{1}{2}}^{\pm}$ are the high-order pointwise approximations of the solution $U(x_{j+\frac{1}{2}},t)$, reconstructed from the nodal values ${U}_j$ using a WENO reconstruction procedure. 
This is different from the classical finite volume schemes, which reconstruct from the cell averages to cell boundary point values, see Appendix \ref{a1} for details.

To ensure the stability of the scheme, we employ the simple and robust Lax-Friedrichs flux, given by
\begin{equation}\label{LF_scalar}
  \widehat {K}(U_{j+\frac{1}{2}}^-,U_{j+\frac{1}{2}}^+) = \frac{1}{2}(K(U_{j+\frac{1}{2}}^-)+K(U_{j+\frac{1}{2}}^+) - \alpha(U_{j+\frac{1}{2}}^+-U_{j+\frac{1}{2}}^-))
\end{equation}
where $\alpha = \max\limits_U |\frac{\partial K}{\partial U}|$ is chosen either globally or locally depending on the problem setting.

To obtain the fully discrete scheme, the semi-discrete formulation \eqref{scheme_scalar} is typically integrated in time using high-order strong stability preserving Runge-Kutta (SSP-RK) methods or multistep time discretization techniques, as discussed in \cite{gottlieb2001strong}. In this work, we use the third-order SSP-RK method for time discretization, which is given by the following steps
  \begin{align}
  & \int_{I_j} {U}^{(1)}  \ dx 
  = \int_{I_j} {U}^n  \ dx 
  + \Delta t  \ \text{RHS}_j({U}^n), \nonumber\\ 
  & \int_{I_j} {U}^{(2)} \ dx 
  = \dfrac{3}{4}\int_{I_j} {U}^{n}  \ dx 
  + \dfrac{1}{4}\left(\int_{I_j}{U}^{(1)}  \ dx 
  + \Delta t  \ \text{RHS}_j({U}^{(1)}) \right),\label{3rk_no}\\ 
  & \int_{I_j} {U}^{n+1}  \ dx 
  = \dfrac{1}{3}\int_{I_j} {U}^{n}  \ dx 
  + \dfrac{2}{3}\left(\int_{I_j} {U}^{(2)}  \ dx 
  + \Delta t  \ \text{RHS}_j({U}^{(2)}) \right), \nonumber
  \end{align}
where $\text{RHS}_j({U})$ represents the spatial operator in \eqref{scheme_scalar}. 

One important consideration is that, although we use the nodal values $\{U_j\}$ at cell midpoints as the degrees of freedom, the cell averages still need to be updated within the finite volume framework. To compute the cell average, we interpolate a high-degree polynomial, incorporating information from both the current and neighboring cells. Specifically, the average is expressed as $$\int_{I_j}{U} \ dx:= \sum_{i\in S_j} a_i {U}_i,$$ where $S_j = \{ \dots, j-1,j,j+1,\dots\}$ is the set of all index tuples of cells, which depends on the precision of the scheme.  

In particular, we provide a detailed implementation of the resulting scheme that achieves fifth-order accuracy and accommodates all the boundary conditions. 
We interpolate a fifth-order polynomial within cell $I_j$ using the nodal values $\{{U}_{j-2},{U}_{j-1},{U}_{j},{U}_{j+1},{U}_{j+2}\}$. 
For cells whose interpolation stencil does not intersect the boundary, namely $j = 3,\cdots,N-2$, the corresponding cell average can then be computed using the following formula, 
\begin{equation}\label{tem_c}
  \int_{I_j} {U} \ dx  := 
  \Delta x\left( \dfrac{-17}{5760}{U}_{j-2} +  \dfrac{77}{1440}{U}_{j-1} + \dfrac{863}{960}{U}_{j} + \dfrac{77}{1440}{U}_{j+1} - \dfrac{17}{5760}{U}_{j+2}\right).
\end{equation}
Near the boundaries, one-sided stencils are employed. The right-biased stencil $\{U_1,\ldots,U_5\}$ is used for cells $I_1,I_2$, whereas the left-biased stencil $\{U_{N-4},\ldots,U_N\}$ is adopted for cells $I_{N-1},I_N$. These cases can be expressed as follows,
\begin{equation}\label{tem_lr}
  \begin{aligned}
    &\int_{I_1} {U} \ dx  := 
    \Delta x\left( \dfrac{6463}{5760}{U}_{1} -  \dfrac{523}{1440}{U}_{2} + \dfrac{383}{960}{U}_{3} - \dfrac{283}{1440}{U}_{4} + \dfrac{223}{5760}{U}_{5}\right),\\
    &\int_{I_2} {U} \ dx  := 
    \Delta x\left( \dfrac{223}{5760}{U}_{1} +  \dfrac{1337}{1440}{U}_{2} + \dfrac{23}{960}{U}_{3} + \dfrac{17}{1440}{U}_{4} - \dfrac{17}{5760}{U}_{5}\right),\\
    &\int_{I_{N-1}} {U} \ dx  := 
    \Delta x\left( \dfrac{-17}{5760}{U}_{N-4} +  \dfrac{17}{1440}{U}_{N-3} + \dfrac{23}{960}{U}_{N-2} + \dfrac{1337}{1440}{U}_{N-1} + \dfrac{223}{5760}{U}_{N}\right),\\
    &\int_{I_{N}} {U} \ dx  := 
    \Delta x\left( \dfrac{223}{5760}{U}_{N-4} -  \dfrac{283}{1440}{U}_{N-3} + \dfrac{383}{960}{U}_{N-2} - \dfrac{523}{1440}{U}_{N-1} + \dfrac{6463}{5760}{U}_{N}\right).
  \end{aligned}
\end{equation}
Substituting these expressions  \eqref{tem_c} and \eqref{tem_lr} for the cell average back into the time discretization \eqref{3rk_no}, the nodal values ${U}_j, j=1,\cdots,N$, which serve as the degrees of freedom, can be obtained by solving a system of linear equations. To demonstrate this process, we use the Euler forward time discretization as an example, referring to the first equation in \eqref{3rk_no}.

We define the vectors
\begin{align*} \widetilde{{U}} = ({U}_1,\cdots,{U}_{N})^T, \quad  \widetilde{\text{RHS}} = (\text{RHS}_1,\cdots,\text{RHS}_{N})^T\end{align*}
which represent the vector of degrees of freedom for the conservative variables and the spatial operator from \eqref{scheme_scalar}, respectively. The time discretization can now be expressed as the following linear system of equations
\begin{align}\label{li_s}
   M \widetilde{{U}}^{n+1} =  M \widetilde{{U}}^{n} +\dfrac{\Delta t}{\Delta x} \ \widetilde{\text{RHS}},
\end{align}
 where $M$ is the mass matrix given by:
\begin{equation*}
   M = \begin{bmatrix}
    \frac{6463}{5760}& \frac{-523}{1440} & \frac{383}{960} & \frac{-283}{1440} & \frac{223}{5760} &  &  &  &   \\
    \frac{223}{5760} & \frac{1337}{1440} & \frac{23}{960} & \frac{17}{1440} & \frac{-17}{5760} &  &  &  &     \\
    \frac{-17}{5760} & \frac{77}{1440} & \frac{863}{960} & \frac{77}{1440} & \frac{-17}{5760} &  &  &  &     \\
     & \ddots & \ddots & \ddots & \ddots & \ddots &  &  &     \\
     &  & \ddots & \ddots & \ddots & \ddots & \ddots &  &     \\
     &  &  & \ddots & \ddots & \ddots & \ddots & \ddots &     \\
     &  &  &  &  \frac{-17}{5760} & \frac{77}{1440} & \frac{863}{960} & \frac{77}{1440} & \frac{-17}{5760} \\
     &  &  &  &  \frac{-17}{5760} & \frac{17}{1440} & \frac{23}{960} & \frac{1337}{1440} & \frac{223}{5760} \\
     &  &  &  &  \frac{223}{5760} & \frac{-283}{1440} & \frac{383}{960} & \frac{-523}{1440} & \frac{6463}{5760}
  \end{bmatrix}.
\end{equation*}
Owing to the fact that the mass matrix $M$ is sparse and remains constant at each time step, its inverse can be precomputed once and reused throughout the computation. Consequently, the computational cost of solving the system is relatively low.

\begin{remark}
For periodic boundary conditions, the interior relation \eqref{tem_c} is applied uniformly to all cells, with indices interpreted in the periodic sense, thereby eliminating the need for the one-sided formulas in \eqref{tem_lr} and rendering $M$  cyclic pentadiagonal.

\end{remark}


\begin{remark}
  The nodal FV-WENO scheme introduced in this section is specifically designed to satisfy the well-balanced property for balance laws, rather than serving as a shock-capturing method for conservation laws. The proposed nodal FV-WENO scheme shares with WCNS-type methods \cite{deng2000developing} the use of WENO interpolation to reconstruct interface states, the underlying discretizations are different. WCNS schemes are finite difference methods in which nodal unknowns are evolved through compact approximations of flux derivatives. In contrast, our method retains a finite volume formulation: cell averages are obtained by integrating polynomial interpolants constructed from the nodal values, and their evolution follows the finite volume framework. Nevertheless, the constructed scheme also exhibits robust behavior for problems involving discontinuous solutions. Since shock-capturing performance is not the main objective of this work, we present only two representative tests, based on the Burgers and Buckley-Leverett equations, in Appendix \ref{c1}. 
\end{remark}

\subsection{Conservation and linear stability analysis}
In this subsection, we establish the conservation property of the scheme \eqref{scheme_scalar}, and also analyze the linear stability of the fifth-order nodal FV-WENO scheme for the linear advection equation with periodic boundary conditions.

\subsubsection{Conservation}
Although the nodal values are used as the computational degrees of freedom, the scheme \eqref{scheme_scalar} updates the cell averages in the standard finite volume form. The nodal FV-WENO scheme can be written as
\begin{equation*}
  \frac{d}{dt}\sum_{i\in S_j} a_i U_i = -\left(\widehat K_{j+\frac{1}{2}} - \widehat K_{j-\frac{1}{2}}\right):=\mathrm{RHS}_j(U).
\end{equation*}
Local conservation follows from the use of the same numerical flux
$\widehat K_{j+\frac12}$ at the interface shared by the two neighboring cells. Summing over all cells, the interior numerical fluxes cancel pairwise. Hence, for periodic boundary conditions,
\begin{equation*}
  \frac{d}{dt} \sum_j U_j = 0.
\end{equation*}
Equivalently, in matrix form, the discrete total mass determined by the corresponding mass matrix is exactly conserved.


\subsubsection{Fourier analysis}
We consider the linear advection problem with periodic boundary conditions:
\begin{equation*}
  U_t+aU_x=0, \quad a>0.
\end{equation*}
The first-order nodal FV-WENO scheme reduces to the traditional first-order finite volume scheme, since $\overline U_j = U_j$, $U_{j+1/2}^- = U_j$ and $U_{j+1/2}^+ = U_{j+1}$. It has the same CFL analysis as the traditional first-order FV scheme, namely, $\frac{a\Delta t}{\Delta x}\le 1$. Therefore, the following analysis is restricted to fifth-order spatial discretization with linear weights and third-order SSP-RK time discretization.

To study the semi-discrete stability, we consider Fourier modes of the form
\begin{equation}\label{Fourier_mode}
  U_j(t)=\widehat U(\kappa,t)e^{\mathrm{i}\kappa x_j}.
\end{equation}
Denote $\theta = \kappa\Delta x$, and from the relation \eqref{tem_c}, we have $\overline U_j=m(\theta)U_j$, with
\begin{equation*}
  m(\theta)=\frac{863}{960} + \frac{77}{720}\cos\theta - \frac{17}{2880}\cos(2\theta).
\end{equation*}
It is easy to evaluate that $$\frac{283}{360}\le m(\theta)\le 1.$$
Substituting the linear reconstruction from Appendix \ref{a1} into the semi-discrete scheme \eqref{scheme_scalar} yields
\begin{equation*}
  \frac{d\widehat U}{dt} = \lambda(\theta)\widehat U,
\end{equation*}
where
\begin{equation*}
  \lambda(\theta) = -\frac{a}{\Delta x}\frac{(1-e^{-\mathrm{i}\theta})r(\theta)}{m(\theta)},
\end{equation*}
and   
\begin{equation*}
  r(\theta) = \frac1{128} \left(3e^{-2\mathrm{i}\theta} -20e^{-\mathrm{i}\theta} +90+60e^{\mathrm{i}\theta} -5e^{2\mathrm{i}\theta}\right).
\end{equation*}
Some tedious manipulation gives 
\begin{equation*}
  \operatorname{Re}\lambda(\theta) = -\frac{a}{\Delta x} \frac{3(1-\cos\theta)^3}{32m(\theta)} \le0.
\end{equation*}
Thus the semi-discrete nodal FV-WENO scheme is spectrally stable.

\begin{figure}[htbp]
  \centering
  \includegraphics[width= 5.cm,scale=1]{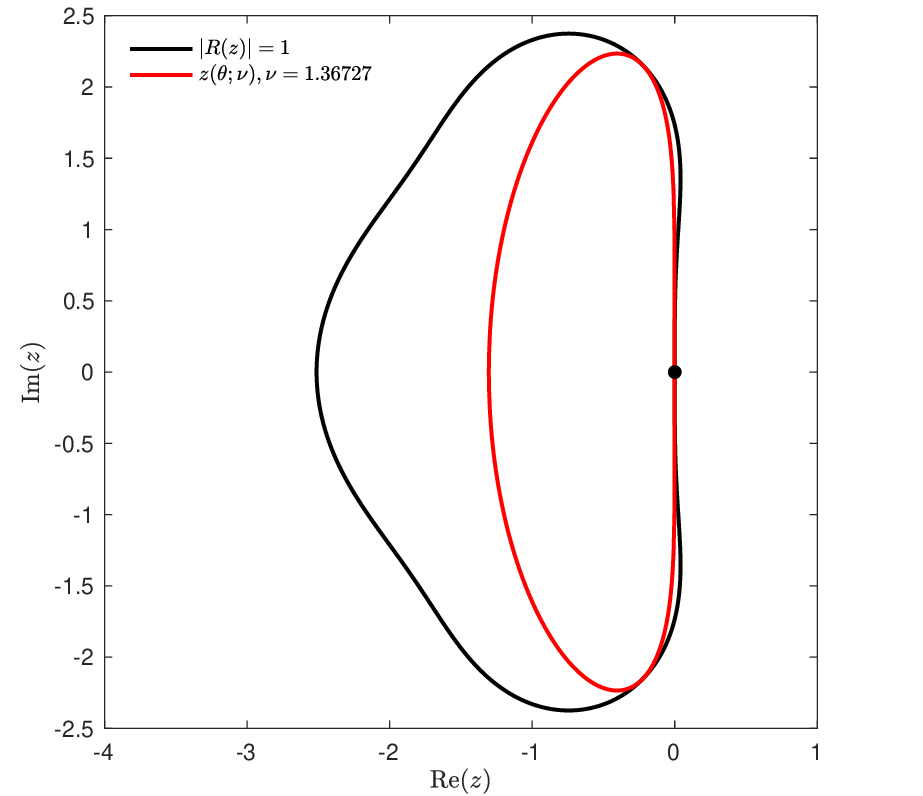}
  \caption{The boundary of the SSP-RK3 absolute stability region, $|\mathcal R(z)|=1$ (black curve); the Fourier-symbol locus at the critical CFL number $\nu =  1.36727$ (red curve).}\label{fig:cfl}
\end{figure}

To the fully-discrete scheme, let $U_j^n = \widehat U^n(\kappa)e^{\mathrm{i}\kappa x_j}$, the amplification factor is
\begin{equation*}
  \widehat U^{n+1}= G(\theta)\widehat U^n, \qquad G(\theta)=R\bigl(\Delta t\lambda(\theta)\bigr),
\end{equation*}
where 
\begin{equation*}
  \mathcal{R}(z) = 1+z+\frac{z^2}{2}+\frac{z^3}{6}.
\end{equation*}
Thus the fully-discrete scheme is linearly stable provided that 
\begin{equation}\label{stability_con}
  \max\limits_{\theta}|\mathcal R(\Delta t \lambda(\theta))| \le 1.
\end{equation} 
To more clearly visualize the corresponding CFL restriction, in Fig. \ref{fig:cfl}, we plot the Fourier-symbol locus $\{\Delta t\lambda(\theta):\theta\in[-\pi,\pi]\}$ together with the boundary of the SSP-RK3 absolute stability region, $|\mathcal R(z)|=1$. The critical CFL number is determined as the largest value for which the entire locus remains within the absolute stability region. Thus the CFL condition reads as $\nu = \frac{a\Delta t}{\Delta x}\le 1.36727$.

For comparison, 
Shadab et al. \cite{shadab2019fifth} showed a maximum linear CFL number of $1.44$ for the fifth-order finite volume WENO scheme in Cartesian coordinates with third-order SSP-RK time discretization. Thus, under the same linear advection setting, the proposed scheme admits a comparable, albeit slightly more restrictive, linear CFL limit.

\begin{remark}
  Under the periodic boundary conditions, the discrete Fourier modes form a basis of eigenvectors of the matrix $M$, and the corresponding eigenvalues $m(\theta)$ are nonzero, therefore the matrix $M$ is invertible.
\end{remark}

\section{Well-balanced nodal FV-WENO scheme}\label{se:fv}
In this section, we introduce a high-order well-balanced nodal FV-WENO scheme for hyperbolic balance laws, designed to preserve all the targeted stationary solutions in one dimension. Since the same third-order SSP-RK time discretization \eqref{3rk_no} applies to the fully-discrete scheme, we focus on the semi-discrete formulation and limit the discussion to the one-dimensional case.

\subsection{The well-balanced semi-discrete scheme with local flux}
Drawing inspiration from the framework presented in \cite{kurganov2023well}, but diverging in key aspects, we propose a well-balanced nodal FV-WENO scheme within the context of the flux localization technique. 
In most existing finite volume scheme, the well-balanced is achieved by reconstructing the perturbations around a locally recovered reference equilibrium. For generalized moving steady states, this recovery may require solving a local ODE problem subject to prescribed cell averages. In contrast, the nodal FV-WENO method adopts nodal values as degrees of freedom and preserves the target stationary states at the nodes. The corresponding cell averages are obtained from the nodal values and updated by a conservative finite volume scheme. This approach maintains the conservative finite volume update of the cell averages and avoids a separate local recovery of the reference equilibrium, thus providing a natural mechanism for preserving the target stationary states.
Specifically, we first integrate the trivially satisfied equation  $w_t=0$ into the system \eqref{ncs}, and then rewrite the system as
\begin{equation}\label{ncs_l}
  \boldsymbol{U}_t + \boldsymbol{F}(\boldsymbol{U})_x + \mathcal{B}(\boldsymbol {U})\boldsymbol {U}_x = \boldsymbol{S}(\boldsymbol{U}).
\end{equation} 
For notational simplicity, we augment the conservative variables $\boldsymbol{u}$ by incorporating the known function $w$, resulting in the vector $\boldsymbol{U} = (\boldsymbol{u}, w)^T$,  where the superscript ${}^T$ denotes the transpose. The structures of the flux function $ \boldsymbol{F}(\boldsymbol{U})$, the matrix $\mathcal{B}(\boldsymbol {U})$ and the source term $\boldsymbol{S}(\boldsymbol{U})$ are defined as follows:
\begin{equation*}
  \boldsymbol{F}(\boldsymbol{U}) = \left[\begin{array}{c}
          \boldsymbol{f}(\boldsymbol{u}) \\
          0
        \end{array}\right], \quad
  \mathcal{B}(\boldsymbol {U}) = \left[\begin{array}{cc}
    \mathcal{G}(\boldsymbol {u}) & -\boldsymbol{r}(\boldsymbol{u})\\
    0& 0 
  \end{array}\right],\quad \boldsymbol{S}(\boldsymbol{U}) = \left[\begin{array}{c}
    \boldsymbol{s}(\boldsymbol{u}) \\
    0
  \end{array}\right].
\end{equation*} 
Subsequently, we express this system in a quasi-conservative form as
\begin{equation}\label{quasi_ncs}
  \boldsymbol{U}_t + \boldsymbol{K}(\boldsymbol{U})_x = 0,
\end{equation}
with the flux 
\begin{equation}\label{glo_flux}
  \boldsymbol{K}(\boldsymbol{U}) = \boldsymbol{F}(\boldsymbol{U}) - \boldsymbol{R}(\boldsymbol{U}),  
  \quad  \boldsymbol{R}(\boldsymbol{U}) = \int_{\hat x}^{x} \boldsymbol{S}(\boldsymbol {U}(\xi,t))-\mathcal{B}(\boldsymbol {U}(\xi,t))\boldsymbol {U}(\xi,t)_{\xi}  \ d\xi.
\end{equation}
Here $\hat x$ is an arbitrary constant.

To preserve more intricate moving equilibrium within the framework of a finite volume scheme, especially for high-order schemes, most existing well-balanced methods require either the recovery of the reference solution or the use of an ODE solver to address equation \eqref{steady}. Inspired by \cite{zhang2025well}, and to overcome these challenges, we adopt the strategy outlined in Section \ref{se:fvnv}, denoting the nodal values $\boldsymbol{U}_j$ to approximate the point values of $\boldsymbol{U}(x_j,t)$. This approach is designed to ensure that the steady state is preserved at all midpoints. The finite volume scheme is given by
\begin{equation}\label{scheme}
  \dfrac{d}{dt}\int_{I_j}\boldsymbol{U} \ dx:= \text{RHS}_j(\boldsymbol{U}) = - \left(\widehat{\boldsymbol{K}}^{j,\text{mod}}_{j+\frac{1}{2}}   - \widehat{\boldsymbol{K}}^{j,\text{mod}}_{j-\frac{1}{2}}\right).
\end{equation}
With flux localization and modification, the numerical fluxes $\widehat{\boldsymbol{K}}^{j,\text{mod}}_{j\pm \frac{1}{2}} $ will be defined in Sec. \ref{sec:flux}.

\subsubsection{Transformation of the conservative and equilibrium variables}
Before delving into the specifics of the well-balanced treatment of our scheme, we first present the conversion between the sets of conservative and equilibrium variables. We denote  $\boldsymbol{v}$ as the equilibrium variables for the system \eqref{ncs}, which remain constant in the exact steady state.

On one hand, we express the transformation from conservative variables $\boldsymbol{u}$ to equilibrium variables $\boldsymbol{v}$ as $\boldsymbol{v} = \mathbb{V}(\boldsymbol{u},w)$. For most hyperbolic balance laws, the equilibrium variables $\boldsymbol{v}$ are given by a simple algebraic relation, which can be formulated as a function of $\boldsymbol{u}$ and $w$ such as $\boldsymbol{v} = \boldsymbol{g}_1(\boldsymbol{u},w)$. In particular, for equations involving friction or Coriolis forces, the steady-state solutions exhibit a global, space-integral quantity, which can be expressed as 
\begin{equation}\label{rela}
  \boldsymbol{v} = \boldsymbol{g}_1(\boldsymbol{u},w) + \mathcal{Q}, \quad \mathcal{Q} = \int_{\hat x}^{x} \boldsymbol{g}_2(\boldsymbol{u}) \ d \xi.
\end{equation}
Regardless of the circumstances, the transformation function $\boldsymbol{v} = \mathbb{V}(\boldsymbol{u},w)$ can be directly obtained from the explicit expression. In this paper, we use this latter form to present the well-balanced treatment, as the former can be viewed as a special case where $\boldsymbol{g}_2 = 0$. 

To ensure that the equilibrium values are preserved accurately at all midpoints throughout the calculation, we evaluate the corresponding equilibrium variables as
\begin{equation}\label{eq_v}
  \boldsymbol{v}_j = \boldsymbol{g}_1(\boldsymbol{u}_j,w_j) + \mathcal{Q}_j, \ \mathcal{Q}_j = \int_{\hat x}^{x_j} \boldsymbol{g}_2(\boldsymbol{u}) \ d \xi, \ j = 1,\cdots,N,
\end{equation}
where $\mathcal{Q}_j$ is the global term which can be calculated recursively as
\begin{equation}\label{q_re}
  \mathcal{Q}_j = \mathcal{Q}_{j-\frac{1}{2}} + \int_{x_{j-\frac{1}{2}}}^{x_j} \boldsymbol{g}_2(\boldsymbol{u}) \ d \xi.
\end{equation}
By taking $\hat x = x_{\frac{1}{2}}$, we have $\mathcal{Q}_{\frac{1}{2}} = 0$. 

The well-balanced property is ensured by interpolating a high-order polynomial within the cell $I_j$ based on the point values $\{\boldsymbol{g}_2(\boldsymbol{u}_i)\}_{i\in S_j}$ in the stencil $S_j$. The integral is then evaluated numerically using a quadrature formula that preserves the global order of accuracy. Thus, the proposed scheme is, in general, exactly well-balanced only with respect to the corresponding discrete steady state defined through $Q_j$, rather than the exact continuous steady state defined through $Q$, unless the quadrature used to compute $Q_j$ is exact for the relevant integrand. We will show that the semi-discrete and fully discrete schemes preserve this ``discrete" steady state exactly.

On the other hand, given $\boldsymbol{v}$, the inverse transformation is given by $\boldsymbol{u} = \mathbb{U}(\boldsymbol{v},w)$. Due to the nonlinear nature of the functional relationship, it must be solved using a nonlinear iterative solver, such as Newton's method. In the numerical examples presented below, we set the tolerance for the Newton solver to TOL = $10^{-14}$.

\subsubsection{Well-balanced treatment}\label{sec:flux}
Now, we define the numerical fluxes $\widehat{\boldsymbol{K}}^{j,\text{mod}}_{j\pm \frac{1}{2}} $ in \eqref{scheme}.

{\noindent \textbf{(a) Flux modification}}
\smallskip

The hydrostatic reconstruction strategy introduced in \cite{audusse2004fast} to achieve the well-balanced property serves as our inspiration. To preserve the well-balanced property, we make a slight modification to the Lax-Friedrichs numerical flux \eqref{LF_scalar}, which is given by
  \begin{align}\label{LF_mod_fv}
    \widehat{\boldsymbol{K}}^{j,\text{mod}}_{j\pm\frac{1}{2}} &= \dfrac{1}{2}\left(\boldsymbol{K}^j(\boldsymbol{U}_{j\pm\frac{1}{2}}^-) +\boldsymbol{K}^j(\boldsymbol{U}_{j\pm\frac{1}{2}}^+) - \alpha_{j\pm\frac{1}{2}} (\boldsymbol{U}_{j\pm\frac{1}{2}}^{*,+} - \boldsymbol{U}_{j\pm\frac{1}{2}}^{*,-}) \right),
  \end{align}
with $\alpha_{j\pm\frac{1}{2}}$ defined as
\begin{equation}\label{LF_alpha_fv}
  \alpha_{j\pm\frac{1}{2}} = \max\limits_{\boldsymbol{U}=\boldsymbol{U}_{j\pm\frac{1}{2}}^{\pm}} \{|\lambda_1(\boldsymbol{U})|,\dots,|\lambda_l(\boldsymbol{U})|\},
\end{equation}
and $\lambda_l$ are the eigenvalues of the Jacobi matrix $ \frac{\partial \boldsymbol{F}}{\partial \boldsymbol{U}}+\mathcal{G}(\boldsymbol {u})$. The maximum is taken locally in our computation.  

Moreover, we define the cell interface values and their modification to be
\begin{equation}\label{u_modify_fv}
  \boldsymbol{U}_{j+\frac{1}{2}}^{\pm} 
  = \left(\mathbb{U}(\boldsymbol{v}_{j+\frac{1}{2}}^{\pm},w_{j+\frac{1}{2}}^{\pm}),w_{j+\frac{1}{2}}^{\pm}\right)^T, \quad
  \boldsymbol{U}_{j+\frac{1}{2}}^{*,\pm} 
  = \left(\mathbb{U}(\boldsymbol{v}_{j+\frac{1}{2}}^{\pm},w_{j+\frac{1}{2}}^{*}),w_{j+\frac{1}{2}}^{*}\right)^T,
\end{equation}
and
\begin{equation}\label{b_modify_fv}
  w_{j+\frac{1}{2}}^{*} = \min(w_{j+\frac{1}{2}}^{-}, w_{j+\frac{1}{2}}^{+}).
\end{equation}
To preserve the equilibrium state, the values $\boldsymbol{v}_{j+\frac{1}{2}}^{\pm}$ in \eqref{u_modify_fv} are obtained from the nodal values $\boldsymbol{v}_j = \mathbb{V}(\boldsymbol{u}_j,w_j)$ via a WENO reconstruction procedure outlined in Appendix \ref{a1}. For robustness and to reduce oscillations at strong discontinuities, we prefer using the local characteristic decomposition \cite{shu1998essentially} over the component-wise version. Since the reconstruction is applied to the local characteristic fields of the equilibrium variables, similar to the approach in \cite{zhang2023moving}, the Jacobian matrix is replaced by
$$\left(\dfrac{\partial {\boldsymbol {v}}}{\partial\boldsymbol u} \right)
\left(\dfrac{\partial \boldsymbol{f}}{\partial \boldsymbol{u}} +\mathcal{G}(\boldsymbol {u}) \right)
\left(\dfrac{\partial {\boldsymbol {v}}}{\partial\boldsymbol u}\right)^{-1}. $$

\bigskip
{\noindent \textbf{(b) Flux localization}}
\smallskip

From the ODE system \eqref{steady}, which the non-trivial stationary solutions satisfy, and its equivalent form \eqref{quasi_ncs}, it is straightforward to show that the steady state can alternatively be expressed as $\boldsymbol{K}(\boldsymbol{U}) = \text{constant}$. From \eqref{glo_flux}, it is evident that $\boldsymbol{R}$ contains an integral quantity that depends on $\hat x$. 
In most existing well-balanced schemes with flux globalization, $\hat x$ is typically chosen as the left boundary of the computational domain, leading to a global, recursive calculation of this integral. In contrast, the key to maintaining the well-balanced property in our scheme lies in evaluating the local flux function, $\boldsymbol{K}^j(\boldsymbol{U}_{j\pm\frac{1}{2}}^{\pm})$, rather than the global flux used in \cite{zhang2025well}. 

All definitions and calculations are performed in cell $I_j$, and we use the superscript ${}^j$ to indicate this. The local flux is defined as
\begin{equation}\label{Rpm_fv} 
  \boldsymbol{K}^{j}(\boldsymbol{U}_{j+\frac{1}{2}}^{\pm}) = \boldsymbol{F}(\boldsymbol{U}_{j+\frac{1}{2}}^{\pm}) - \boldsymbol{R}_{j+\frac{1}{2}}^{j,\pm},\quad
 \boldsymbol{K}^j(\boldsymbol{U}_{j-\frac{1}{2}}^{\pm}) = \boldsymbol{F}(\boldsymbol{U}_{j-\frac{1}{2}}^{\pm}) - \boldsymbol{R}_{j-\frac{1}{2}}^{j,\pm}.
\end{equation}

To handle nonconservative products in scenarios involving discontinuous bottom topographies, we follow the general framework presented in \cite{zhang2024arxivwell}, which allows us to compute the nonconservative products in a well-balanced manner. Before delving further, we introduce a family of paths that is essential for developing and analyzing path-conservative numerical schemes. The DLM theory, established by Dal Maso, LeFloch, and Murat in \cite{dal1995definition}, represents nonconservative products as Borel measures, with their definition depending on the selection of a family of paths connecting two arbitrary states. Subsequently, Par{\'e}s introduced the concept of path-conservative numerical schemes in \cite{pares2006numerical}.

Here we define a Lipschitz continuous path function in the phase space $\Omega\subset\mathbb{R}^l$ 
\begin{equation}\label{phi_fv}
  \boldsymbol{\phi}_{j+\frac{1}{2}}(\tau) = \boldsymbol{\phi}(\tau; \boldsymbol{U}_{j+\frac{1}{2}}^{-}, \boldsymbol{U}_{j+\frac{1}{2}}^{+})
\end{equation} 
for the well-balanced nodal FV-WENO scheme, which satisfies
\begin{equation*}
  \boldsymbol{\phi}: [0,1]\times \Omega \times \Omega \mapsto \Omega, \ 
  \boldsymbol{\phi}(0;\boldsymbol{U}_{j+\frac{1}{2}}^{-},\boldsymbol{U}_{j+\frac{1}{2}}^{+}) = \boldsymbol{U}_{j+\frac{1}{2}}^{-}, \ 
  \boldsymbol{\phi}(1;\boldsymbol{U}_{j+\frac{1}{2}}^{-},\boldsymbol{U}_{j+\frac{1}{2}}^{+}) = \boldsymbol{U}_{j+\frac{1}{2}}^{+}.
\end{equation*}

For many interesting systems of balance laws, the steady-state solutions of \eqref{quasi_ncs} can generally be expressed as
\begin{equation}\label{re_x_fv}
  \boldsymbol{K}(\boldsymbol{U})_x = \boldsymbol{F}(\boldsymbol{U})_x + \mathcal{B}(\boldsymbol{U})\boldsymbol{U}_x -\boldsymbol{S}(\boldsymbol{U})= \mathcal{L}(\boldsymbol {U})\widetilde{\boldsymbol{v}}(\boldsymbol {U})_x,
\end{equation}
where $\mathcal{L}(\boldsymbol {U})$ is a matrix, and $\widetilde{\boldsymbol{v}} = (\boldsymbol{v},0)^T$ is the vector of equilibrium variables. 

We compute the local integral term in each cell. Designating $\hat x$ as $x_{j-\frac{1}{2}}$, we have $\boldsymbol{R}_{j-\frac{1}{2}}^{j,-}=0$. 
Utilizing the fact that 
\begin{equation*}
  \lim_{\delta\to0}\int_{x_0-\delta}^{x_0+\delta}\boldsymbol{S}(\boldsymbol{U}(x))\ dx=0, 
\end{equation*}
and the definition of the path function \eqref{phi_fv}, we can obtain 
\begin{equation}\label{re_int}
  \boldsymbol{R}_{j-\frac{1}{2}}^{j,+} = \boldsymbol{R}_{j-\frac{1}{2}}^{j,-} + 
  \int_0^1  -\mathcal{B}(\boldsymbol{\phi}_{j-\frac{1}{2}}(\tau))  \frac{\partial }{\partial \tau}\boldsymbol{\phi}_{j-\frac{1}{2}}(\tau) \ d\tau.
\end{equation}
At the cell interface $x_{j-\frac{1}{2}}$, the following relation is satisfied,
  \begin{equation}\label{re_x_fv_int}
   \frac{\partial }{\partial \tau} \boldsymbol{F}(\boldsymbol{\phi}_{j-\frac{1}{2}}(\tau)) + \mathcal{B}(\boldsymbol{\phi}_{j-\frac{1}{2}}(\tau))  \frac{\partial }{\partial \tau}\boldsymbol{\phi}_{j-\frac{1}{2}}(\tau) = \mathcal{L}(\boldsymbol{\phi}_{j-\frac{1}{2}}(\tau)) \frac{\partial }{\partial \tau} \widetilde{\boldsymbol{v}}(\boldsymbol{\phi}_{j-\frac{1}{2}}(\tau)).
  \end{equation}
Plugging the relation \eqref{re_x_fv_int} back into the integral \eqref{re_int}, we have the evaluation,
\begin{align}
  \boldsymbol{R}_{j-\frac{1}{2}}^{j,+} \overset{\eqref{re_x_fv_int}}{=} & \boldsymbol{R}_{j-\frac{1}{2}}^{j,-} +\int_0^1 \frac{\partial }{\partial \tau} \boldsymbol{F}(\boldsymbol{\phi}_{j-\frac{1}{2}}(\tau)) \ d\tau - \int_0^1 \mathcal{L}(\boldsymbol{\phi}_{j-\frac{1}{2}}(\tau)) \frac{\partial }{\partial \tau} \widetilde{\boldsymbol{v}}(\boldsymbol{\phi}_{j-\frac{1}{2}}(\tau)) \ d\tau  \notag \\
  =& \boldsymbol{R}_{j-\frac{1}{2}}^{j,-} + 
  \boldsymbol{F}(\boldsymbol{U}_{j-\frac{1}{2}}^+) - \boldsymbol{F}(\boldsymbol{U}_{j-\frac{1}{2}}^-) - \int_0^1 \mathcal{L}(\boldsymbol{\phi}_{j-\frac{1}{2}}(\tau)) \frac{\partial }{\partial \tau} \widetilde{\boldsymbol{v}}(\boldsymbol{\phi}_{j-\frac{1}{2}}(\tau)) \ d\tau. \label{Rplusm_fv}
\end{align}
Similarly, $\boldsymbol{R}_{j+\frac{1}{2}}^{j,\pm}$ proceeds as follows:
\begin{align}
  \boldsymbol{R}_{j+\frac{1}{2}}^{j,-} =& \boldsymbol{R}_{j-\frac{1}{2}}^{j,+} + \int_{I_j} \boldsymbol{S}(\boldsymbol{U})-\mathcal{B}(\boldsymbol{U}) \boldsymbol{U}_x\ dx \notag\\
  \overset{\eqref{re_x_fv}}{=}& \boldsymbol{R}_{j-\frac{1}{2}}^{j,+} + \int _{I_j} \boldsymbol{F}(\boldsymbol{U})_x - \mathcal{L}(\boldsymbol {U})\widetilde{\boldsymbol{v}}(\boldsymbol {U})_x \ dx \notag\\
  =& \boldsymbol{R}_{j-\frac{1}{2}}^{j,+} + \boldsymbol{F}(\boldsymbol{U}_{j+\frac{1}{2}}^{-}) - \boldsymbol{F}(\boldsymbol{U}_{j-\frac{1}{2}}^{+})  - \int _{I_j}\mathcal{L}(\boldsymbol {U})\widetilde{\boldsymbol{v}}(\boldsymbol {U})_x \ dx, \label{Rminus_fv}\\
  \boldsymbol{R}_{j+\frac{1}{2}}^{j,+} =& \boldsymbol{R}_{j+\frac{1}{2}}^{j,-} + 
  \int_0^1  -\mathcal{B}(\boldsymbol{\phi}_{j+\frac{1}{2}}(\tau))  \frac{\partial }{\partial \tau}\boldsymbol{\phi}_{j+\frac{1}{2}}(\tau) \ d\tau \notag \\
  \overset{\eqref{re_x_fv_int}}{=} & \boldsymbol{R}_{j+\frac{1}{2}}^{j,-} +\int_0^1 \frac{\partial }{\partial \tau} \boldsymbol{F}(\boldsymbol{\phi}_{j+\frac{1}{2}}(\tau)) \ d\tau - \int_0^1 \mathcal{L}(\boldsymbol{\phi}_{j+\frac{1}{2}}(\tau)) \frac{\partial }{\partial \tau} \widetilde{\boldsymbol{v}}(\boldsymbol{\phi}_{j+\frac{1}{2}}(\tau)) \ d\tau   \notag \\
  =& \boldsymbol{R}_{j+\frac{1}{2}}^{j,-} + 
  \boldsymbol{F}(\boldsymbol{U}_{j+\frac{1}{2}}^+) - \boldsymbol{F}(\boldsymbol{U}_{j+\frac{1}{2}}^-) - \int_0^1 \mathcal{L}(\boldsymbol{\phi}_{j+\frac{1}{2}}(\tau)) \frac{\partial }{\partial \tau} \widetilde{\boldsymbol{v}}(\boldsymbol{\phi}_{j+\frac{1}{2}}(\tau)) \ d\tau. \label{Rplus_fv}
\end{align}
Considering the well-balanced property, several special treatments are required for calculating the integral term in  \eqref{Rplusm_fv}-\eqref{Rplus_fv}. We first define a linear segment path to connect the left and right states of the equilibrium functions $\boldsymbol{V} = (\boldsymbol{v},w)^T$ in the phase space,  given by 
\begin{equation}\label{linear_path_V_fv}
  \boldsymbol{\phi}_{j-\frac{1}{2}}(\tau) = {\boldsymbol{V}}_{j-\frac{1}{2}}^- + \tau({\boldsymbol{V}}_{j-\frac{1}{2}}^+ - {\boldsymbol{V}}_{j-\frac{1}{2}}^-).
\end{equation}
In this case, the integral term in \eqref{Rplusm_fv} reduces to
  \begin{equation}\label{v_tilde_fv}
    \int_0^1 \mathcal{L}(\boldsymbol{\phi}_{j-\frac{1}{2}}(\tau)) \left(\widetilde{\boldsymbol{v}}_{j-\frac{1}{2}}^{+} - \widetilde{\boldsymbol{v}}_{j-\frac{1}{2}}^{-}\right)\ d\tau,
  \end{equation}
which is evaluated numerically using a quadrature formula that preserves the global order of accuracy. The choice of this path connecting equilibrium variables, as given by \eqref{linear_path_V_fv}, is crucial for achieving a well-balanced treatment.

Next, we demonstrate how to compute the integral term in the last line of \eqref{Rminus_fv} while maintaining high-order accuracy, taking fifth-order accuracy as an example. We obtain the reconstruction values at four Gauss-Lobatto points in each cell $I_j$ by reconstructing the nodal values of the equilibrium variables $\boldsymbol{v}_j =  \mathbb{V}(\boldsymbol{u}_j,w_j)$. The integral is then calculated by interpolating a fifth-order polynomial with these four reconstruction values. Since the steady states are maintained exactly at all midpoints, the obtained reconstruction values remain constant in equilibrium, causing both integrals to vanish.

\begin{remark}
  The discontinuous Galerkin method with flux globalization, as introduced in \cite{zhang2025well}, can also be treated under the flux localization framework. It produces nearly identical performance.
\end{remark}

\begin{remark}
  The integral term defined in \eqref{glo_flux} is evaluated locally within each cell, with the reference value $\hat x$ reset at the left cell interface. While the spatially integral quantity $\mathcal{Q}$ appearing in the equilibrium variable \eqref{rela} is still computed by a global recursive quadrature. Therefore, the claimed locality concerns the flux-localization step only, and does not mean that the complete construction is local when the target equilibrium itself contains a global integral.
\end{remark}

\subsection{Well-balanced property}
The following theorem establishes the well-balanced property of the nodal FV-WENO scheme.
\begin{theorem}\label{wb_semi_fv}
  If the discrete data is at a steady state, which means if ${\boldsymbol{v}}_{j} = \text{constant}:= c$ at all the nodal values, then 
  the semi-discrete nodal FV-WENO scheme \eqref{scheme} is exact for the ``discrete" one-dimensional steady states at all the nodal values.
\end{theorem}

\begin{proof}
  
 
Assume that $\boldsymbol{v}_{j} = \mathbb{V}(\boldsymbol{u}_{j},w_{j})$  are the equilibrium variables at the midpoints  $x_{j}$, where $j$ ranges from 1 to $N$, and satisfy ${\boldsymbol{v}}_{j} = \text{constant}:= c$. Since the WENO reconstruction \ref{a1} preserves constant data, the reconstructed values for the equilibrium variables are also constant. 

Based on the definition of the updated boundary values in \eqref{u_modify_fv}, we conclude that $\boldsymbol{U}_{j+\frac{1}{2}}^{*,-}$ is equal to $\boldsymbol{U}_{j+\frac{1}{2}}^{*,+}$. As a result, the modified numerical flux, as specified in \eqref{LF_mod_fv}, simplifies to
  \begin{equation*}
    \widehat{\boldsymbol{K}}^{j,\text{mod}}_{j\pm\frac{1}{2}} = \dfrac{1}{2}\left(\boldsymbol{K}^j(\boldsymbol{U}_{j\pm\frac{1}{2}}^-) +\boldsymbol{K}^j(\boldsymbol{U}_{j\pm\frac{1}{2}}^+) \right).
  \end{equation*}
  
  Since ${\widetilde{\boldsymbol{v}}} = (\boldsymbol{v},0)^T$, it follows that ${\widetilde{\boldsymbol{v}}}_{j-\frac{1}{2}}^- = \widetilde{\boldsymbol{v}}_{j-\frac{1}{2}}^+$. Substituting this into \eqref{v_tilde_fv}, equation \eqref{Rplusm_fv} can be expressed as
  \begin{equation}\label{k1}
    \boldsymbol{R}_{j-\frac{1}{2}}^{j,+}
    = \boldsymbol{R}_{j-\frac{1}{2}}^{j,-} + \boldsymbol{F}(\boldsymbol{U}_{j-\frac{1}{2}}^{+}) - \boldsymbol{F}(\boldsymbol{U}_{j-\frac{1}{2}}^{-})\Longrightarrow \boldsymbol{K}^j(\boldsymbol{U}_{j-\frac{1}{2}}^{+})=\boldsymbol{K}^j(\boldsymbol{U}_{j-\frac{1}{2}}^{-}). 
  \end{equation} 
  
 Similarly, following the same derivation process for \eqref{Rplus_fv}, we obtain $\boldsymbol{K}^j(\boldsymbol{U}_{j+\frac{1}{2}}^{+})=\boldsymbol{K}^j(\boldsymbol{U}_{j+\frac{1}{2}}^{-})$.
  
Furthermore, due to the constant reconstructed values in each cell, the interpolation polynomials are also constant. This implies that the integral in equation \eqref{Rminus_fv} also evaluates to zero. Consequently, equation \eqref{Rminus_fv} simplifies to
  \begin{equation}\label{k2}
    \boldsymbol{R}_{j+\frac{1}{2}}^{j,-} 
    = \boldsymbol{R}_{j-\frac{1}{2}}^{j,+} + \boldsymbol{F}(\boldsymbol{U}_{j+\frac{1}{2}}^{-}) - \boldsymbol{F}(\boldsymbol{U}_{j-\frac{1}{2}}^{+})\Longrightarrow \boldsymbol{K}^j(\boldsymbol{U}_{j+\frac{1}{2}}^{-})=\boldsymbol{K}^j(\boldsymbol{U}_{j-\frac{1}{2}}^{+}). 
  \end{equation}
Substituting \eqref{k1} and \eqref{k2} back into \eqref{scheme}, we obtain
  \begin{align*}
    \text{RHS}_j(\boldsymbol{U})  = 
    -\left(\widehat{\boldsymbol{K}}^{j,\text{mod}}_{j+\frac{1}{2}}   - \widehat{\boldsymbol{K}}^{j,\text{mod}}_{j-\frac{1}{2}}\right)
    =  - \left(\boldsymbol{K}^j(\boldsymbol{U}_{j+\frac{1}{2}}^-)  - \boldsymbol{K}^j(\boldsymbol{U}_{j-\frac{1}{2}}^+) \right) = 0.
  \end{align*}
  This completes the proof of the well-balanced property.
\end{proof}

The fully-discrete scheme retains the same property, which can be demonstrated through a straightforward inductive argument. Thus, the detailed proof is omitted here.

\begin{theorem}\label{wb_fully_fv}
The fully-discrete nodal FV-WENO scheme, which combines the spatial discretization in \eqref{scheme} with the third-order SSP-RK time discretization in equation \eqref{3rk_no}, preserves exactness for the ``discrete" one-dimensional steady state at all the nodal values.
\end{theorem}

\section{Extension to two-dimensional case}\label{se:e2d}
Although the primary focus of this paper is on one spatial dimension, the approach can be extended to higher dimensions. In this section, we discuss the generalization of the well-balanced nodal FV-WENO method to two dimensions.

One important distinction is that, unlike the one-dimensional case, the steady state in two dimensions cannot be described by simple algebraic equations. Instead, it requires solving partial differential equations and accommodates a variety of steady states. To date, there has been limited research on genuinely two-dimensional moving equilibrium. In this paper, we focus on quasi one-dimensional steady states and use the single-layer shallow water equations as an example to present the detailed implementation. The equations in two dimensions are given by
\begin{equation}\label{2Dsha}
  \left\{\begin{array}{l}
  h_t+(hu)_{x}+(hv)_{y}=0, \\
  (hu)_t+(hu^2+\frac{1}{2}gh^2)_x+(huv)_y = -ghb_x,\\
  (hv)_t+(huv)_x+(hv^2+\frac{1}{2}gh^2)_y = -ghb_y,
  \end{array}\right.
\end{equation}
where $h,u,v,g,b$ represent the water height, two components of the horizontal velocity, the gravitational constant, and the bottom topography function, respectively.

\subsection{Steady state solution}
To describe the nodal FV-WENO scheme with flux localization, we rewrite \eqref{2Dsha} in quasi-conservative form as
\begin{equation}\label{reswe2d}
  \boldsymbol{U}_t + \boldsymbol{K}(\boldsymbol{U})_x + \boldsymbol{H}(\boldsymbol{U})_y = 0,
\end{equation}
with
  \begin{align*}
    &\boldsymbol{K}(\boldsymbol{U}) = \boldsymbol{F}(\boldsymbol{U}) - \boldsymbol{R}^x(\boldsymbol{U}),\ \boldsymbol{R}^x(\boldsymbol{U}) = \int_{\hat x}^x \boldsymbol{S}^x(\boldsymbol{U}(\xi,y,t)) \ d\xi, \\
    &\boldsymbol{H}(\boldsymbol{U}) = \boldsymbol{G}(\boldsymbol{U}) - \boldsymbol{R}^y(\boldsymbol{U}), \ \boldsymbol{R}^y(\boldsymbol{U}) = \int_{\hat y}^y \boldsymbol{S}^y(\boldsymbol{U}(x,\eta,t)) \ d\eta,
  \end{align*}
where $\hat x,\hat y$ are arbitrary numbers, and 
\begin{align*}
  & \boldsymbol{U} = (h,hu,hv,b)^T,\\
  & \boldsymbol{F}(\boldsymbol{U}) = (hu,hu^2+\frac{1}{2}gh^2,huv,0)^T,\ \boldsymbol{S}^x(\boldsymbol{U}) =(0,-ghb_x,0,0)^T, \\
  & \boldsymbol{G}(\boldsymbol{U}) = (hv,huv,hv^2+\frac{1}{2}gh^2,0)^T,\ \boldsymbol{S}^y(\boldsymbol{U}) =(0,0,-ghb_y,0)^T.
\end{align*}
The steady state of the system \eqref{reswe2d} satisfies 
\begin{equation}\label{ss_2d}
  \boldsymbol{K}(\boldsymbol{U})_x + \boldsymbol{H}(\boldsymbol{U})_y = 0.
\end{equation}
Analogous to one-dimensional case, \eqref{ss_2d} can be equivalently written as 
\begin{equation}\label{sse_2d}
  \mathcal{L}^x(\boldsymbol {U})\widetilde{\boldsymbol{v}^x}(\boldsymbol {U})_x+\mathcal{L}^y(\boldsymbol {U})\widetilde{\boldsymbol{v}^y}(\boldsymbol {U})_y = 0,
\end{equation}
with
\begin{equation*}
  \mathcal{L}^x = \left[\begin{array}{ccccc}
    0 & 1 & 0 & 0 \\
    h & u & 0 & 0 \\ 
    0 & v & hu & 0 \\
    0 & 0 & 0 & 1 \\
  \end{array}\right], \
  \mathcal{L}^y = \left[\begin{array}{ccccc}
    0 & 0 & 1 & 0 \\
    0 & hv & u & 0 \\ 
    h & 0 & v & 0 \\
    0 & 0 & 0 & 1 \\
  \end{array}\right], \
\end{equation*}
and the equilibrium variable vectors in the $x$- and $y$-direction are defined as 
\begin{equation*}
  \widetilde{\boldsymbol{v}^x} = (\boldsymbol{v}^x,0)^T = (\mathcal{E},hu,v,0)^T,\quad \widetilde{\boldsymbol{v}^y} = (\boldsymbol{v}^y,0)^T = (E,u,hv,0)^T.
\end{equation*} 
Here, $\mathcal{E} = \frac{1}{2}u^2+g(h+b)$ and $E = \frac{1}{2}v^2+g(h+b)$.

Although \eqref{sse_2d} can be reformulated into a more compact and simpler form
\begin{equation}\label{ssg_swe}
  (hu)_x+(hv)_y = 0, \quad \mathcal{E}_x + vu_y = 0,\quad uv_x + E_y = 0, 
\end{equation}
the explicit solution of this system is difficult to obtain, making it more challenging to design a well-balanced scheme capable of preserving genuinely two-dimensional equilibrium. However, simpler requirements can be applied to achieve the following steady state
\begin{equation}\label{ss_quasi}
  \begin{aligned}
  & \boldsymbol{K}(\boldsymbol{U})_x =\boldsymbol{F}(\boldsymbol{U})_x - \boldsymbol{S}^x(\boldsymbol{U})= \mathcal{L}^x(\boldsymbol {U})\widetilde{\boldsymbol{v}^x}(\boldsymbol {U})_x=0, \\
  & \boldsymbol{H}(\boldsymbol{U})_y =\boldsymbol{G}(\boldsymbol{U})_y - \boldsymbol{S}^y(\boldsymbol{U})= \mathcal{L}^y(\boldsymbol {U})\widetilde{\boldsymbol{v}^y}(\boldsymbol {U})_y=0.
  \end{aligned} 
\end{equation}
Specifically, the steady-state conditions are given by
\begin{equation}\label{ss_px}
  (hu)_x=0,\quad v=0,\quad \mathcal{E}_x=0,
\end{equation}
and
\begin{equation}\label{ss_py}
  u=0,\quad  (hv)_y=0,\quad E_y=0.
\end{equation}
The well-balanced nodal finite volume scheme presented in this section is designed to precisely maintain the quasi one-dimensional steady state defined by \eqref{ss_px} and \eqref{ss_py}.

\subsection{The well-balanced scheme with flux localization}\label{subse:wb2d}
We consider rectangular mesh cells, with the computational domain partitioned into a set of uniform cells:
$$\mathcal {D} = \bigcup_{i=1}^{Nx}\bigcup_{j=1}^{Ny} I_{ij} = \bigcup_{i=1}^{Nx}\bigcup_{j=1}^{Ny}[x_{i-\frac{1}{2}}, x_{i+\frac{1}{2}}]\times[y_{j-\frac{1}{2}}, y_{j+\frac{1}{2}}],$$
where $\Delta x = x_{i+\frac{1}{2}}-x_{i-\frac{1}{2}}$ and $\Delta y = y_{j+\frac{1}{2}}-y_{j-\frac{1}{2}}$. The nodal values $\boldsymbol{U}_{ij}$ approximate the pointwise values of $\boldsymbol{U}(x_i,y_j,t)$. The well-balanced finite volume scheme is given by
\begin{equation}\label{scheme_2d}
  \dfrac{d}{dt}\int_{I_{ij}}\boldsymbol{U} \ dx dy:= \text{RHS}_{ij}(\boldsymbol{U}) = - \left(\widehat{\boldsymbol{K}}^{i,\text{mod}}_{i+\frac{1}{2},j}   - \widehat{\boldsymbol{K}}^{i,\text{mod}}_{i-\frac{1}{2},j}\right) - \left(\widehat{\boldsymbol{H}}^{j,\text{mod}}_{i,j+\frac{1}{2}}   - \widehat{\boldsymbol{H}}^{j,\text{mod}}_{i,j-\frac{1}{2}}\right).
\end{equation}
The numerical fluxes in $x$- and $y$-direction are defined as
  \begin{align*}
  &\widehat{\boldsymbol{K}}^{i,\text{mod}}_{i\pm\frac{1}{2},j} 
  = \sum_{\beta} \omega_{\beta} \widehat{\boldsymbol{K}}^{i,\text{mod}}_{i\pm\frac{1}{2},j_{\beta}}, \quad 
  \widehat{\boldsymbol{H}}^{j,\text{mod}}_{i,j\pm\frac{1}{2}} = \sum_{\alpha} \omega_{\alpha}\widehat{\boldsymbol{H}}^{j,\text{mod}}_{i_{\alpha},j\pm\frac{1}{2}},\\
  &\widehat{\boldsymbol{K}}^{i,\text{mod}}_{i\pm\frac{1}{2},j_{\beta}} = \frac{1}{2}\left( \boldsymbol{K}^{i}(\boldsymbol{U}^{-}_{i\pm\frac{1}{2},j_{\beta}}) + \boldsymbol{K}^{i}(\boldsymbol{U}^{+}_{i\pm\frac{1}{2},j_{\beta}}) - \alpha_{i+\frac{1}{2},j_{\beta}} (\boldsymbol{U}^{*,+}_{i\pm\frac{1}{2},j_{\beta}} -\boldsymbol{U}^{*,-}_{i\pm\frac{1}{2},j_{\beta}} )\right),\\
  &\widehat{\boldsymbol{H}}^{j,\text{mod}}_{i_{\alpha},j\pm\frac{1}{2}} = \frac{1}{2}\left( \boldsymbol{H}^{j}(\boldsymbol{U}^{-}_{i_{\alpha},j\pm\frac{1}{2}}) + \boldsymbol{H}^{j}(\boldsymbol{U}^{+}_{i_{\alpha},j\pm\frac{1}{2}}) - \alpha_{i_{\alpha},j\pm\frac{1}{2}} (\boldsymbol{U}^{*,+}_{i_{\alpha},j\pm\frac{1}{2}} -\boldsymbol{U}^{*,-}_{i_{\alpha},j\pm\frac{1}{2}} )\right),
  \end{align*}
where $\omega_\beta,y_{j_{\beta}}$  are Gaussian quadrature weights and nodes on $[y_{j-\frac{1}{2}},y_{j+\frac{1}{2}}]$, and $\omega_\alpha,x_{i_{\alpha}}$ defined on $[x_{i-\frac{1}{2}},x_{i+\frac{1}{2}}]$, respectively.

As in the one-dimensional case, the cell interface values and their modifications are defined by
  \begin{align*}
  &\boldsymbol{U}^{\pm}_{i+\frac{1}{2},j_{\beta}} = \left(\mathbb{U}(\boldsymbol{v}_{i+\frac{1}{2},j_{\beta}}^{x,\pm},b_{i+\frac{1}{2},j_{\beta}}^{\pm}),b_{i+\frac{1}{2},j_{\beta}}^{\pm}\right)^T, \quad
  \boldsymbol{U}_{i+\frac{1}{2},j_{\beta}}^{*,\pm} 
  = \left(\mathbb{U}(\boldsymbol{v}_{i+\frac{1}{2},j_{\beta}}^{x,\pm},b_{i+\frac{1}{2},j_{\beta}}^{*}),b_{i+\frac{1}{2},j_{\beta}}^{*}\right)^T,\\
  &\boldsymbol{U}^{\pm}_{i_{\alpha},j\pm\frac{1}{2}} = \left(\mathbb{U}(\boldsymbol{v}_{i_{\alpha},j\pm\frac{1}{2}}^{y,\pm},b_{i_{\alpha},j\pm\frac{1}{2}}^{\pm}),b_{i_{\alpha},j\pm\frac{1}{2}}^{\pm}\right)^T, \quad
  \boldsymbol{U}_{i_{\alpha},j\pm\frac{1}{2}}^{*,\pm} 
  = \left(\mathbb{U}(\boldsymbol{v}_{i_{\alpha},j\pm\frac{1}{2}}^{y,\pm},b_{i_{\alpha},j\pm\frac{1}{2}}^{*}),b_{i_{\alpha},j\pm\frac{1}{2}}^{*}\right)^T,
  \end{align*}
with 
\begin{equation*}
  b_{i+\frac{1}{2},j_{\beta}}^{*} = \min(b_{i+\frac{1}{2},j_{\beta}}^{-},b_{i+\frac{1}{2},j_{\beta}}^{+}),\quad 
  b_{i_{\alpha},j\pm\frac{1}{2}}^{*} = \min(b_{i_{\alpha},j\pm\frac{1}{2}}^{-},b_{i_{\alpha},j\pm\frac{1}{2}}^{+}).
\end{equation*}
The high-order approximation values $\boldsymbol{v}_{i+\frac{1}{2},j_{\beta}}^{x,\pm}$ and $\boldsymbol{v}_{i_{\alpha},j\pm\frac{1}{2}}^{y,\pm}$  
are obtained from the nodal values $\{\boldsymbol{v}^x_{kl},\boldsymbol{v}^y_{kl}\}_{k,l\in S_{ij}}$ by WENO reconstruction procedure, where $S_{ij}$ denotes the index set of cells used in the reconstruction. 
The localized fluxes are constituted by
\begin{equation*}
  \boldsymbol{K}^{i}(\boldsymbol{U}^{\pm}_{i\pm\frac{1}{2},j_{\beta}}) = \boldsymbol{F}(\boldsymbol{U}^{\pm}_{i\pm\frac{1}{2},j_{\beta}}) - \boldsymbol{R}^{i,\pm}_{i\pm\frac{1}{2},j_{\beta}}, \quad 
  \boldsymbol{H}^{j}(\boldsymbol{U}^{\pm}_{i_{\alpha},j\pm\frac{1}{2}}) = \boldsymbol{G}(\boldsymbol{U}^{\pm}_{i_{\alpha},j\pm\frac{1}{2}}) - \boldsymbol{R}^{j,\pm}_{i_{\alpha},j\pm\frac{1}{2}},
\end{equation*}
Using the path functions $\boldsymbol{\phi}_{i-\frac{1}{2},j_{\beta}}(\tau)$ in $x$-direction and $\boldsymbol{\phi}_{i_{\alpha},j-\frac{1}{2}}(s)$ in $y$-direction, along with the relationship \eqref{ss_quasi}, the local integral terms are evaluated analogously to the one-dimensional case \eqref{Rplusm_fv}-\eqref{Rplus_fv},
  \begin{align*}
    &\boldsymbol{R}^{i,-}_{i-\frac{1}{2},j_{\beta}} = 0,\\
    &\boldsymbol{R}^{i,+}_{i-\frac{1}{2},j_{\beta}} = \boldsymbol{R}^{i,-}_{i-\frac{1}{2},j_{\beta}} + 
    \boldsymbol{F}(\boldsymbol{U}_{i-\frac{1}{2},j_{\beta}}^+) - \boldsymbol{F}(\boldsymbol{U}_{i-\frac{1}{2},j_{\beta}}^-) - \int_0^1 \mathcal{L}^x(\boldsymbol{\phi}_{i-\frac{1}{2},j_{\beta}}(\tau)) \frac{\partial }{\partial \tau} \widetilde{\boldsymbol{v}^x}(\boldsymbol{\phi}_{i-\frac{1}{2},j_{\beta}}(\tau)) \ d\tau, \\
    &\boldsymbol{R}^{i,-}_{i+\frac{1}{2},j_{\beta}} = \boldsymbol{R}^{i,+}_{i-\frac{1}{2},j_{\beta}}  + \boldsymbol{F}(\boldsymbol{U}_{i+\frac{1}{2},j_{\beta}}^{-}) - \boldsymbol{F}(\boldsymbol{U}_{i-\frac{1}{2},j_{\beta}}^{+})  - \int _{x_{i-\frac{1}{2}}}^{x_{i+\frac{1}{2}}}\mathcal{L}^x(\boldsymbol {U})\widetilde{\boldsymbol{v}^x}(\boldsymbol {U})_x \ dx, \\
    &\boldsymbol{R}^{i,+}_{i+\frac{1}{2},j_{\beta}} = \boldsymbol{R}^{i,-}_{i+\frac{1}{2},j_{\beta}} + 
    \boldsymbol{F}(\boldsymbol{U}_{i+\frac{1}{2},j_{\beta}}^+) - \boldsymbol{F}(\boldsymbol{U}_{i+\frac{1}{2},j_{\beta}}^-) - \int_0^1 \mathcal{L}^x(\boldsymbol{\phi}_{i+\frac{1}{2},j_{\beta}}(\tau)) \frac{\partial }{\partial \tau} \widetilde{\boldsymbol{v}^x}(\boldsymbol{\phi}_{i+\frac{1}{2},j_{\beta}}(\tau)) \ d\tau, \\
  \end{align*}
and 
  \begin{align*}
    &\boldsymbol{R}^{j,-}_{i_{\alpha},j-\frac{1}{2}} = 0,\\
    &\boldsymbol{R}^{j,+}_{i_{\alpha},j-\frac{1}{2}} = \boldsymbol{R}^{j,-}_{i_{\alpha},j-\frac{1}{2}} + 
    \boldsymbol{G}(\boldsymbol{U}_{i_{\alpha},j-\frac{1}{2}}^+) - \boldsymbol{G}(\boldsymbol{U}_{i_{\alpha},j-\frac{1}{2}}^-) - \int_0^1 \mathcal{L}^y(\boldsymbol{\phi}_{i_{\alpha},j-\frac{1}{2}}(s)) \frac{\partial }{\partial s} \widetilde{\boldsymbol{v}^y}(\boldsymbol{\phi}_{i_{\alpha},j-\frac{1}{2}}(s)) \ ds, \\
    &\boldsymbol{R}^{j,-}_{i_{\alpha},j+\frac{1}{2}} = \boldsymbol{R}^{j,+}_{i_{\alpha},j-\frac{1}{2}}  + \boldsymbol{G}(\boldsymbol{U}_{i_{\alpha},j+\frac{1}{2}}^{-}) - \boldsymbol{G}(\boldsymbol{U}_{i_{\alpha},j-\frac{1}{2}}^{+})  - \int _{y_{j-\frac{1}{2}}}^{y_{j+\frac{1}{2}}}\mathcal{L}^y(\boldsymbol {U})\widetilde{\boldsymbol{v}^y}(\boldsymbol {U})_y \ dy, \\
    &\boldsymbol{R}^{j,+}_{i_{\alpha},j+\frac{1}{2}} = \boldsymbol{R}^{j,-}_{i_{\alpha},j+\frac{1}{2}} + 
    \boldsymbol{G}(\boldsymbol{U}_{i_{\alpha},j+\frac{1}{2}}^+) - \boldsymbol{G}(\boldsymbol{U}_{i_{\alpha},j+\frac{1}{2}}^-) - \int_0^1 \mathcal{L}^y(\boldsymbol{\phi}_{i_{\alpha},j+\frac{1}{2}}(s)) \frac{\partial }{\partial s} \widetilde{\boldsymbol{v}^y}(\boldsymbol{\phi}_{i_{\alpha},j+\frac{1}{2}}(s)) \ ds. 
  \end{align*}
We define the path functions as linear interpolants between the left and right equilibrium states:
  \begin{align*}
  &\boldsymbol{\phi}_{i-\frac{1}{2},j_{\beta}}(\tau) = \boldsymbol{V}_{i-\frac{1}{2},j_{\beta}}^{x,-} + \tau(\boldsymbol{V}_{i-\frac{1}{2},j_{\beta}}^{x,+} - \boldsymbol{V}_{i-\frac{1}{2},j_{\beta}}^{x,-}),\\
  &\boldsymbol{\phi}_{i_{\alpha},j-\frac{1}{2}}(s) = \boldsymbol{V}_{i_{\alpha},j-\frac{1}{2}}^{y,-} + s(\boldsymbol{V}_{i_{\alpha},j-\frac{1}{2}}^{y,+} - \boldsymbol{V}_{i_{\alpha},j-\frac{1}{2}}^{y,-}),
  \end{align*}
where $\boldsymbol{V}^x = (\boldsymbol{v}^x,b)^T$ and $\boldsymbol{V}^y = (\boldsymbol{v}^y,b)^T$  are the equilibrium variables for the shallow water equations. As in the one-dimensional case, applying these paths simplifies the evaluation of the flux correction integrals. For brevity, we omit the analogous derivations.

\begin{theorem}
  Both the semi-discrete nodal FV-WENO scheme \eqref{scheme_2d} and the corresponding fully-discrete scheme for the shallow water equations are exact at all nodal values for the two-dimensional steady states defined in \eqref{ss_px} and \eqref{ss_py}. 
\end{theorem}
\begin{proof}
  As the proof closely parallels that of Theorem \ref{wb_semi_fv}, it is omitted here for brevity.
\end{proof}

We conclude this section with several remarks regarding the proposed well-balanced nodal FV-WENO schemes.

\begin{remark}
For the fully-discrete scheme, the nodal values $\boldsymbol{U}_{ij}$ are also updated by solving the same linear equation system as in \eqref{li_s} with the mass matrix $M$ here being the block diagonal matrix.
\end{remark}


\begin{remark}
We use the single-layer shallow water equations as a representative example to illustrate the detailed implementation of our proposed nodal FV-WENO scheme in higher dimensions. The approach is readily extendable to other complex balance laws such as Euler equations with gravity and Ripa model, see Appendix \ref{b1} for more details.
\end{remark}

\section{Numerical examples}\label{se:nu}
In this section, we present several numerical tests to evaluate the performance of the well-balanced nodal FV-WENO scheme for a range of hyperbolic balance laws. As illustrative examples, we implement a fifth-order nodal FV-WENO scheme combined with a third-order SSP-RK time discretization. In all numerical experiments, the CFL number is set to 0.5. Specifically, the time step is chosen as
$\Delta t = \text{CFL}\frac{ \Delta x }{\max_{j}{\alpha_{j+\frac{1}{2}}}}$ where $\alpha_{j+\frac{1}{2}}$ is defined in \eqref{LF_alpha_fv}. 
For accuracy tests, the time step is taken sufficiently small so that spatial discretization errors dominate.

\subsection{Euler equations with gravitation}\label{subse:euler}
We first consider the compressible Euler equations with gravitation:
\begin{equation}\label{euler}
  \boldsymbol{U} = \left[\begin{array}{c}
    \rho \\
    \rho u \\
    E \\ 
    \phi
  \end{array}\right], 
  \boldsymbol{F}(\boldsymbol{U}) = \left[\begin{array}{c}
    \rho u \\
    \rho u^2 + p \\
    u(E+p) \\ 
    0
  \end{array}\right], 
  \mathcal{B}(\boldsymbol {U}) = \left[\begin{array}{ccccc}
    0 & 0 & 0 & 0 \\
    0 & 0 & 0 & \rho \\ 
    0 & 0 & 0 & \rho u \\
    0 & 0 & 0 & 0 \\
  \end{array}\right],
  \boldsymbol{S}(\boldsymbol{U}) = 0,
\end{equation}
where $\rho,u,p,E$ represent the mass density, gas velocity, pressure, and total energy. The gravitational potential $\phi$ is the time-independent gravitational and specified explicitly. For the ideal gas, the pressure is given by the equation of state
\begin{equation*}
  p = (\gamma-1)\left(E - \dfrac{1}{2}\rho u^2\right).
\end{equation*}
We adopt the well-balanced strategy described in Section~\ref{se:fv} to preserve steady adiabatic flows, which satisfy the following equilibrium conditions:
\begin{equation}\label{adi:euler}
  s = \text{constant}, \quad \rho u= \text{constant}, \quad  \varepsilon := \dfrac{1}{2}u^2 + h + \phi = \text{constant}.
\end{equation}
where $s$ denotes the specific entropy and $h$ is the specific enthalpy.

\begin{example}{\bf Well-balanced property test for the steady adiabatic flow}\label{euler:wb_adi}
\end{example}
To demonstrate the well-balanced property for the steady adiabatic flow \eqref{adi:euler}, we follow the setup in \cite{grosheintz2020well} and consider the domain $[0,2]$ with initial conditions
\begin{equation}\label{euler:adi_ini}
(\rho,u,p)(x,0) = (\rho_{eq}(x), u_{eq}(x), p_{eq}(x)),
\end{equation}
and a linear gravitational potential $\phi(x) = x$.
The equilibrium profile $(\rho_{eq}, u_{eq}, p_{eq})$ is determined by specifying the point values
\begin{equation*}
(\rho_0,u_0,p_0) = (1, -M c_{s,0},1)
\end{equation*}
at $x_0 = 0$, where $\gamma = 5/3$ is the ratio of specific heats, and $c_{s,0} = \sqrt{\gamma}$ is the sound speed at $x_0$. We examine three cases: a hydrostatic flow ($M=0$), a subsonic flow ($M=0.01$), and a supersonic flow ($M=2.5$). The values at the ghost points are set equal to the corresponding initial conditions. For the hydrostatic and subsonic cases, the simulation runs until $t = 4$, while for the supersonic case, the final time is set to $t = 1$.
Table~\ref{euler:adi_1D} reports the $L^1$ and $L^\infty$ errors using 100 uniform cells, confirming that the well-balanced property is indeed preserved.

\begin{table}[htb]
  \centering
  \caption{ Example \ref{euler:wb_adi}: $L^1$ and $L^{\infty}$ errors for the steady adiabatic flow of Euler equations.}\label{euler:adi_1D}
  \begin{tabular}{c c c c c c c c  }
    \toprule
    \multirow{2}{*} {case} &\multicolumn{3}{c}{$L^1$ error}&\multicolumn{3}{c}{$L^{\infty}$ error}\\
     \cmidrule(lr){2-4} \cmidrule(lr){5-7}
     &{$\rho$}&{$\rho u$}&{$E$}&{$\rho$}&{$\rho u$}&{$E$}\\
     \midrule
     $M=0$    & 3.60E-15 & 2.64E-15	& 4.65E-15 & 5.69E-15	& 4.93E-15 & 7.55E-15 \\
     $M=0.01$ & 1.26E-14 & 9.71E-15 &	1.93E-14 & 8.77E-15 &	7.64E-15 & 1.78E-14 \\
     $M=2.5$  & 8.50E-15 & 2.21E-14 & 4.79E-14 & 1.51E-14 & 3.55E-14 & 5.95E-14 \\ 
     \bottomrule
   \end{tabular}
\end{table}

\begin{example}{\bf Small perturbation of the steady adiabatic flow}\label{euler:per_adi}
\end{example}
Using almost the same setup as in Example~\ref{euler:wb_adi}, we introduce a localized perturbation to the initial pressure:
\begin{equation}\label{euler:adi_peradd}
(\rho,u,p)(x,0) = (\rho_{eq}(x), u_{eq}(x), p_{eq}(x)+A\exp(-100(x-\bar x)^2))
\end{equation}
where the location of the disturbance $\bar{x}$ is chosen as
\begin{equation*}
\bar x =\left\{\begin{array}{lll}
    1.0,&\text{for} \  M=0, \\
    1.1,&\text{for} \  M=0.01,\\
    1.5,&\text{for} \  M=2.5.
    \end{array}\right.
\end{equation*}
We consider a small perturbation amplitude of $A = 10^{-6}$ and set the final time to $t = 0.45$ for the hydrostatic and subsonic cases, and $t = 0.25$ for the supersonic case.

\begin{figure}[htb!]
  \centering
  \subfigure[pressure perturbation]{
  \centering
  \includegraphics[width=5.cm,scale=1]{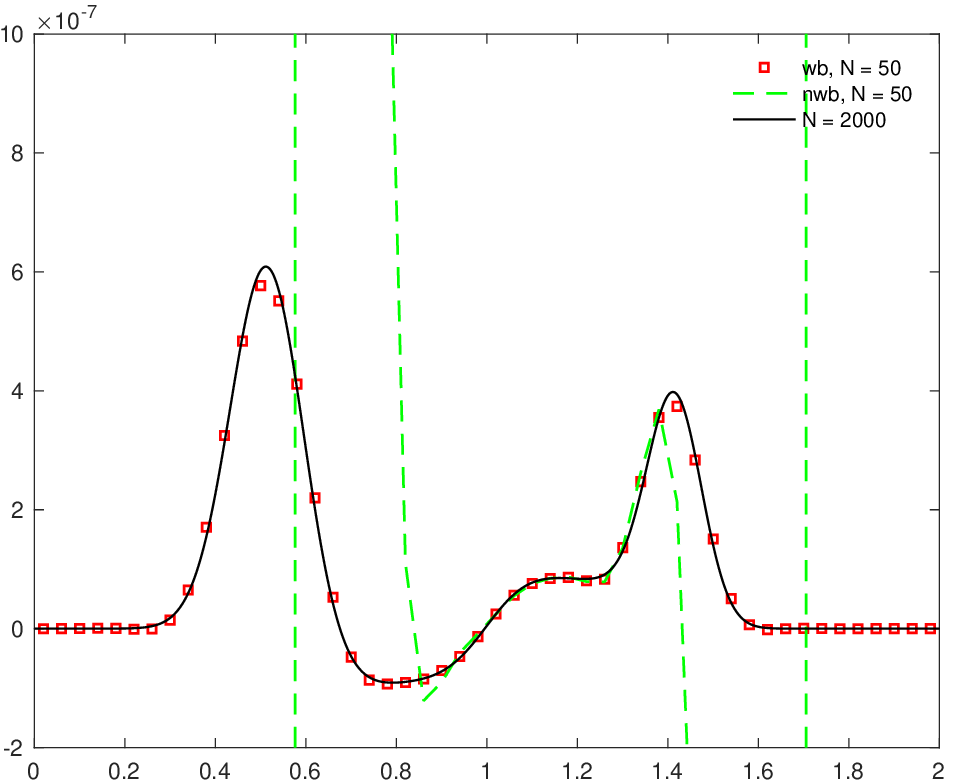}
  }
  \subfigure[velocity perturbation]{
  \centering
  \includegraphics[width=5.cm,scale=1]{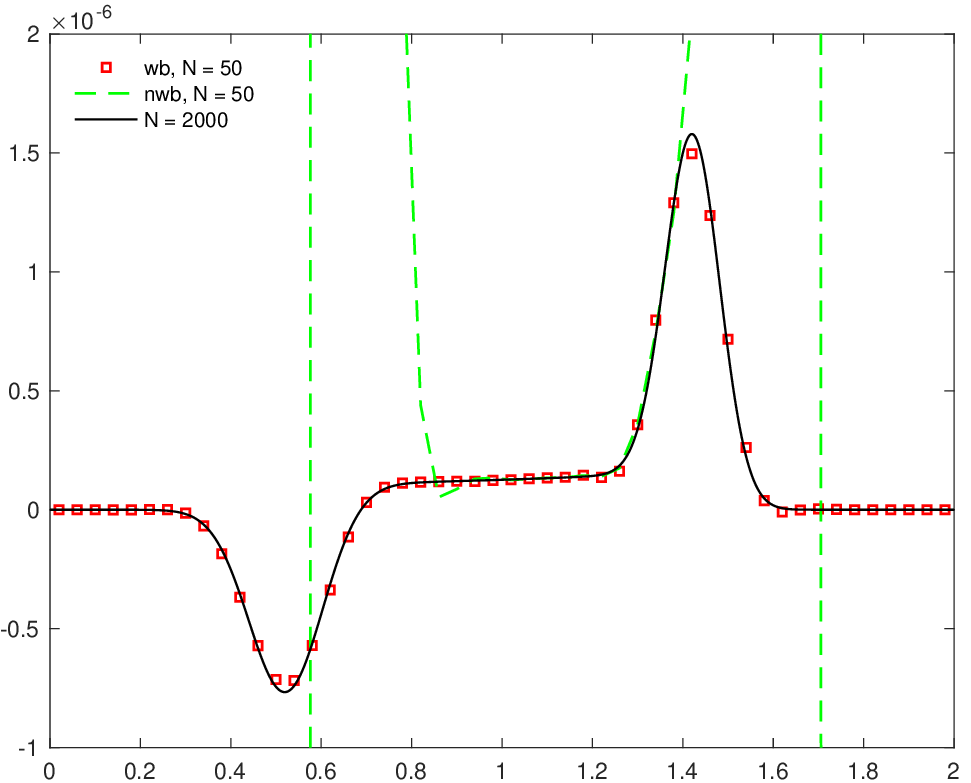}
  }
  \caption{Example \ref{euler:per_adi}: Pressure perturbation and velocity perturbation of the hydrostatic atmosphere ($M=0$) for the small amplitude wave propagation $A=10^{-6}$  in (\ref{euler:adi_peradd}). }\label{euler:persmall_hydro}
\end{figure}

\begin{figure}[htb!]
  \centering
  \subfigure[pressure perturbation]{
  \centering
  \includegraphics[width=5.cm,scale=1]{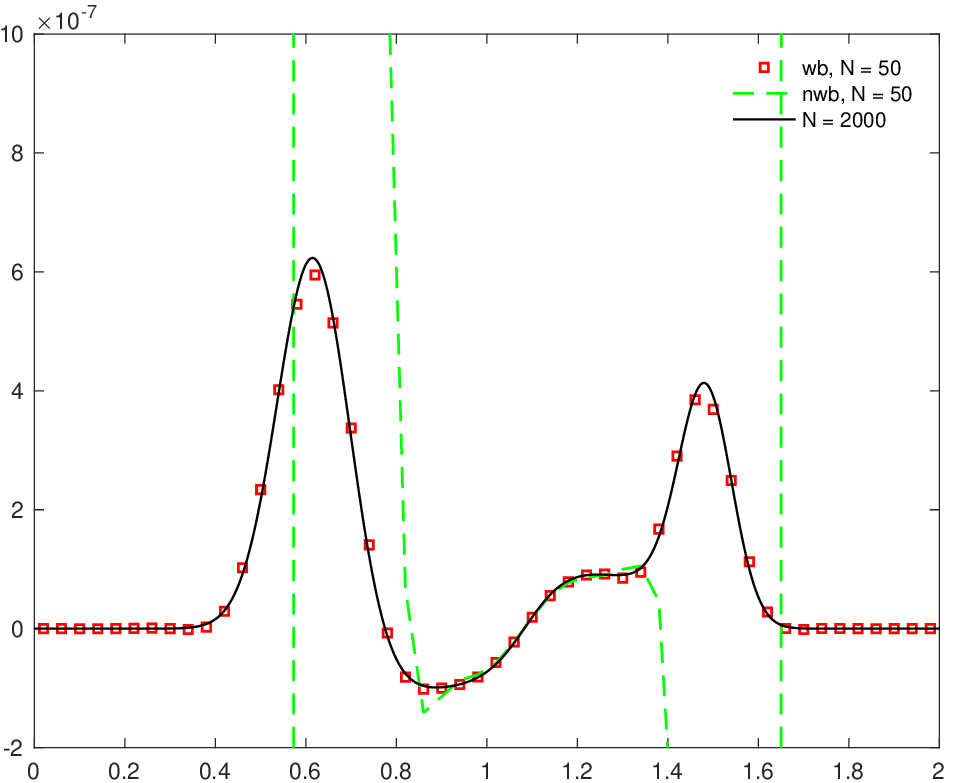}
  }
  \subfigure[velocity perturbation]{
  \centering
  \includegraphics[width=5.cm,scale=1]{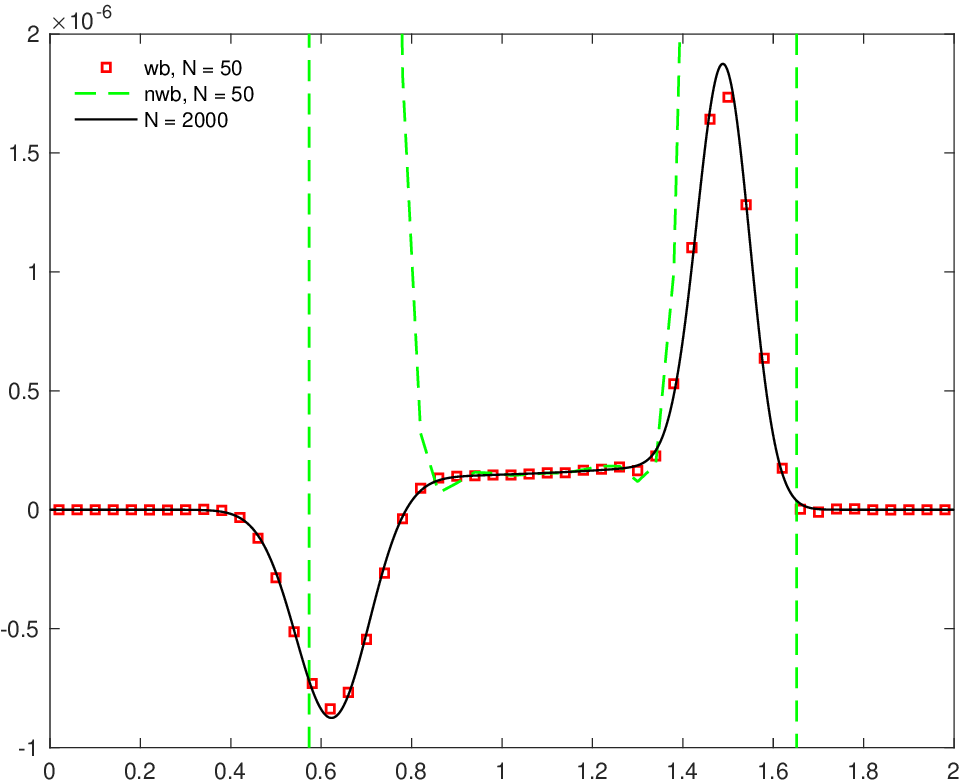}
  }
  \caption{Example \ref{euler:per_adi}: Pressure perturbation and velocity perturbation of the subsonic atmosphere ($M=0.01$) for the small amplitude wave propagation $A=10^{-6}$ in (\ref{euler:adi_peradd}). }\label{euler:persmall_sub}
\end{figure}

\begin{figure}[htb!]
  \centering
  \subfigure[pressure perturbation]{
  \centering
  \includegraphics[width=5.cm,scale=1]{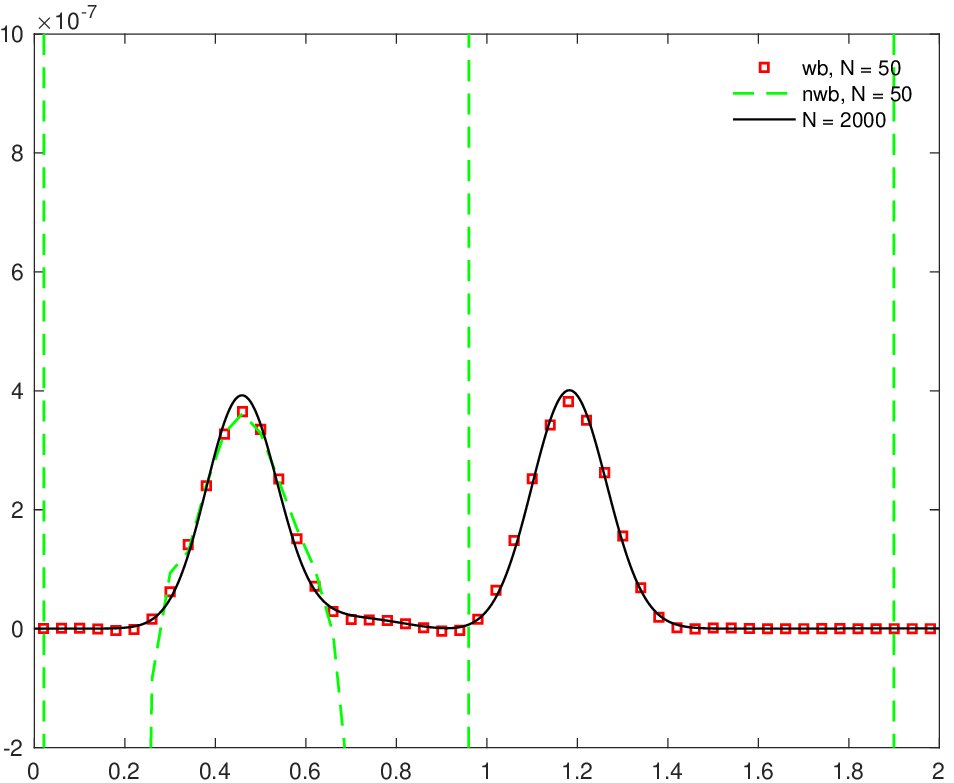}
  }
  \subfigure[velocity perturbation]{
  \centering
  \includegraphics[width=5.cm,scale=1]{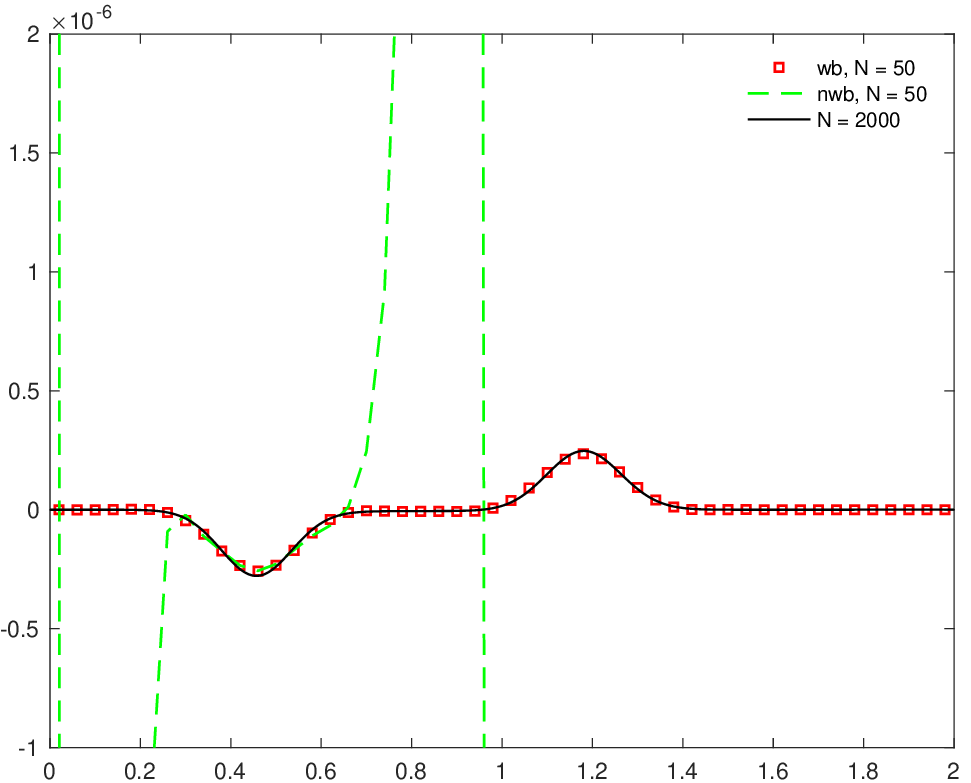}
  }
  \caption{Example \ref{euler:per_adi}: Pressure perturbation and velocity perturbation of the supersonic atmosphere ($M=2.5$) for the small amplitude wave propagation $A=10^{-6}$ in (\ref{euler:adi_peradd}). }\label{euler:persmall_sup}
\end{figure}

Numerical results computed using the proposed well-balanced nodal FV-WENO scheme, as well as a non-well-balanced classical FV-WENO scheme, are presented in Figs.~\ref{euler:persmall_hydro}, \ref{euler:persmall_sub}, and \ref{euler:persmall_sup}. These simulations are performed on a coarse mesh of 50 uniform cells for each Mach number. The reference solution is obtained using the third-order equilibrium-preserving DG method proposed in \cite{2023arXivequilibrium}, computed with a much finer mesh of 2000 cells. The results clearly demonstrate that the non-well-balanced scheme fails to capture the solution accurately, while our well-balanced method remains robust and accurate, even for non-hydrostatic equilibrium states on relatively coarse meshes.

{\begin{example}{\bf Shock tube problem}\label{euler:shock}
\end{example}
We consider the classical shock tube problem under a gravitational field. The initial conditions are given by
\begin{equation}\label{euler:sod_ini}
(\rho,u,p)(x,0) = \left\{\begin{array}{lll}
    (1,0,1),&  \text{if} \  x \leq 0.5, \\
    (0.125,0,0.1),&\text{otherwise},\\
    \end{array}\right.
\end{equation}
with a linear gravitational field $\phi(x) = x$ on the computational domain $[0,1]$. Solid walls boundary conditions are imposed and the ratio of specific heats is set as $\gamma=1.4$. We simulate the solutions up to $t=0.2$ with 200 uniform meshes, and displays the density, velocity, energy and pressure in Fig. \ref{euler:st200}. For comparison, the corresponding results computed by the well-balanced finite volume WENO scheme \cite{li2016highfve} on a refined 2000 meshes are treated as the reference solution. The gravitational force drives the density distribution toward the left, leading to the appearance of negative velocities in certain regions. It can be observed that our proposed nodal FV-WENO scheme achieves good agreement with the reference solution and remains essentially non-oscillatory.
\begin{figure}[htb!]
  \centering
  \subfigure[density]{
  \centering
  \includegraphics[width=5cm,scale=1]{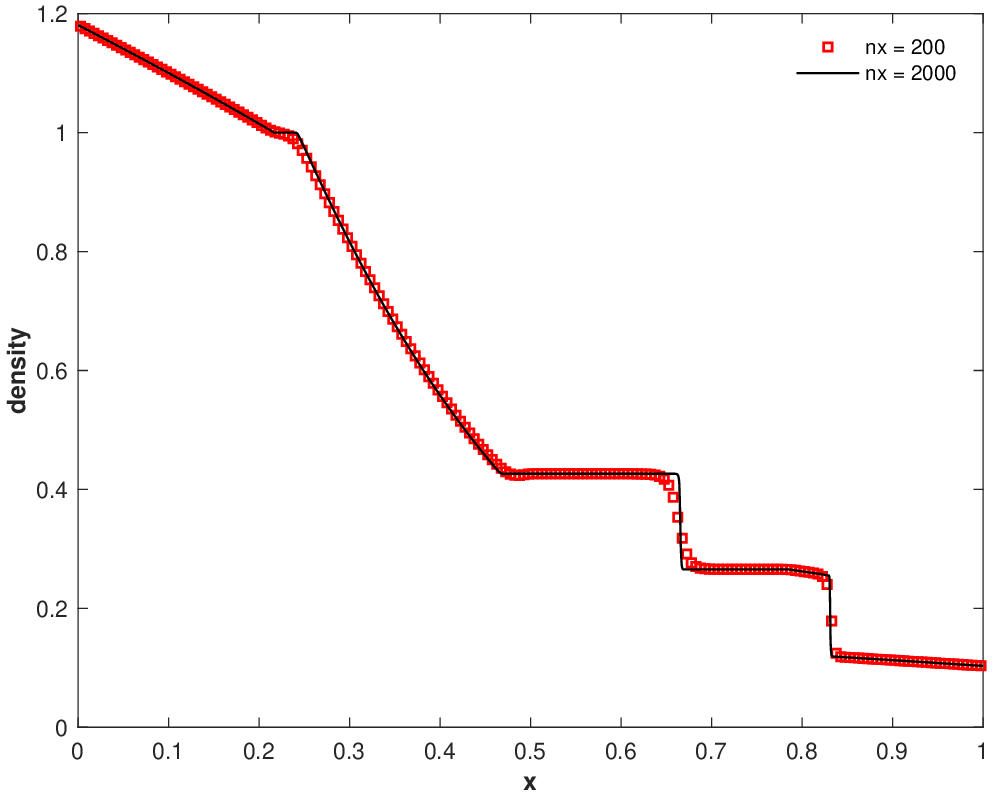}
  }
  \subfigure[velocity]{
  \centering
  \includegraphics[width=5cm,scale=1]{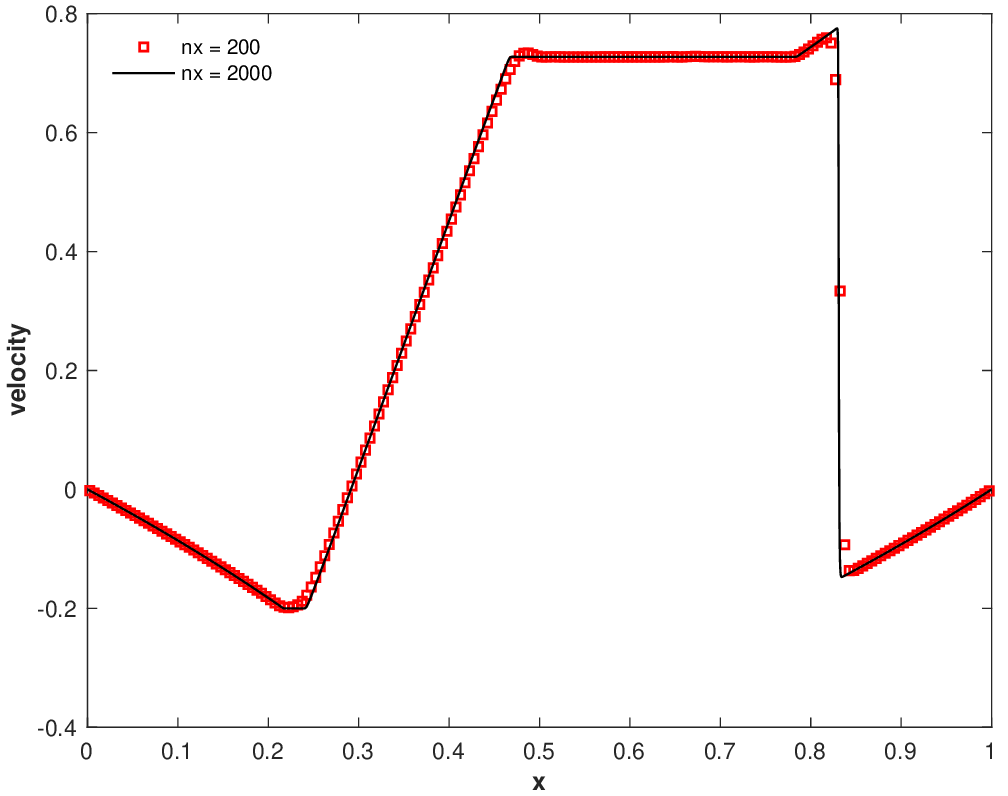}
  }
  \subfigure[energy]{
  \centering
  \includegraphics[width=5cm,scale=1]{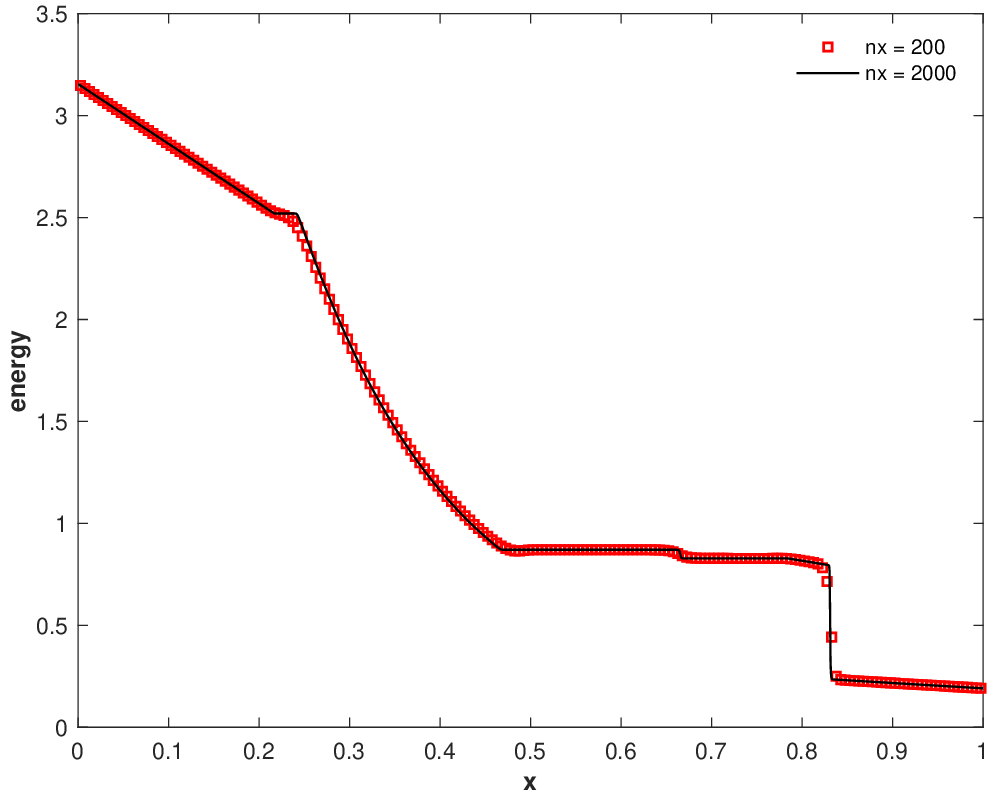}
  }
  \subfigure[pressure]{
  \centering
  \includegraphics[width=5cm,scale=1]{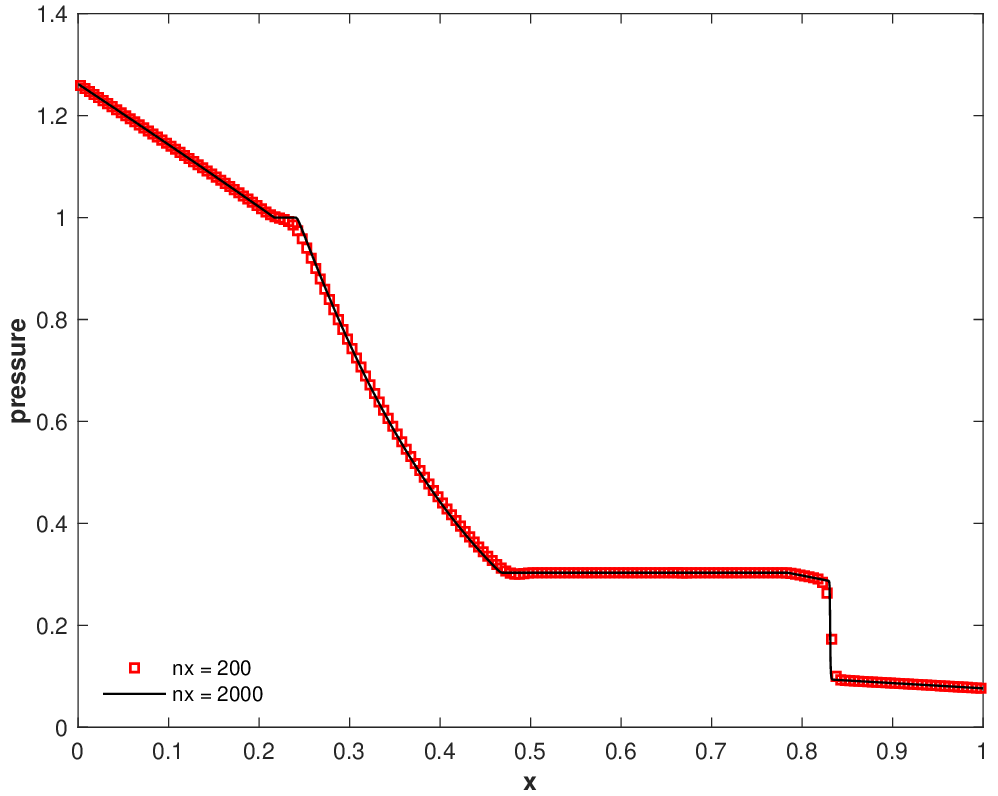}
  }
  \caption{Example \ref{euler:shock}: Shock tube problem with initial value (\ref{euler:sod_ini}) at time $t=0.2$.}\label{euler:st200}
\end{figure}
}

\subsection{Single-layer thermal rotating shallow water equations}\label{subse:trswe}
As a representative example of hyperbolic balance laws, the single-layer thermal rotating shallow water equations can be expressed in the form \eqref{ncs_l}, with the specific definitions given by
\begin{equation}\label{rtswe}
  \boldsymbol{U} = \left[\begin{array}{c}
    h \\
    hu \\
    hv \\
    h\theta \\
    b
  \end{array}\right], 
  \boldsymbol{F}(\boldsymbol{U}) = \left[\begin{array}{c}
    hu \\
    hu^2 + \frac{1}{2}gh^2\theta \\
    huv \\
    hu\theta\\
    0
  \end{array}\right], 
  \mathcal{B}(\boldsymbol {U}) = \left[\begin{array}{ccccc}
    0 & 0 & 0 & 0  & 0\\
    0 & 0 & 0 & 0  & gh\theta\\ 
    0 & 0 & 0 & 0  & 0\\
    0 & 0 & 0 & 0  & 0\\
    0 & 0 & 0 & 0  & 0\\
  \end{array}\right], 
  \boldsymbol{S}(\boldsymbol{U}) = \left[\begin{array}{c}
    0 \\
    fhv \\
    -fhu \\
    0 \\
    0
  \end{array}\right].
\end{equation}
Here, $h,u,v,\theta,g,b$ denote the water height, the two components of horizontal velocity, the potential temperature, the gravitational constant, and the bottom topography function, respectively. The Coriolis parameter $f$ is assumed to vary linearly with $x$. 
This system admits non-trivial steady-state solutions of the form
\begin{equation}\label{mov:sltrswe}
  \begin{aligned}
  \theta &= \text{constant},\quad E = \frac{1}{2}u^2+g\theta(h+b) - \int_{\hat x}^x fv \ d \xi = \text{constant}, \\
  hu &= \text{constant}, \quad \mu = v+\int_{\hat x}^x f \ d \xi = \text{constant}.
  \end{aligned}
\end{equation}

\begin{figure}[!ht]
  \centering
  \subfigure[surface level $h+b$]{
  \centering
  \includegraphics[width=4.3cm,scale=1]{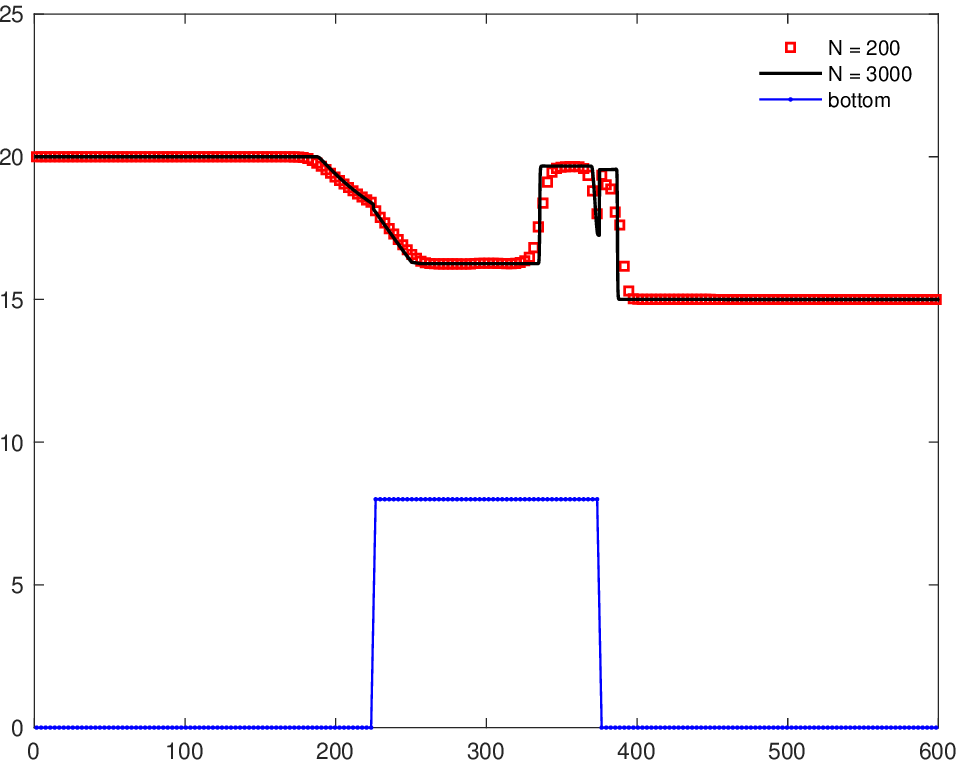}
  }
  \subfigure[discharge $hu$]{
  \centering
  \includegraphics[width=4.3cm,scale=1]{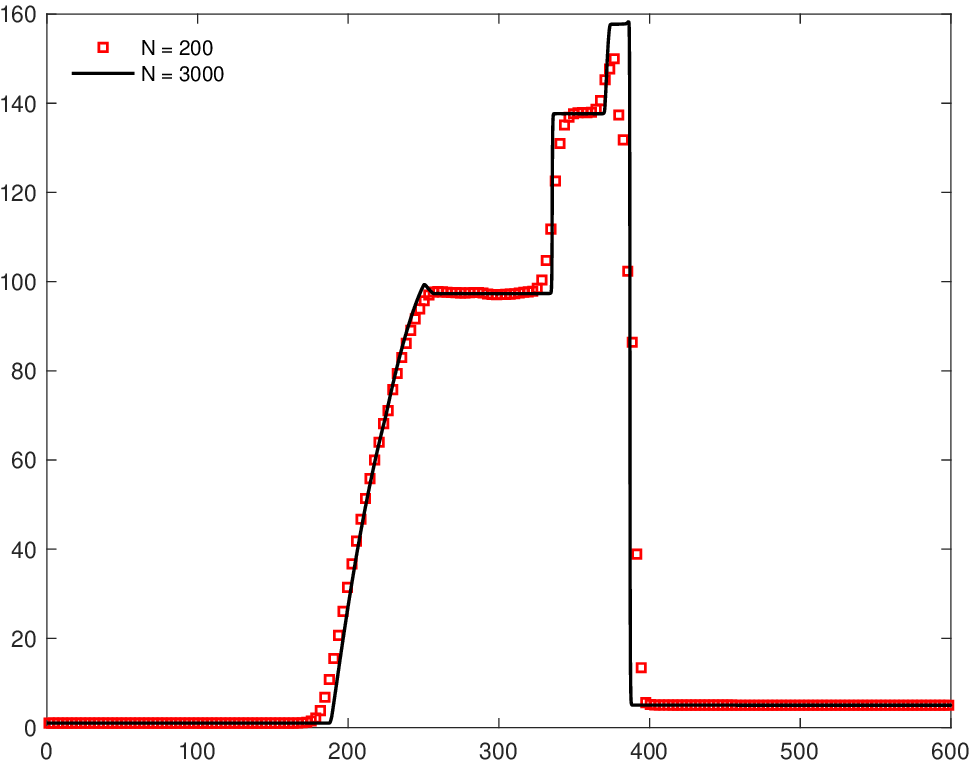}
  }
  \subfigure[velocity $u$]{
  \centering
  \includegraphics[width=4.3cm,scale=1]{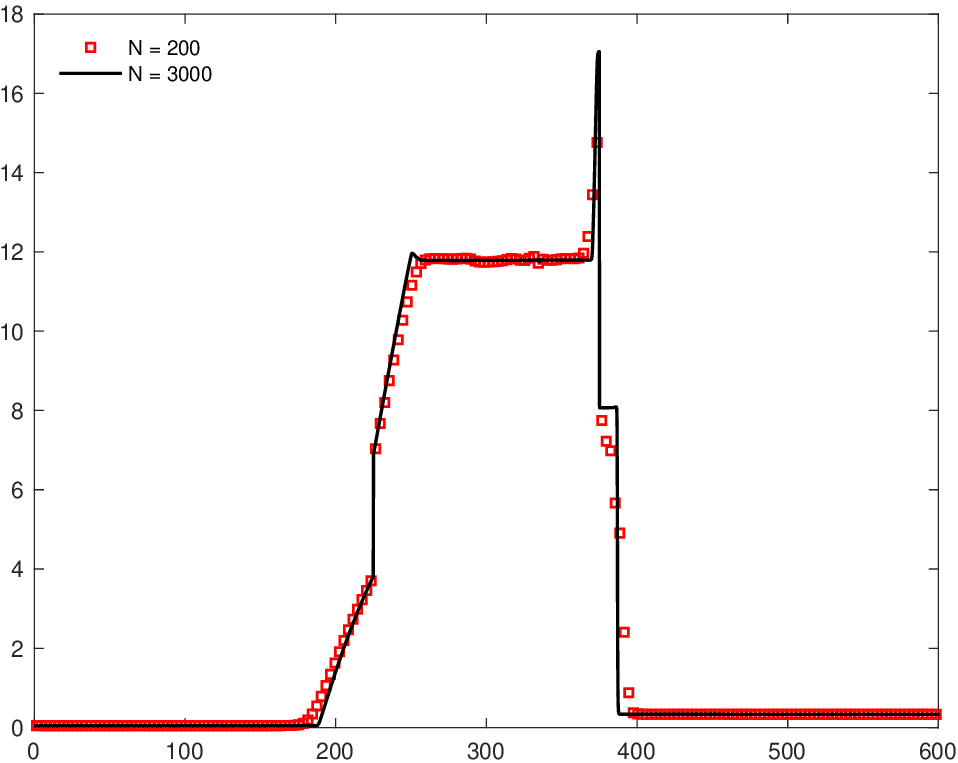}
  }
  \subfigure[ $h\theta$]{
  \centering
  \includegraphics[width=4.3cm,scale=1]{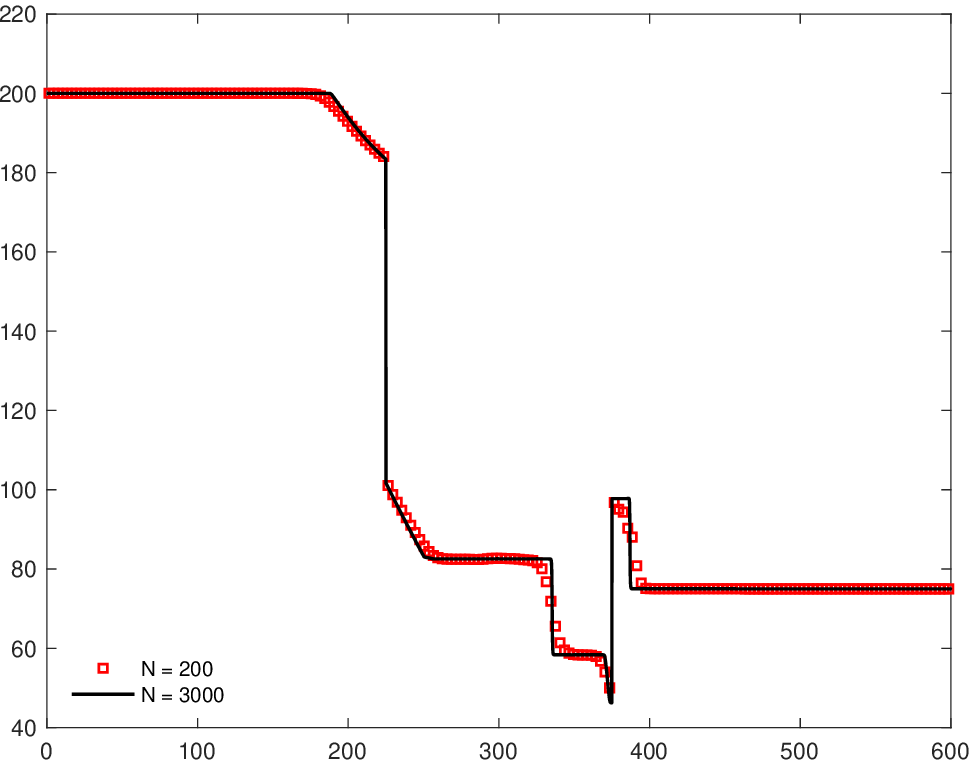}
  }
  \subfigure[temperature $\theta$]{
  \centering
  \includegraphics[width=4.3cm,scale=1]{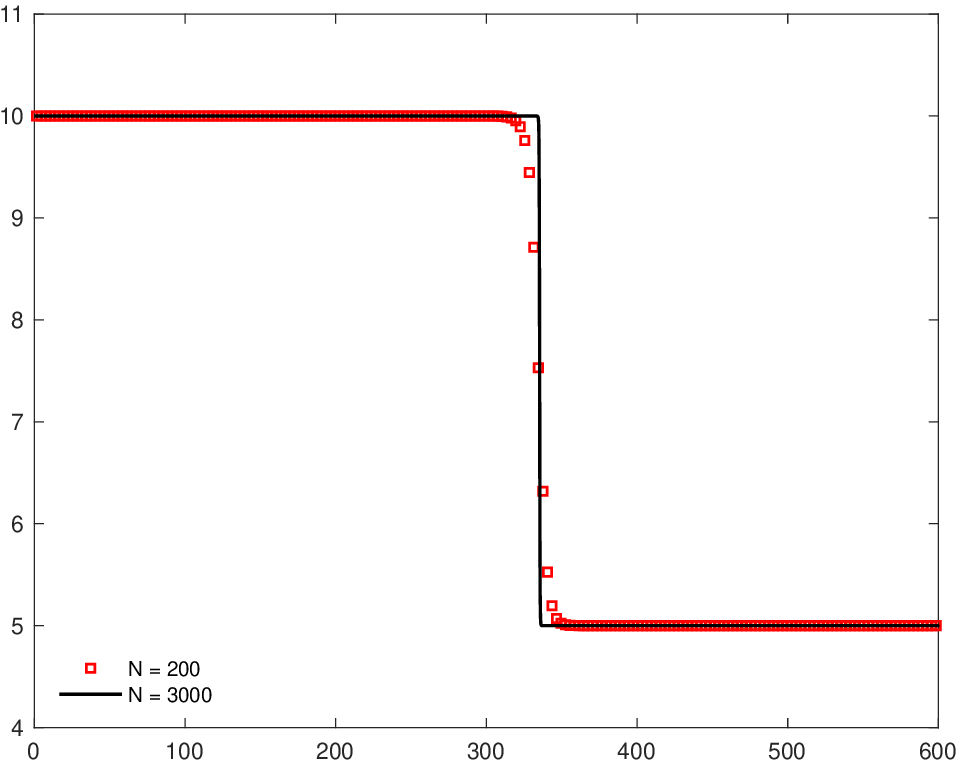}
  }
  \subfigure[pressure $\frac{1}{2}gh^2\theta$]{
  \centering
  \includegraphics[width=4.3cm,scale=1]{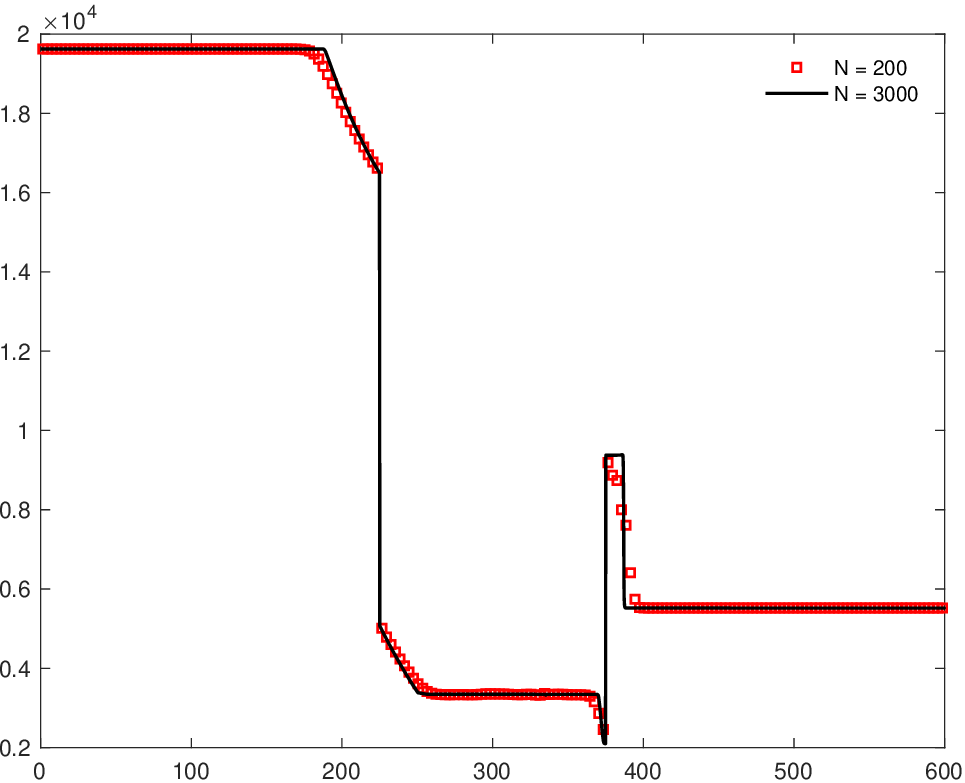}
  }
  \caption{Example \ref{Ripa:dam}: Numerical solutions for the dam breaking problem over a discontinuous bottom topography (\ref{Ripa:bdam}) with the initial value (\ref{Ripa:inidam}) at time $t=3$.}\label{rtswe:rdu}
\end{figure}

\begin{example}{\bf The dam breaking problem over a discontinuous bottom}\label{Ripa:dam}
\end{example}
To assess the performance of the nodal FV-WENO scheme over a nontrivial bottom topography, we first consider a non-rotational case, meaning zero Coriolis forces. We implement a dam-breaking test case, similar to the one used in \cite{2023arXivequilibrium}, over a rectangular bump. This case involves rapidly varying flow over a discontinuous bottom topography on the domain $[0,600]$, with the bottom topography defined as
\begin{equation}\label{Ripa:bdam}
  b(x)=\left\{\begin{array}{lll}
    8,&   \text{if} \ \   225\leq x \leq 375, \\
    0,&\text{otherwise},\\
    \end{array}\right.
\end{equation}
and the initial conditions are given by
\begin{equation}\label{Ripa:inidam}
 (h,u,\theta)(x,0) = \left\{\begin{array}{lll}
    (20-b(x),1,10),&  \text{if} \  x \leq 300, \\
    (15-b(x),5,5),&\text{otherwise}.\\
    \end{array}\right.
\end{equation}
We simulate the scenario up to $t=3$ with transmissive boundary conditions. The results, using 200 uniform cells, are shown in Fig. \ref{rtswe:rdu}. For comparison, we take the solution obtained by the third-order moving equilibrium-preserving DG scheme \cite{2023arXivequilibrium} on a much finer grid of 3000 uniform cells as the reference solution. As noted in \cite{2023arXivequilibrium}, capturing the correct water discharge is more challenging than capturing the water height. The results demonstrate that the well-balanced nodal FV-WENO scheme yields well-resolved, oscillation-free solutions.

\begin{example}{\bf Test for the well-balanced property}\label{wbmov_rtswe}
\end{example}
This example aims to validate the exact preservation of the moving steady state \eqref{mov:sltrswe}, which includes the global integral term. The initial conditions on the computational domain $[0,1]$ are in exact equilibrium:
\begin{equation}\label{mov_ini:rswe}
  (h,u,v,\theta)(x,0) = (h_{\text{eq}}(x),u_{\text{eq}}(x),v_{\text{eq}}(x),\theta_{\text{eq}}(x)) = \left(e^{2x},\frac{1}{2}e^{-2x},-fx,2 \right),
\end{equation}
together with the bottom topography 
\begin{equation*}
  b(x) = -\dfrac{f^2}{4}x^2 - e^{2x} - \dfrac{1}{16}e^{-4x}.
\end{equation*}
We set the Coriolis parameter and the gravitational constant to $f = g = 1$. The $L^1$ and $L^{\infty}$ errors for the conservative variables at $t=5$ with 20 and 100  relatively coarse meshes are shown in Table \ref{mov_wb:rtswe}. It is clear that the well-balanced nodal finite volume scheme can preserve the steady state within machine precision.

\begin{figure}[!htb]
  \centering 
  \subfigure[$h-h_{\text{eq}}$]{
  \centering
     \includegraphics[width= 5.cm,scale=1]{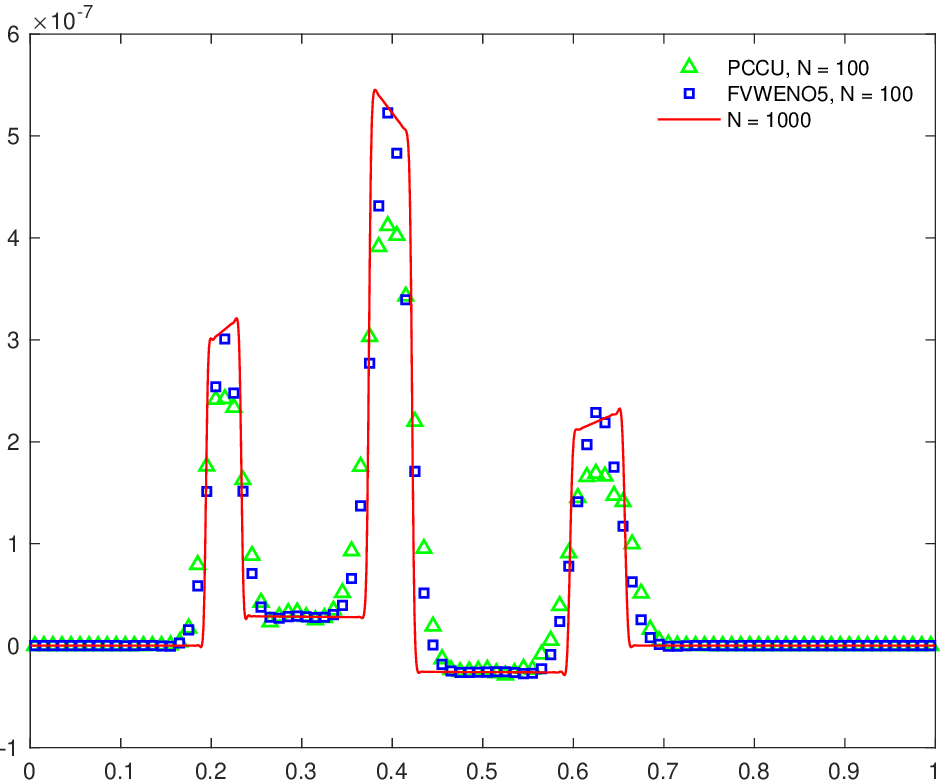}
  }
  \subfigure[$hu-hu_{\text{eq}}$]{
  \centering
     \includegraphics[width= 5.cm,scale=1]{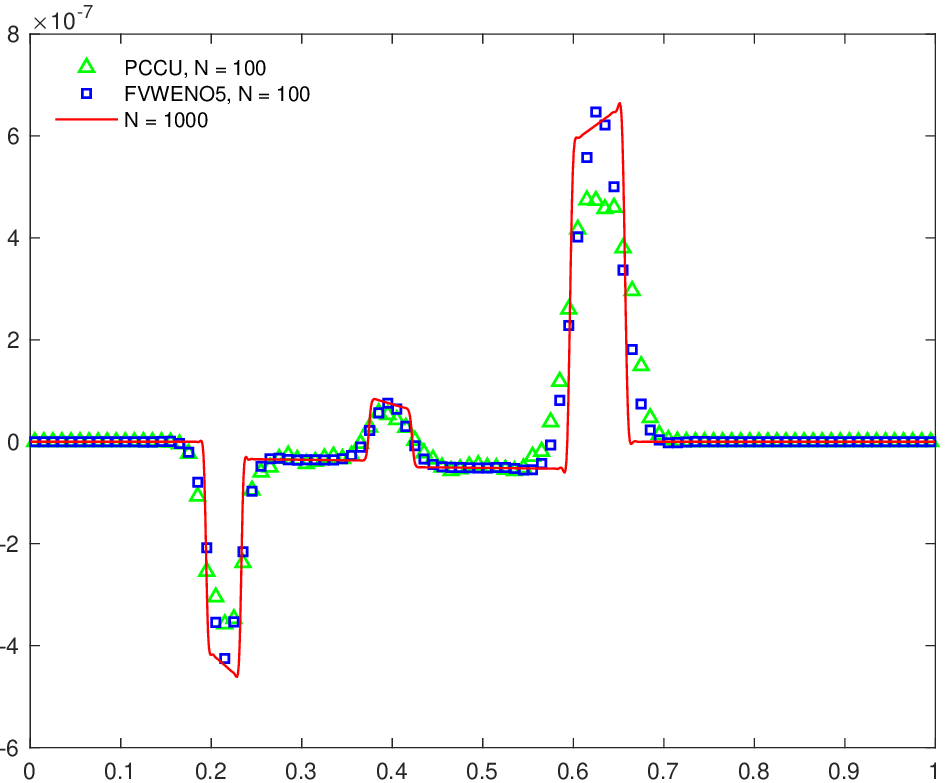}
  } 
  \subfigure[$hv-hv_{\text{eq}}$]{
  \centering
     \includegraphics[width= 5.cm,scale=1]{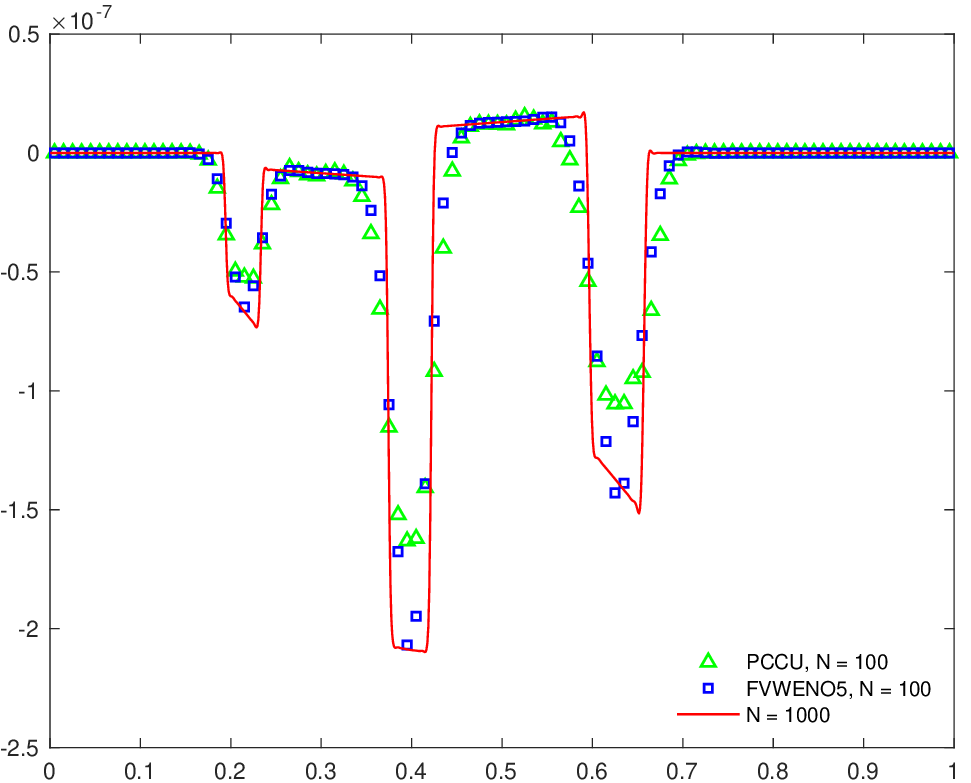}
  }
  \subfigure[$h\theta-h\theta_{\text{eq}}$]{
  \centering
     \includegraphics[width= 5.cm,scale=1]{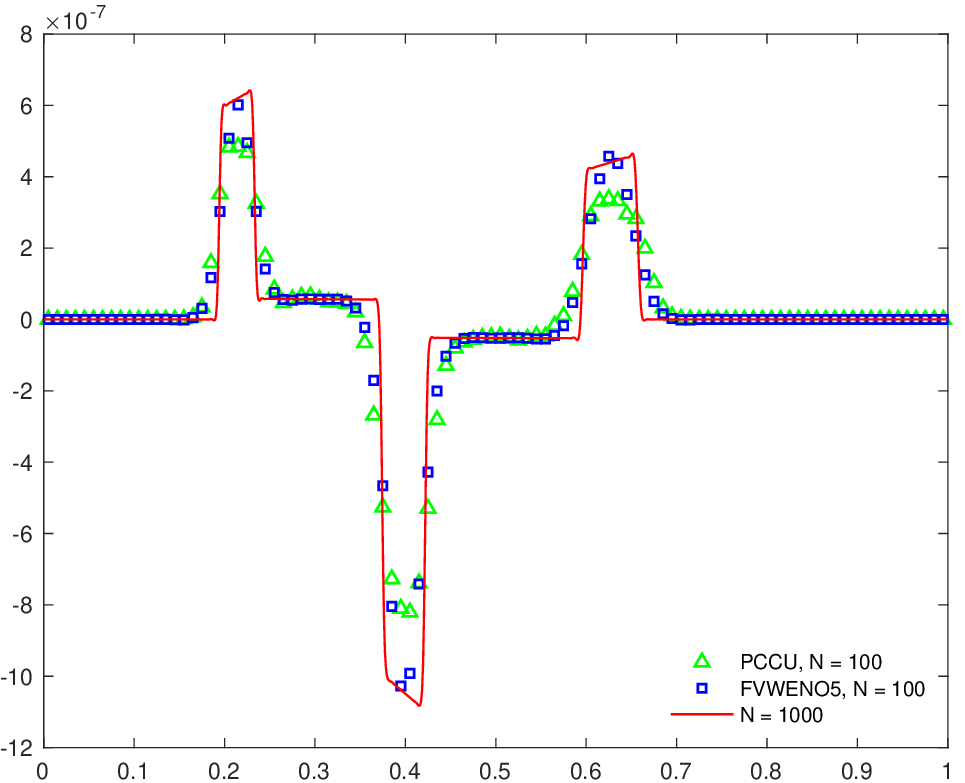}
  } 
  \caption{Example \ref{wbmov_rtswe}: The difference between the moving equilibrium and the numerical solutions computed using the fifth-order nodal FV-WENO scheme and the second-order PCCU scheme.}\label{movper:rtswe}
\end{figure}

Next, we introduce a small disturbance of  $A = 10^{-6}$ to both the water height $h$ and temperature $\theta$, while keeping the velocities unchanged. The initial conditions are given by
\begin{align}\label{movper_ini:rswe}
  (h,u,v,\theta)(x,0) = \left\{\begin{array}{lll}
    (h_{\text{eq}}(x)+A,u_{\text{eq}}(x),v_{\text{eq}}(x),\frac{h_{\text{eq}}(x)\theta_{\text{eq}}(x)}{h_{\text{eq}}(x)+A}) &  \text{if} \ 0.35\leq x \leq 0.4, \\
    (h_{\text{eq}}(x),\qquad u_{\text{eq}}(x),v_{\text{eq}}(x),\theta_{\text{eq}}(x)) &  \text{otherwise}.\\
  \end{array}\right.
\end{align}
Fig. \ref{movper:rtswe} presents the difference between the numerical solutions at $t=0.1$ and the original equilibrium state, with 100 and 1000 uniform cells for comparison. We also compare the performance of our scheme with that of the second-order path-conservative central-upwind (PCCU) scheme developed by Castro D\'iaz et al. \cite{diaz2019path}. As expected, no oscillations are generated, and the fifth-order well-balanced nodal FV-WENO scheme achieves higher resolution compared to the second-order PCCU scheme.

\begin{table}[htb!]
  \centering
  \caption{Example \ref{wbmov_rtswe}: $L^1$ and $L^{\infty}$ errors for the moving equilibrium state \eqref{mov_ini:rswe}.}\label{mov_wb:rtswe}
  \begin{tabular}{c c c c c c c  }
   \toprule
    $N$ &{error}&{$h$}&{$hu$}&{$hv$}&{$h\theta$}\\
    \midrule
    \multirow{2}{*}{20} &$L^1$        & 1.02E-14 & 7.76E-14 & 8.97E-15 & 1.16E-14  \\
                  ~     &$L^{\infty}$ & 2.82E-14 & 8.58E-14 & 4.09E-14 & 4.97E-14  \\
    \cmidrule(lr){1-6}
    \multirow{2}{*}{100}&$L^1$        & 3.87E-14 & 1.05E-13 & 2.38E-14 & 5.63E-14 \\
                  ~     &$L^{\infty}$ & 8.44E-14 & 1.50E-13 & 1.06E-13 & 2.43E-13  \\
    \bottomrule
  \end{tabular}
 \end{table}

\subsection{Two-layer rotating shallow water equations}\label{subse:2tlrswe}
The two-layer rotating shallow water equations fall into the category of non-conservative hyperbolic systems and take the form of \eqref{ncs_l}, with the specific definition
\begin{scriptsize}
\begin{equation}\label{2lrswe}
  \boldsymbol{U} = \left[\begin{array}{c}
    h_1 \\
    h_1u_1 \\
    h_1v_1 \\
    h_2 \\
    h_2u_2 \\
    h_2v_2 \\
    b
  \end{array}\right], 
  \boldsymbol{F}(\boldsymbol{U}) = \left[\begin{array}{c}
    h_1u_1 \\
    h_1u_1^2 + \frac{1}{2}gh_1^2 \\
    h_1u_1v_1 \\
    h_2u_2 \\
    h_2u_2^2 + \frac{1}{2}gh_2^2 \\
    h_2u_2v_2 \\
    0
  \end{array}\right], 
  \mathcal{B}(\boldsymbol {U}) = \left[\begin{array}{ccccccc}
    0 & 0 & 0 & 0  & 0 & 0  & 0\\
    0 & 0 & 0 & gh_1& 0 & 0  & gh_1\\ 
    0 & 0 & 0 & 0  & 0 & 0  & 0\\
    0 & 0 & 0 & 0  & 0 & 0  & 0\\
    grh_2 & 0 & 0 & 0  & 0 & 0  & gh_2\\
    0 & 0 & 0 & 0  & 0 & 0  & 0\\
    0 & 0 & 0 & 0  & 0 & 0  & 0\\
  \end{array}\right], 
  \boldsymbol{S}(\boldsymbol{U}) = \left[\begin{array}{c}
    0 \\
    fh_1v_1 \\
    -fh_1u_1 \\
    0 \\
    fh_2v_2 \\
    -fh_2u_2 \\
    0
  \end{array}\right],
\end{equation}
\end{scriptsize}
whose non-trivial steady states are given by 
  \begin{align*}
    h_1u_1 &= \text{constant}, \quad  \mu_1 = v_1 +\int_{\hat x}^x f \ d \xi  = \text{constant}, \\
    h_2u_2&= \text{constant}, \quad \mu_2 = v_2 +\int_{\hat x}^x f \ d \xi  = \text{constant},\\
      E_1:&= \frac{1}{2}u_1^2 + g(h_1 + h_2 + b)  -\int_{\hat x}^x fv_1 \ d \xi= \text{constant},  \\
  E_2:&= \frac{1}{2}u_2^2 + g(r h_1 + h_2 + b) -\int_{\hat x}^x fv_2 \ d \xi= \text{constant}.
  \end{align*}
Here, $h_i,u_i,v_i,E_i$ represent the thickness, velocities and energy in each layer, with subscript 1 referring to the upper layer and subscript 2 referring to the lower one. The function $b$ is the nonflat bottom topography and $g$ is the gravitational constant. Each layer is assumed to have a constant density $\rho_i$ and $r:={\rho_1}/{\rho_2}<1$ is the density ratio.

\begin{example}{\bf Accuracy test}\label{acc1D_2lrswe}
\end{example}
We use this example to validate the high-order accuracy of our developed finite volume WENO scheme. To this end, a manufactured solution is constructed by introducing additional source terms into the governing equation \eqref{ncs_l}:
\begin{equation}\label{add_so}
  \boldsymbol{U}_t + \boldsymbol{F}(\boldsymbol{U})_x + \mathcal{B}(\boldsymbol {U})\boldsymbol {U}_x = \boldsymbol{S}(\boldsymbol{U}) + \widetilde{\boldsymbol{S}}(x,t).
\end{equation} 
Periodic boundary conditions are applied here, with the computational domain set as [0,2]. We take the gravitational constant $g=1$, Coriolis parameter $f=1$ and the density ratio $r = 0.7$. The bottom topography is given by
\begin{equation*}
  b(x) = 1.5 + \sin(\pi x).
\end{equation*}
The following exact solutions are considered
 \begin{align}
   &h_1(x,t) = 6+\cos(\pi x)\cos(\pi t), \quad h_1u_1(x,t) = \sin(\pi x)\sin(\pi t), \nonumber\\
   &h_1v_1(x,t) = 6+\cos(\pi x)\cos(\pi t),\quad h_2(x,t) = 4+\cos(\pi x)\cos(\pi t), \label{iniadd:2lrswe}\\
   &h_2u_2(x,t) = \sin(\pi x)\sin(\pi t),\quad h_2v_2(x,t) = 4+\cos(\pi x)\cos(\pi t). \nonumber
 \end{align}
We compute until $t=0.1$ and the $L^1$ error and the order of accuracy for the conservative variables are displayed in Table \ref{FVacc1D_add_2lrswe}. From the results, we can observe the anticipated fifth-order optimal convergence rates.
\begin{table}[!htb]
  \centering
  \caption{Example \ref{acc1D_2lrswe}:  $L^1$ errors and numerical orders of accuracy with initial condition (\ref{iniadd:2lrswe}), $t = 0.1$.}\label{FVacc1D_add_2lrswe}
  \begin{tabular}{c c c c c c c c c }
  \toprule
    \multirow{2}{*} {$N$} & \multicolumn{2}{c}{$h_1$}&\multicolumn{2}{c}{$ h_1u_1 $}& \multicolumn{2}{c}{$h_1v_1$} \\
   \cmidrule(lr){2-3} \cmidrule(lr){4-5} \cmidrule(lr){6-7} 
   ~ &$L^1$ error &order&$L^1$ error &order&$L^1$ error &order\\
   \midrule
   25  & 1.25E-04 &   -- & 4.71E-04 &   -- & 1.44E-04 &   --  \\
   50  & 3.64E-06 & 5.11 & 1.47E-05 & 5.00 & 4.03E-06 & 5.16  \\
   100 & 1.10E-07 & 5.05 & 4.49E-07 & 5.03 & 1.18E-07 & 5.10  \\
   200 & 3.57E-09 & 4.94 & 1.40E-08 & 5.01 & 3.75E-09 & 4.97  \\
   400 & 1.18E-10 & 4.91 & 4.51E-10 & 4.95 & 1.23E-10 & 4.93  \\
   800 & 3.69E-12 & 5.00 & 1.44E-11 & 4.97 & 4.02E-12 & 4.94  \\
   \midrule
   \multirow{2}{*} {$N$} & \multicolumn{2}{c}{$h_2$}&\multicolumn{2}{c}{$ h_2u_2 $}& \multicolumn{2}{c}{$h_2v_2$}\\
   \cmidrule(lr){2-3} \cmidrule(lr){4-5} \cmidrule(lr){6-7}
   ~ &$L^1$ error &order&$L^1$ error &order&$L^1$ error &order\\
   \midrule
   25  & 2.12E-04 &   -- & 4.10E-04 &   -- & 2.05E-04 &   --  \\
   50  & 6.17E-06 & 5.11 & 1.22E-05 & 5.07 & 6.02E-06 & 5.09  \\
   100 & 1.99E-07 & 4.95 & 3.52E-07 & 5.12 & 1.96E-07 & 4.94  \\
   200 & 6.49E-09 & 4.94 & 1.08E-08 & 5.02 & 6.38E-09 & 4.94  \\
   400 & 2.13E-10 & 4.93 & 3.69E-10 & 4.87 & 2.09E-10 & 4.93  \\
   800 & 6.92E-12 & 4.94 & 1.26E-11 & 4.87 & 6.56E-12 & 5.00  \\
   \bottomrule
  \end{tabular}
 \end{table}

\begin{figure}[!ht]
  \centering  
  \subfigure[water height $h_1$]{
    \centering
    \includegraphics[width=4.3cm,scale=1]{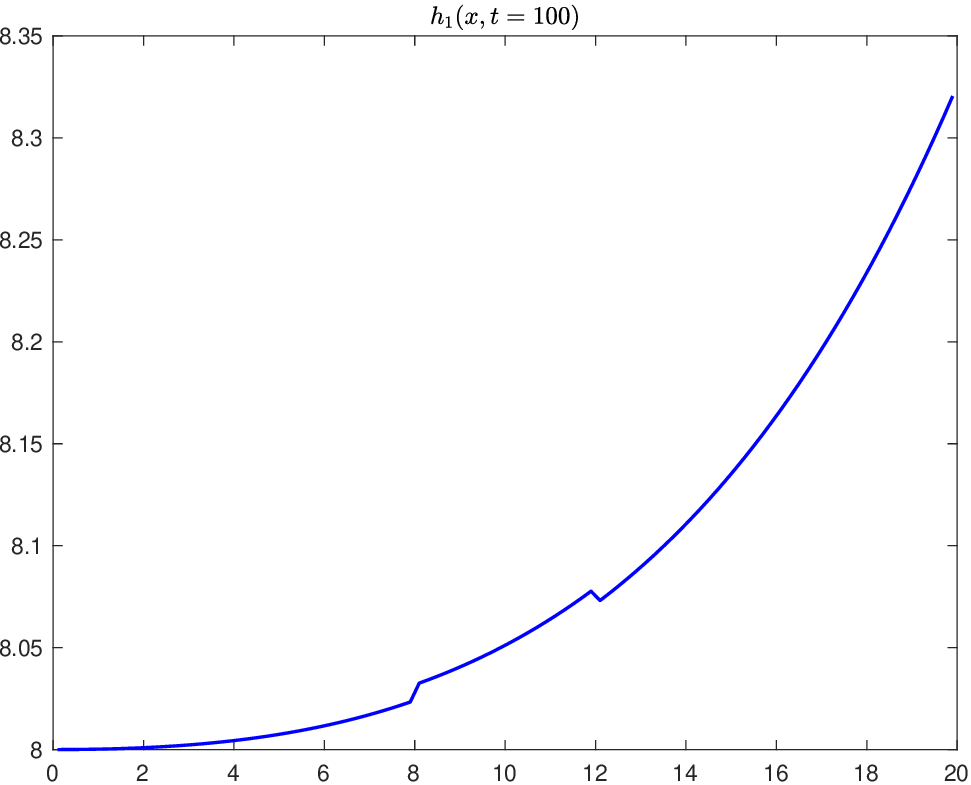}
  }
   \subfigure[discharge $h_1u_1$]{
    \centering
    \includegraphics[width= 4.3cm,scale=1]{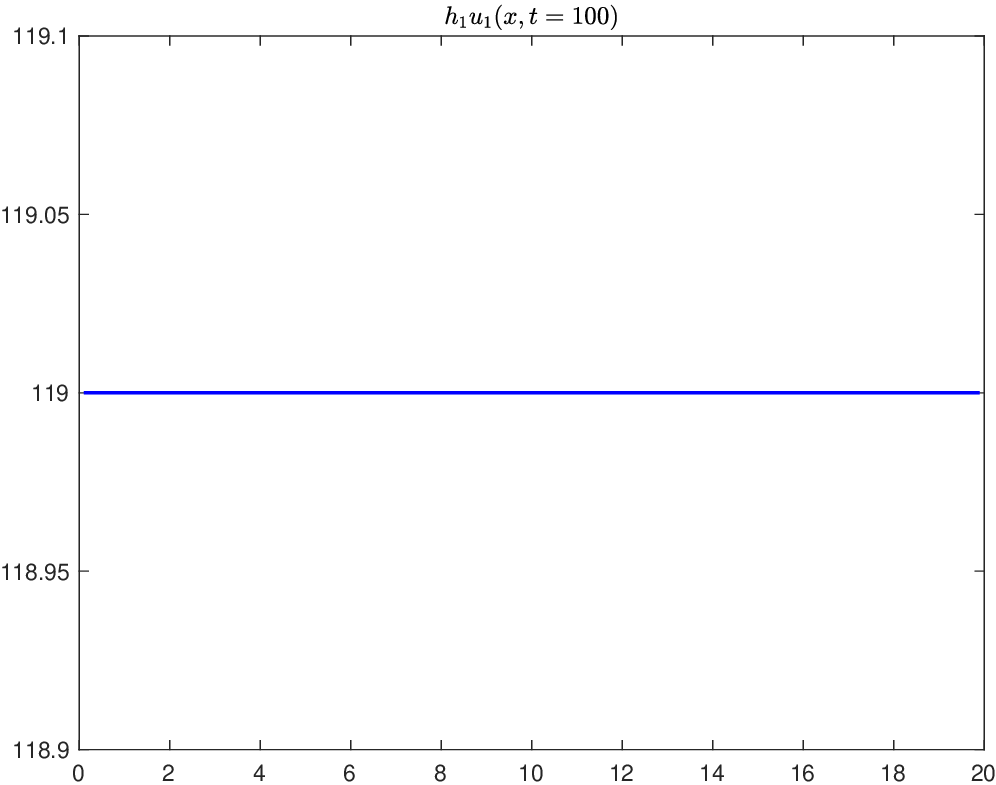}
  }
  \subfigure[energy $E_1$]{
    \centering
    \includegraphics[width= 4.3cm,scale=1]{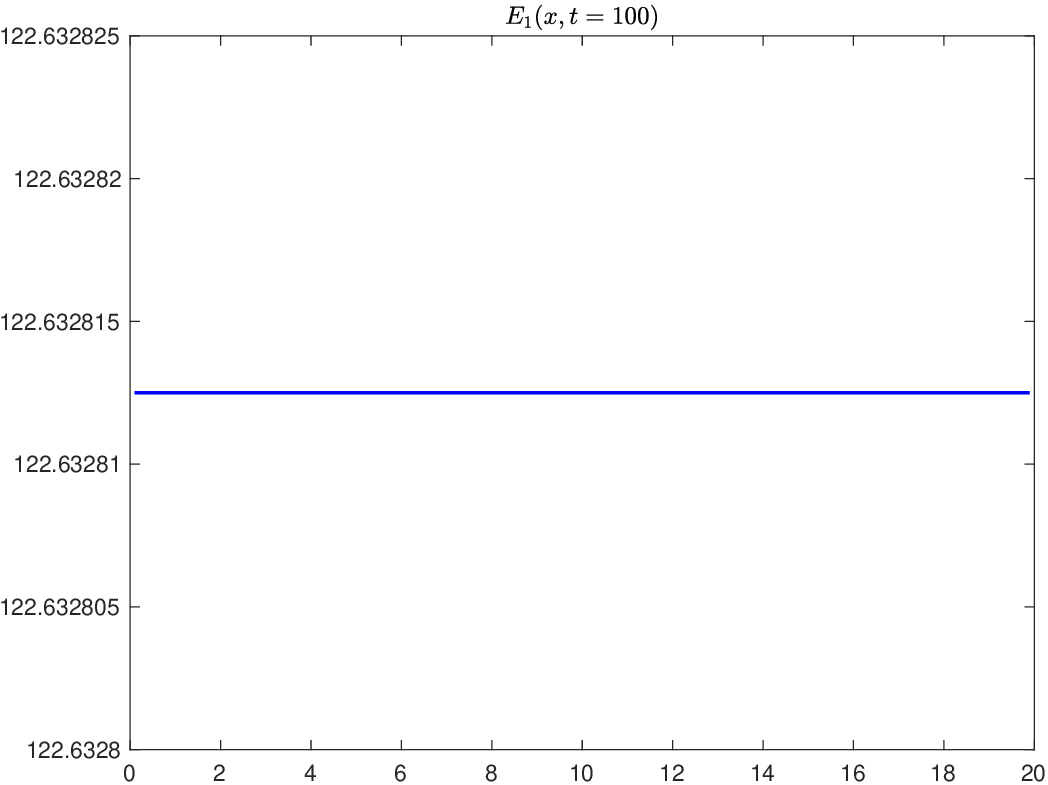}
  }
  \subfigure[water height $h_2$]{
    \centering
    \includegraphics[width= 4.3cm,scale=1]{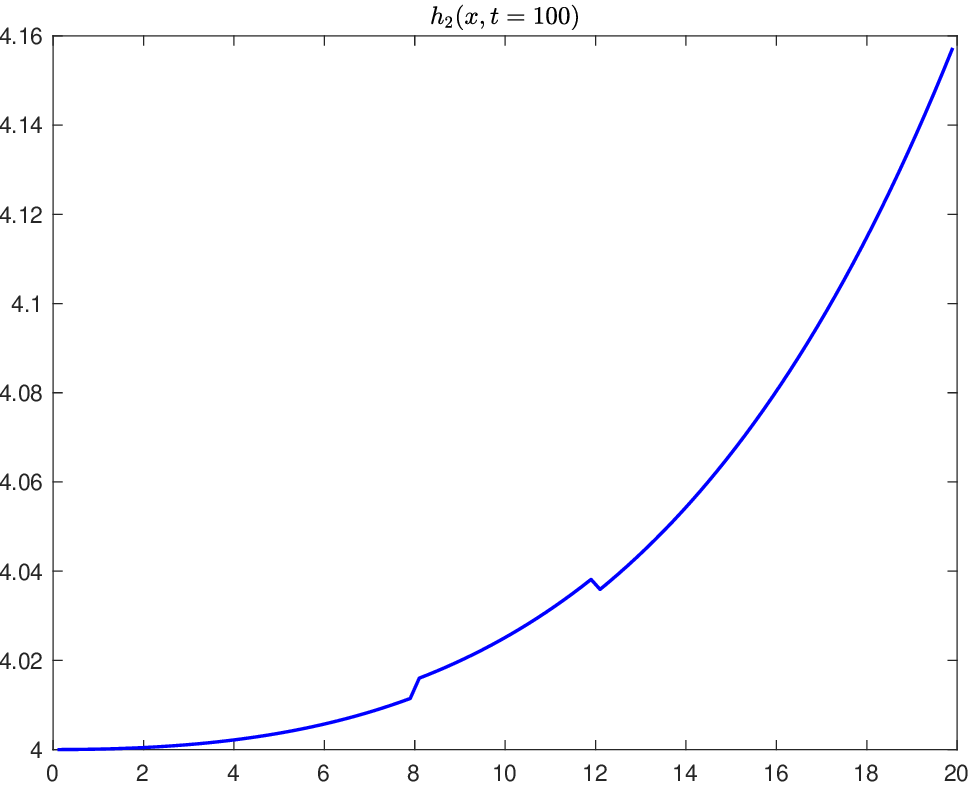}
  }
  \subfigure[discharge $h_2u_2$]{
    \centering
    \includegraphics[width= 4.3cm,scale=1]{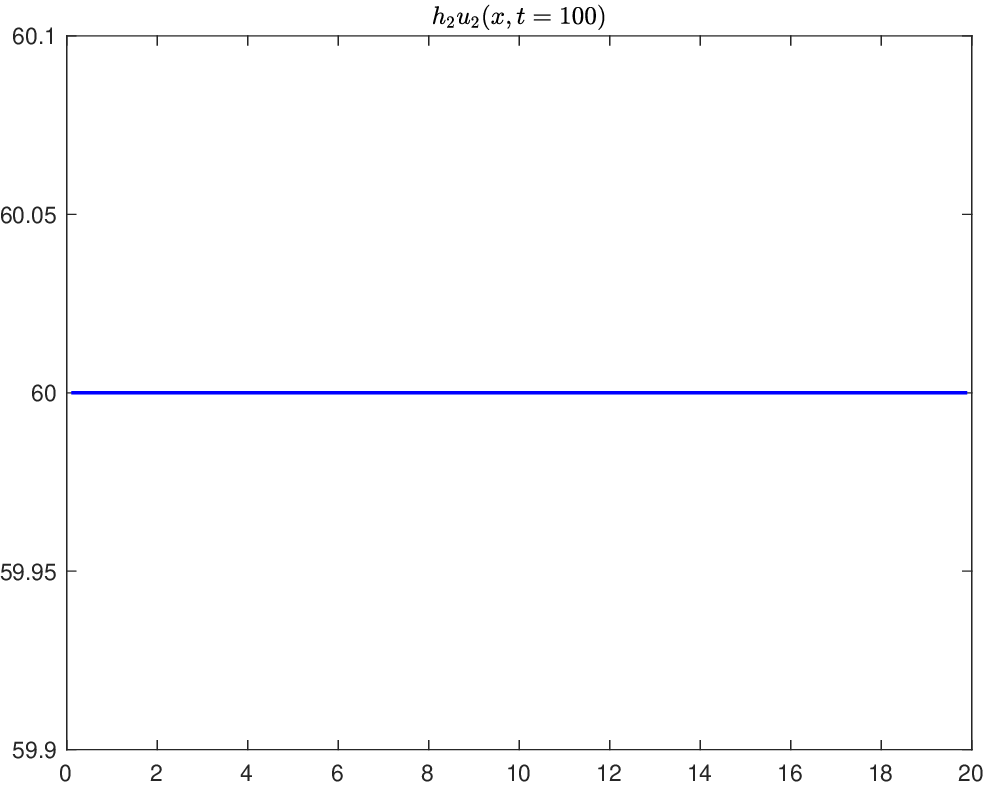}
  }
   \subfigure[energy $E_2$]{
    \centering
    \includegraphics[width= 4.3cm,scale=1]{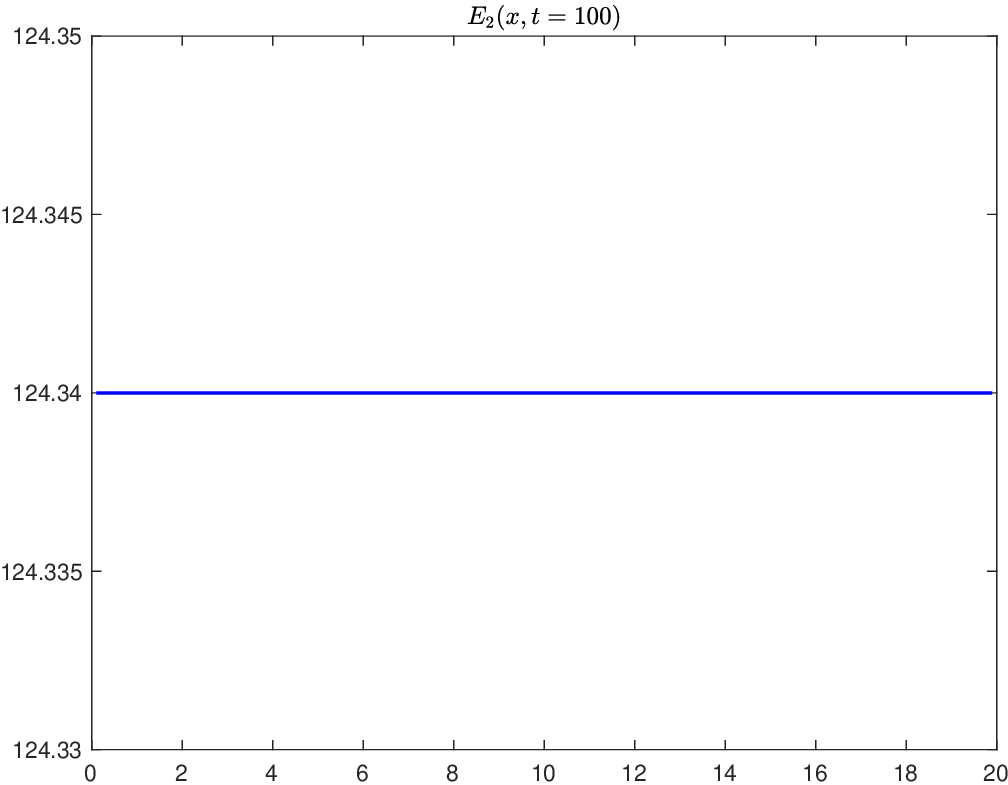}
  }
  \subfigure[velocity $v_1$]{
    \centering
    \includegraphics[width= 4.3cm,scale=1]{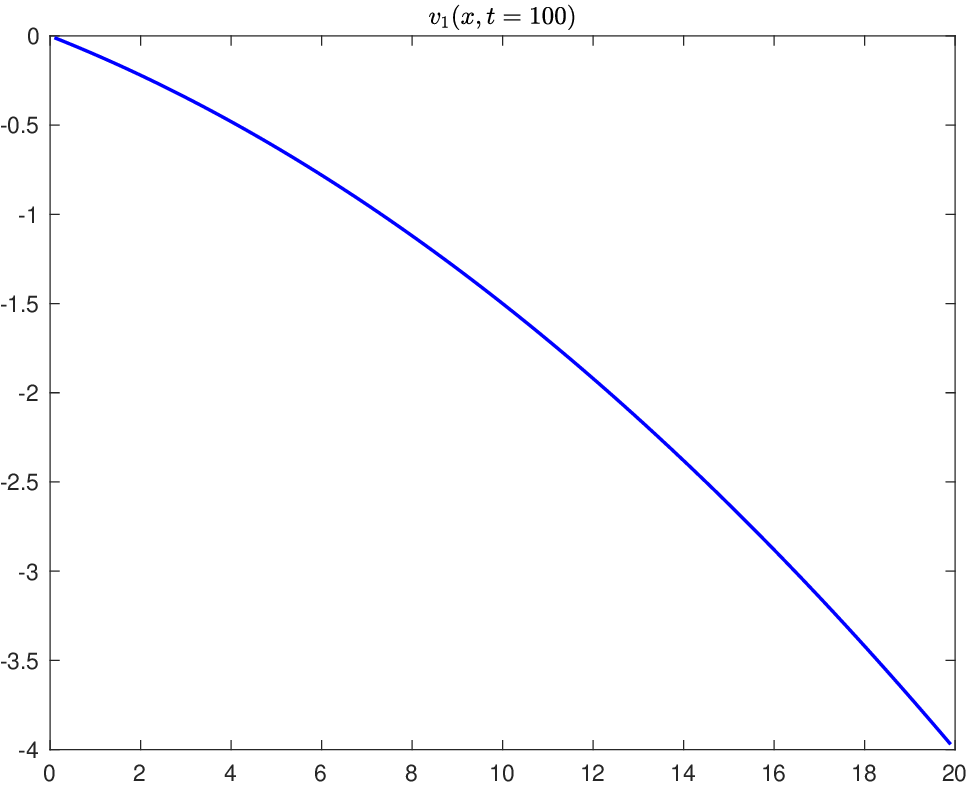}
  }
  \subfigure[velocity $v_2$]{
    \centering
    \includegraphics[width= 4.3cm,scale=1]{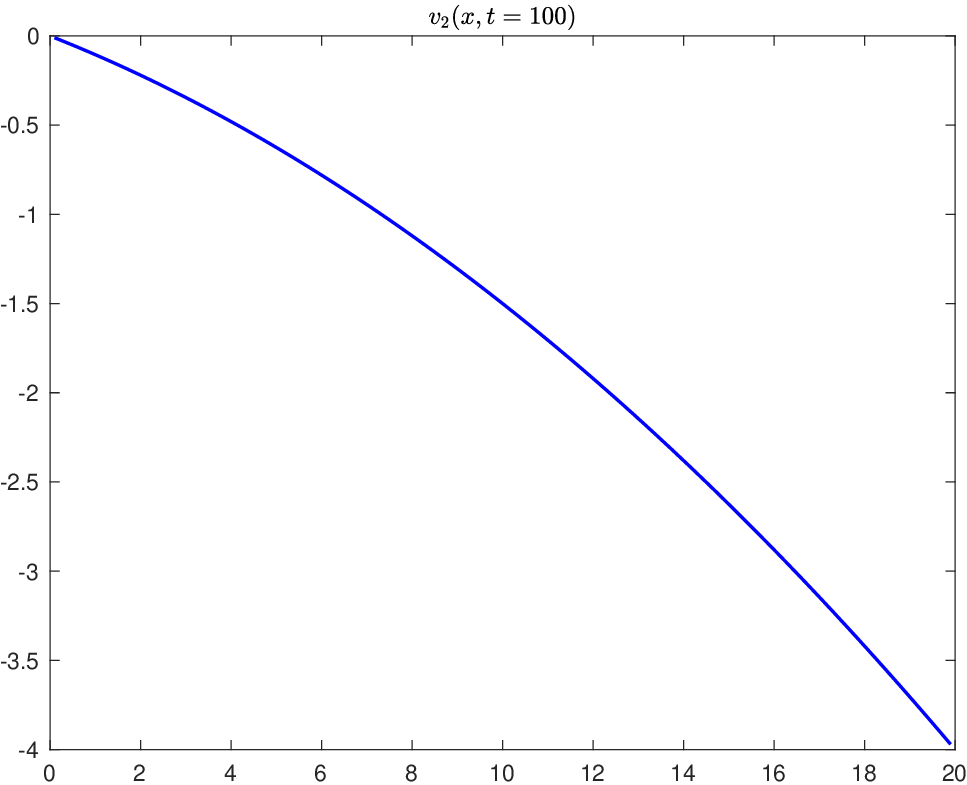}
  }
  \caption{Example \ref{conv_2lrswe}: The discrete steady state of water height, discharge, energy and velocity in the upper layer (left) and lower layer (right), with 100 cells.}\label{conv_dis}
\end{figure}

\begin{example}{\bf Covergence to steady states}\label{conv_2lrswe}
\end{example}
In this test, we adopt a similar example to that in \cite{kurganov2023well} to study the convergence towards steady flows over a discontinuous bottom:
\begin{equation*}
  b(x) = \left\{\begin{array}{lll}
    0.2, &  \text{if} \  8\leq x \leq 12, \\
    0,   &  \text{otherwise},\\
    \end{array}\right.
\end{equation*}
on the computational domain [0,20]. The initial conditions are given by
\begin{equation}\label{conv_ini:rswe}
  h_1(x,0) = 8,\ h_2(x,0) = 4,\ h_1u_1(x,0) = h_2u_2(x,0) = h_1v_1(x,0) =h_2v_2(x,0)= 0,
\end{equation}
with homogeneous Neumann boundary conditions at $x=20$ and Dirichlet boundary conditions at $x=0$,
\begin{equation*}
  \begin{aligned}
  &h_1(0,t) = 8,\ h_1u_1(0,t)= 119, \ h_1v_1(0,t) = 0,\\
  &h_2(0,t) = 4,\ h_2u_2(0,t)= 60, \ h_2v_2(0,t) = 0.
  \end{aligned}
\end{equation*}
We set the Coriolis parameter $f = 0.1+0.01x$, the density ratio $r=0.98$ and the gravitational constant $g=1$. Fig. \ref{conv_dis} shows the numerical solutions computed using our well-balanced nodal FV-WENO scheme based on nodal values, with 100 uniform mesh cells, over a large final time $t=100$. By this time, the solutions have already converged to the corresponding discrete moving steady state.

\begin{figure}[htb]
  \centering 
  \subfigure[$h_1-(h_1)_{\text{eq}}$, $t=0.3$]{
    \centering
    \includegraphics[width= 4.3cm,scale=1]{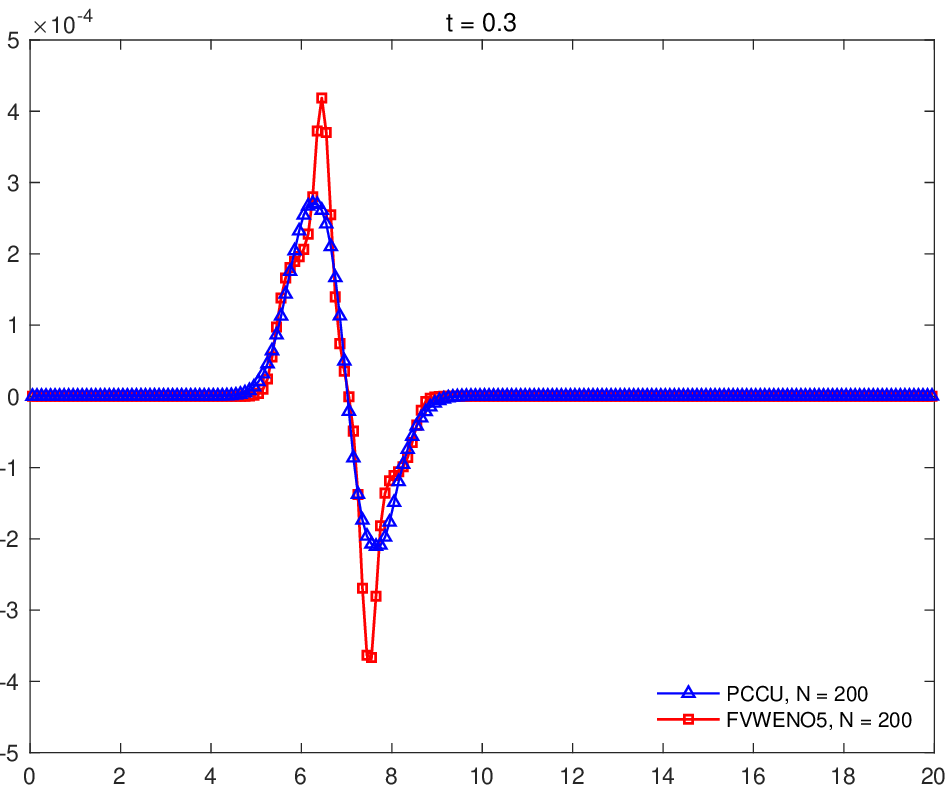}
  }
  \subfigure[$h_1u_1-(h_1u_1)_{\text{eq}}$, $t=0.3$]{
    \centering
    \includegraphics[width= 4.3cm,scale=1]{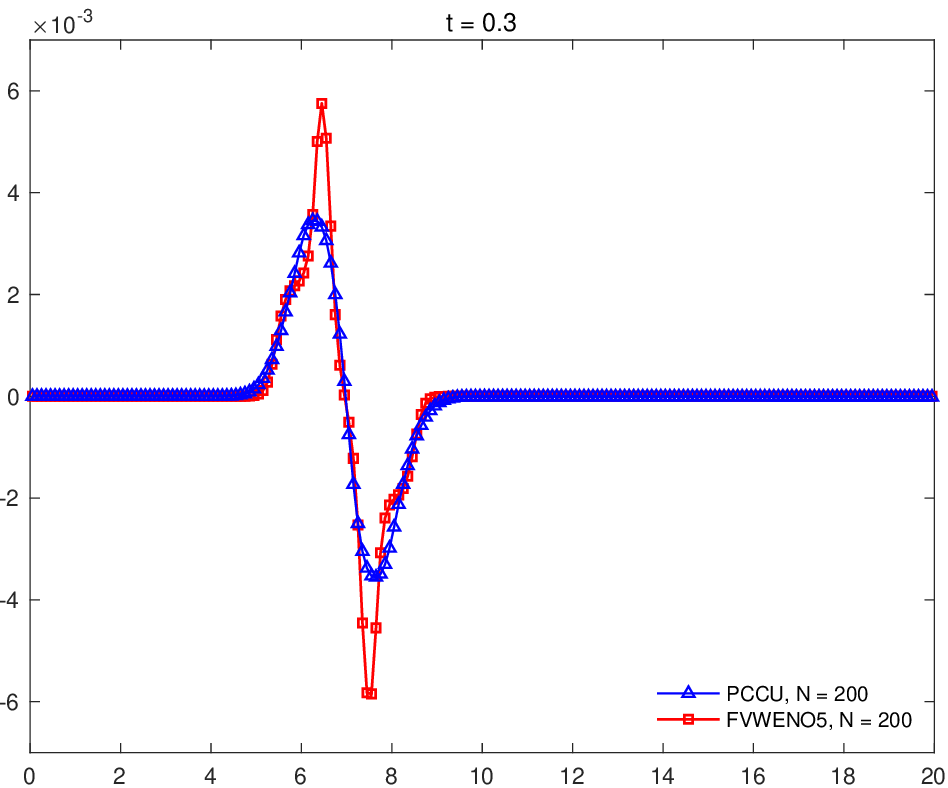}
  }
  \subfigure[$h_1v_1-(h_1v_1)_{\text{eq}}$, $t=0.3$]{
    \centering
    \includegraphics[width= 4.3cm,scale=1]{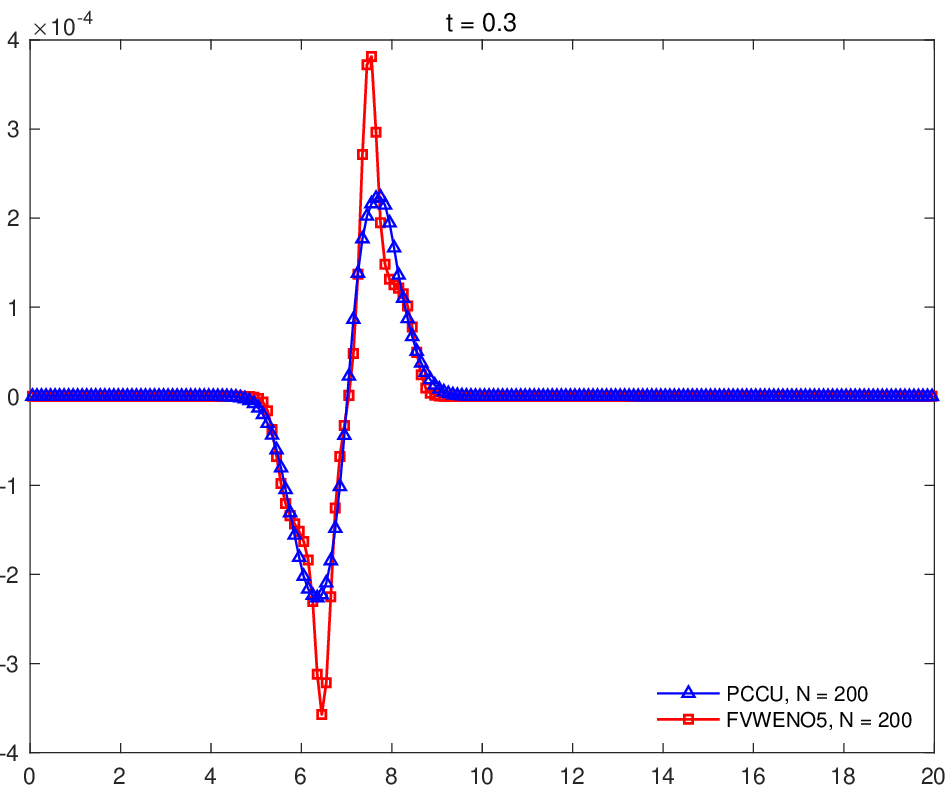}
  }
  \subfigure[$h_1-(h_1)_{\text{eq}}$, $t=0.6$]{
    \centering
    \includegraphics[width= 4.3cm,scale=1]{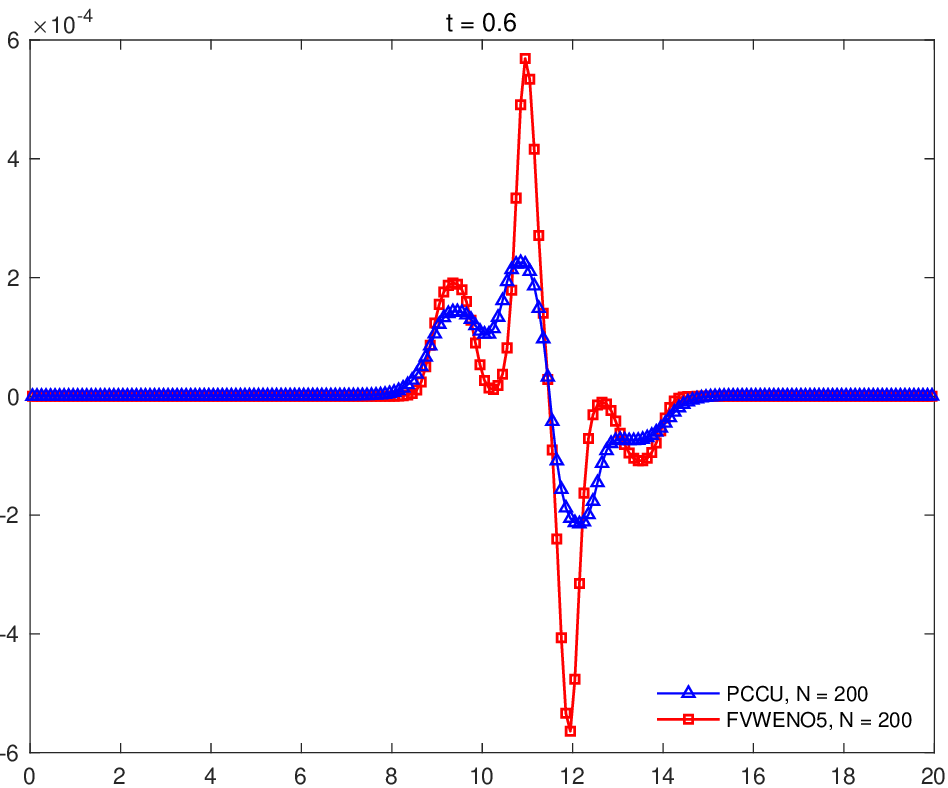}
  }
  \subfigure[$h_1u_1-(h_1u_1)_{\text{eq}}$, $t=0.6$]{
    \centering
    \includegraphics[width= 4.3cm,scale=1]{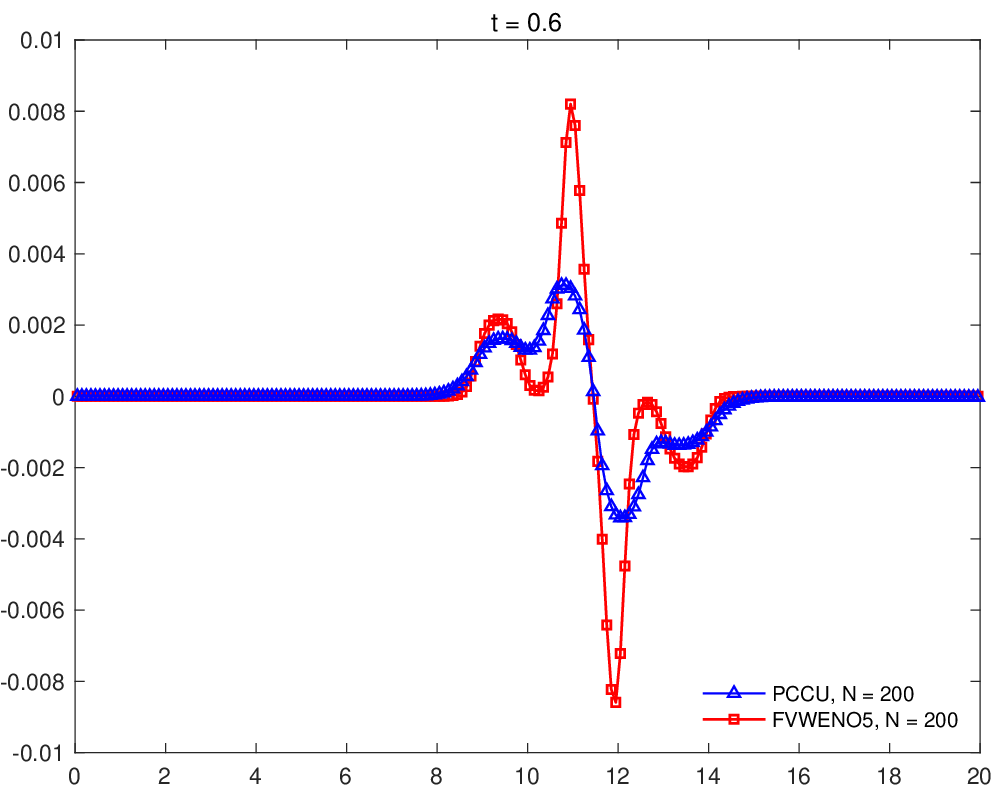}
  }
  \subfigure[$h_1v_1-(h_1v_1)_{\text{eq}}$, $t=0.6$]{
    \centering
    \includegraphics[width= 4.3cm,scale=1]{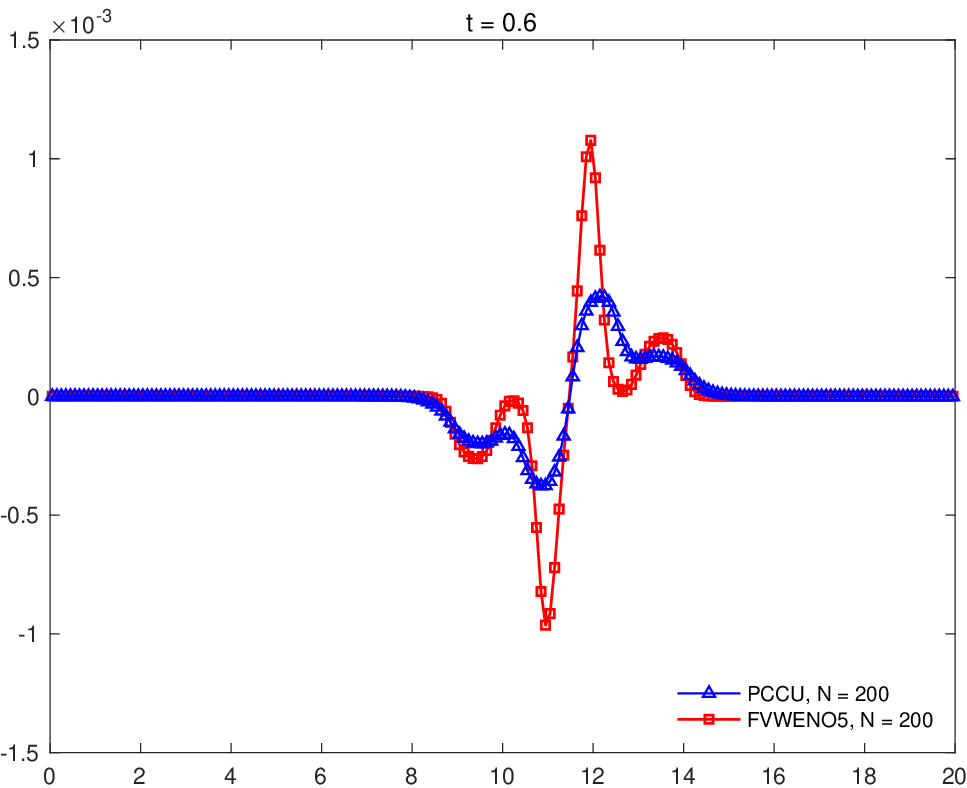}
  } 
  \caption{Example \ref{conv_2lrswe}: The difference between the numerical solutions and the original moving equilibrium at different times $t=0.3$ (top) and $t=0.6$ (bottom). }\label{conv_com}
\end{figure}

To investigate the ability to capture the propagation of the equilibrium state, we initialize the discrete steady-state solution using the numerically converged results and introduce a small perturbation $A = 10^{-4}$ to the upper layer:
\begin{equation*}
  h_1(x,0) = (h_1)_{\text{eq}}(x) +
  \left\{\begin{array}{lll}
  A, &  \text{if} \ 2\leq x \leq 3, \\
  0,      &  \text{otherwise}.\\
  \end{array}\right.
\end{equation*}
We evaluate the performance of our developed fifth-order scheme and the second-order PCCU scheme, both discretized with 200 uniform cells, by comparing their solutions at different times $t = 0.3$ and $0.6$. Fig. \ref{conv_com} shows the temporal evolution of the difference between the numerical solutions and the original moving equilibrium. Due to the similar results for the lower layer, we only present the corresponding results for the upper layer.
Our scheme successfully achieves high resolution and effectively captures small features of the flow, owing to its high-order accuracy and well-balanced property.

\subsection{Two-dimensional shallow water equations}
\begin{figure}[htbp]
  \centering  
  {\centering
     \includegraphics[width= 5.5cm,scale=1]{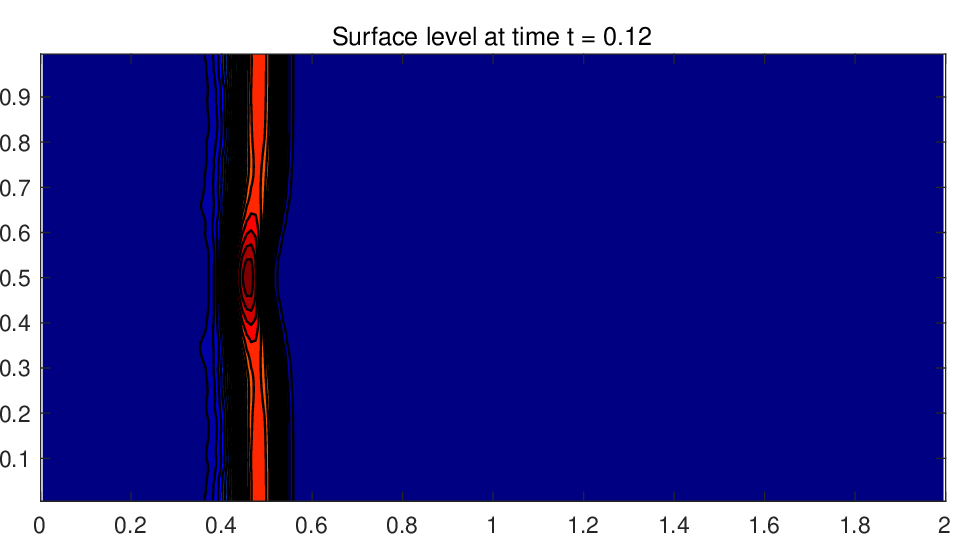}
  }
  {\centering
     \includegraphics[width= 5.5cm,scale=1]{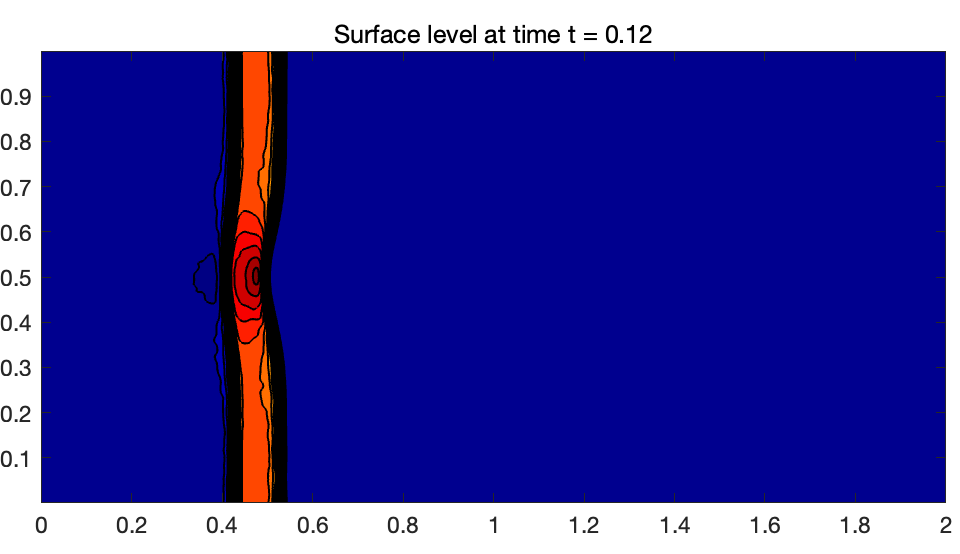}
  } 
  {\centering
     \includegraphics[width= 5.5cm,scale=1]{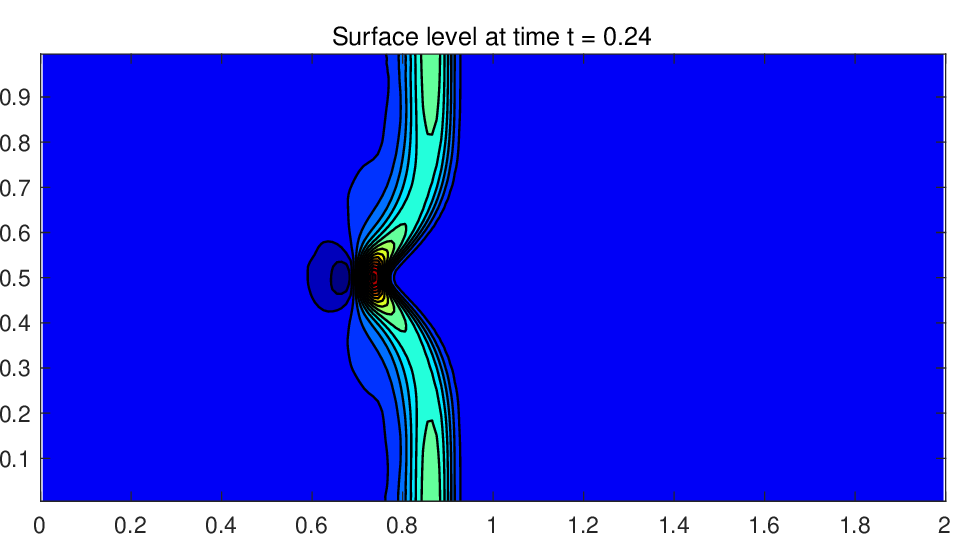}
  }
  {\centering
     \includegraphics[width= 5.5cm,scale=1]{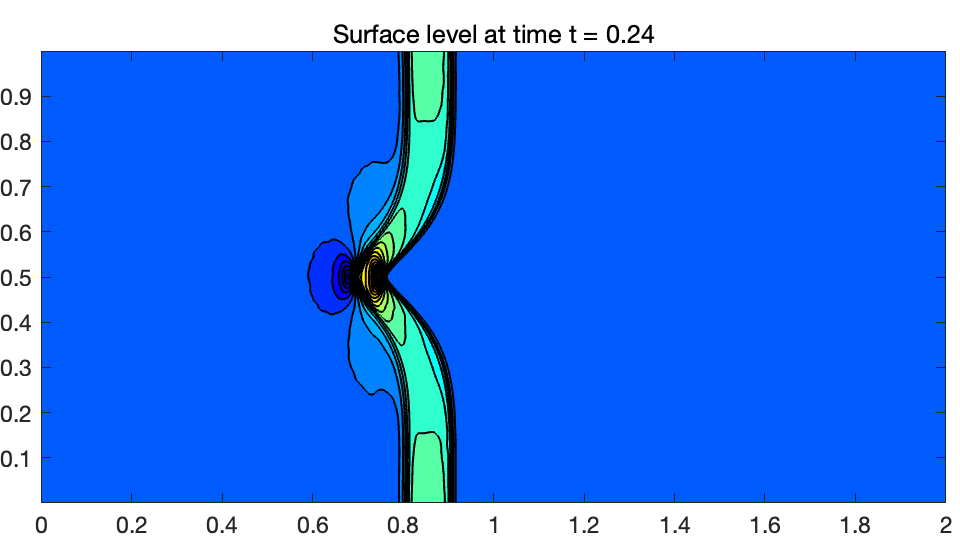}
  } 
  {\centering
     \includegraphics[width= 5.5cm,scale=1]{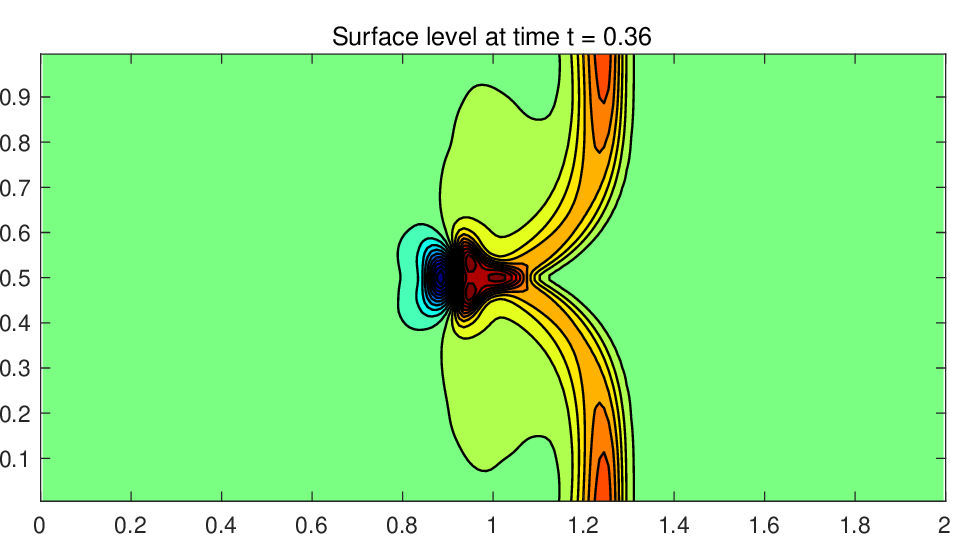}
  }
  {\centering
     \includegraphics[width= 5.5cm,scale=1]{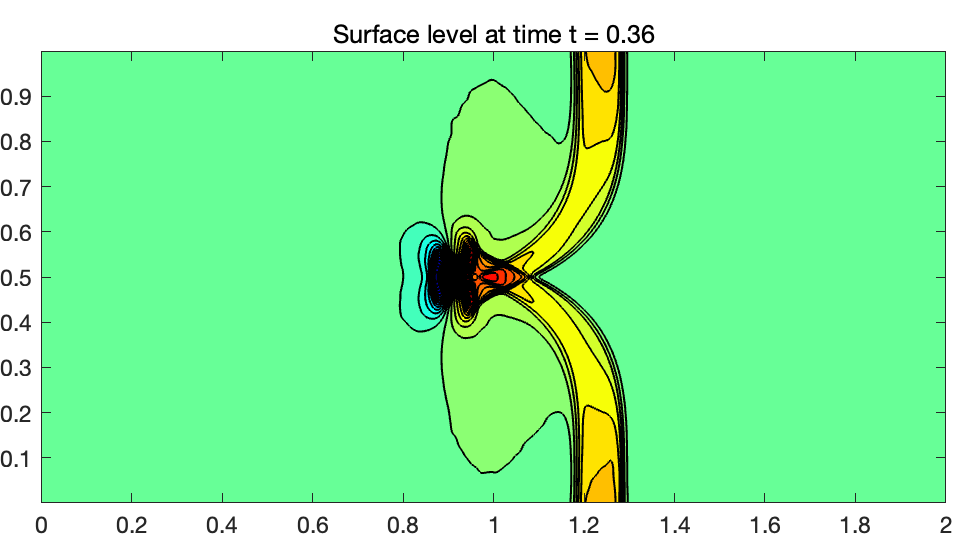}
  }
  {\centering
     \includegraphics[width= 5.5cm,scale=1]{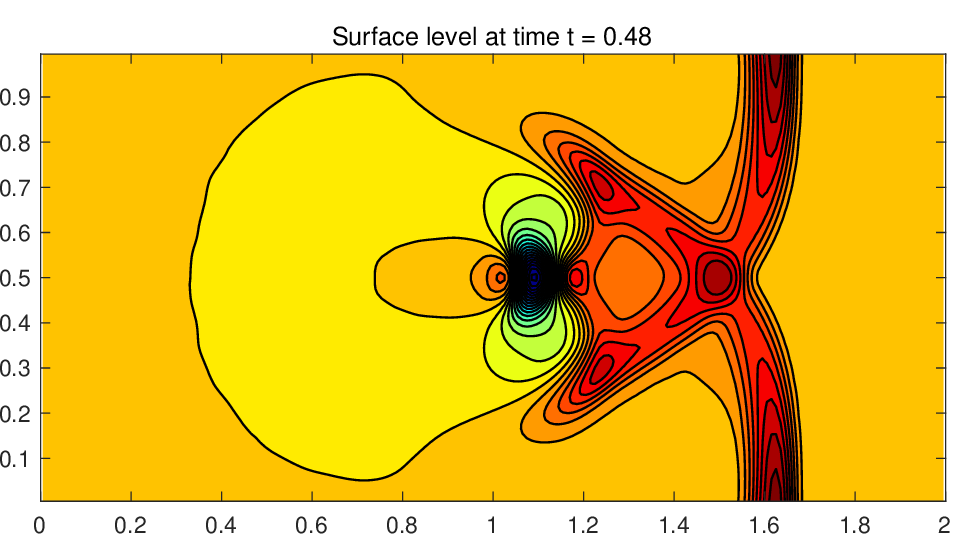}
  }
  {\centering
     \includegraphics[width= 5.5cm,scale=1]{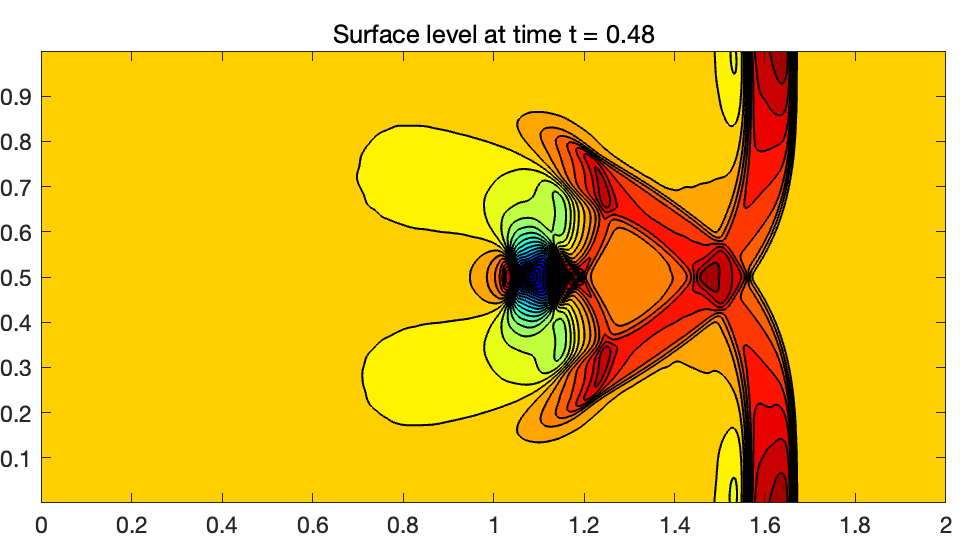}
  }
  {\centering
     \includegraphics[width= 5.5cm,scale=1]{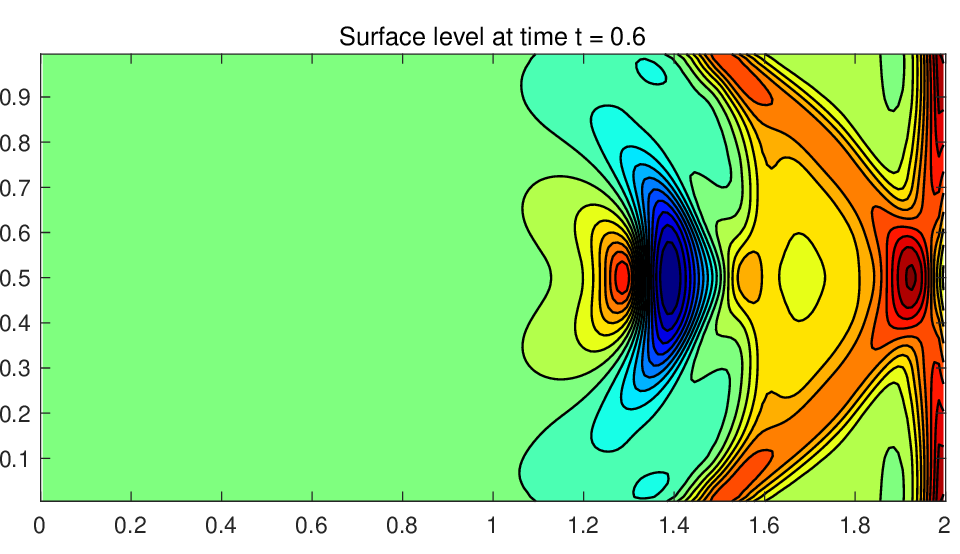}
  }
  {\centering
     \includegraphics[width= 5.5cm,scale=1]{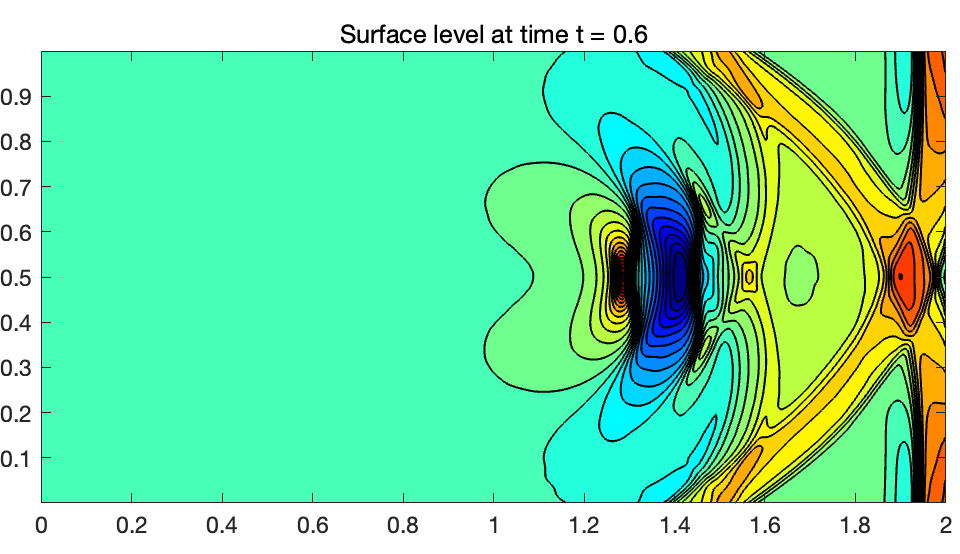}
  }
 \caption{Example \ref{2dswe:still}: The contours of the surface level with 30 uniformly spaced contour lines. From top to bottom: at time $t=0.12$ from 0.99942 to 1.00656;  at time $t=0.24$ from 0.99318 to 1.01659;  at time $t=0.36$ from 0.98814 to 1.01161;  at time $t=0.48$ from 0.99023 to 1.00508;  at time $t=0.6$ from 0.99503 to 1.00629.  Left: 200$\times$100 cells. Right: 600$\times$300 cells.}\label{per_still2d}
\end{figure}

\begin{example}{\bf A small perturbation of the still-water equilibrium}\label{2dswe:still}
\end{example}
We adopt the classical example provided by LeVeque \cite{leveque1998balancing} to demonstrate the capability of our developed well-balanced finite volume WENO scheme in handling perturbations of the still-water steady state in the two-dimensional case. The initial conditions on the computational domain $[0,2]\times [0,1]$ are given by 
  \begin{align*}
    &h(x,y,0) = \left\{\begin{array}{lll}
    1-b(x,y)+A, & \text{if} \ \ 0.05\leq x \leq 0.15, \\
    1-b(x,y),&\text{otherwise},\\
    \end{array}\right. \\
    &(hu)(x,y,0)=(hv)(x,y,0)=0,
  \end{align*}
together with the elliptical hump bottom topography
$$b(x,y)=0.8e^{-5(x-0.9)^2-50(y-0.5)^2}.$$
We simulate for a small perturbation $A= 0.01$ at different times, with 200$\times$100 and 600$\times$300 uniform mesh cells for comparison.  Fig. \ref{per_still2d} illustrates the right-going disturbance as it propagates past the hump, with the contour lines of the surface level $h+b$ revealing intricate small features of the flow with remarkable clarity.

\begin{example}{\bf Pseudo two-dimensional moving-water equilibrium}\label{2dswe:moving}
\end{example} 
In this example, we test the well-balanced property of the nodal FV-WENO scheme with flux localization in the two-dimensional case. To this end, we consider the same pseudo two-dimensional moving-water equilibrium as in \cite{noelle2007high}, where all variables along the $y$-direction are constant ($v=0$). 
The bottom topography defined on the computational domain [0,25]$\times$[0,25] is only a function of $x$, given by
\begin{equation*}
  b(x)=\left\{\begin{array}{lll}
    0.2 - 0.05(x-10)^2,&   \text{if} \ \   8\leq x \leq 12, \\
    0,&\text{otherwise}.\\
    \end{array}\right.
\end{equation*}
The initial steady state solution is computed from 
\begin{equation*}
  \frac{1}{2}u^2+g(h+b)=22.06605,\ hu(x,y,0) = 4.42,\ hv(x,y,0) = 0.
\end{equation*}
We impose the transmissive boundary condition and simulate up to time $t=5$ with 20$\times$20 and 100$\times$100 uniform cells. The $L^1$ and $L^{\infty}$ errors at double precision for the water height $h$ and discharges $hu,hv$ are listed in Table \ref{pse_wb2d}. It is evident that the numerical errors are at the round-off level, which verify the exact preservation of the steady state solution for the nodal FV-WENO scheme we developed. Similar to the results present in \cite{xu2025well}, the $L^1$ errors appear somewhat larger since the area chosen here is large. 

\begin{table}[!htb]
  \centering
  \caption{Example \ref{2dswe:moving}: $L^1$ and $L^{\infty}$ errors for the pseudo two-dimensional moving water equilibrium.}\label{pse_wb2d}
  \begin{tabular}{c c c c c c c c c c c c}
   \toprule
    \multirow{2}{*} {Mesh}&\multicolumn{3}{c}{$L^1$ error}&\multicolumn{3}{c}{$L^{\infty}$ error}\\
    \cmidrule(lr){2-4} \cmidrule(lr){5-7}
    &{$h+b$}&{$hu$}&{$hv$}&{$h+b$}&{$hu$}&{$hv$}\\
    \midrule
    20$\times$20   & 3.26e-13 & 4.22e-12 & 1.18e-12 & 1.55e-15 & 1.24e-14 & 3.13e-15 \\
    100$\times$100 & 4.09e-12 & 1.90e-11 & 2.60e-12 & 1.02e-14 & 6.22e-14 & 1.85e-14 \\
     \bottomrule
  \end{tabular}
\end{table}

\begin{figure}[htb]
  \centering 
  \includegraphics[width=4.8cm,scale=1]{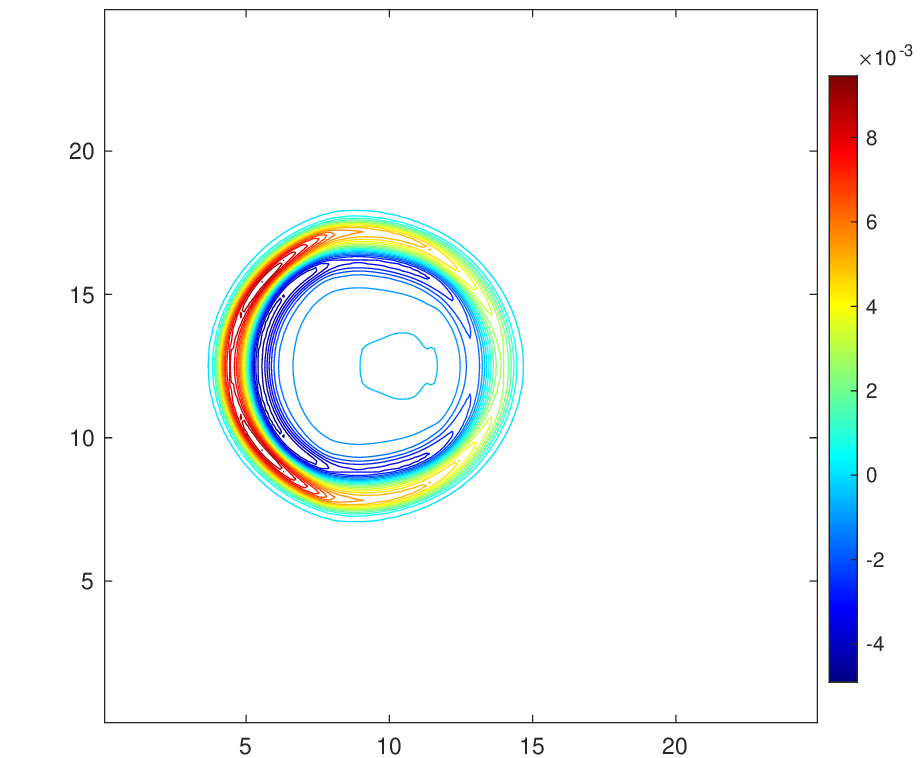} \quad
  \includegraphics[width=4.8cm,scale=1]{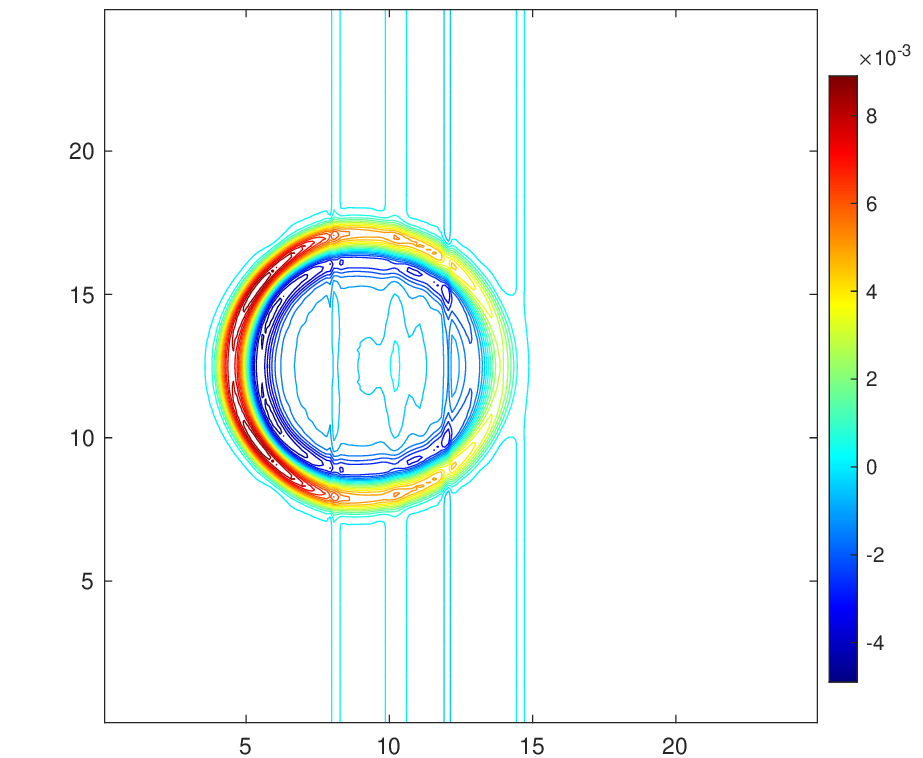}
  \caption{Example \ref{2dsper:moving}: The contours of the difference between the numerical water height and the original moving equilibrium, with 30 uniformly spaced contour lines from -0.006 to 0.01. Left: well-balanced scheme, right: non well-balanced scheme.}\label{2dper:contour}
\end{figure}

\begin{figure}[htb]
  \centering 
  \includegraphics[width=4.8cm,scale=1]{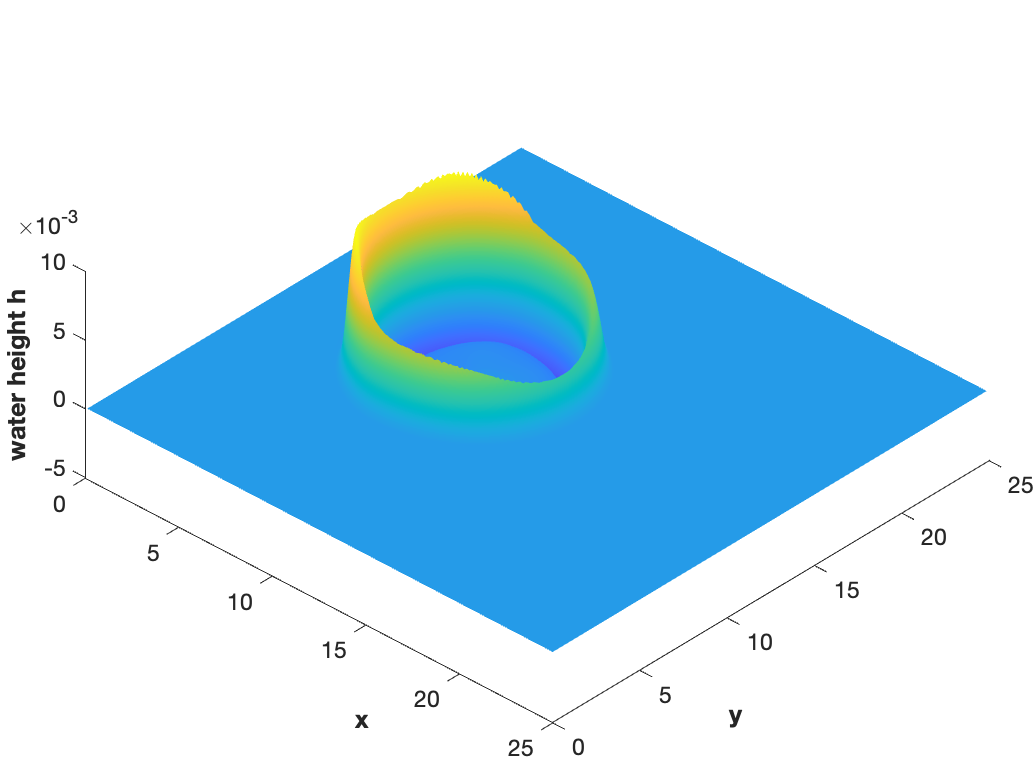} \quad
  \includegraphics[width=4.8cm,scale=1]{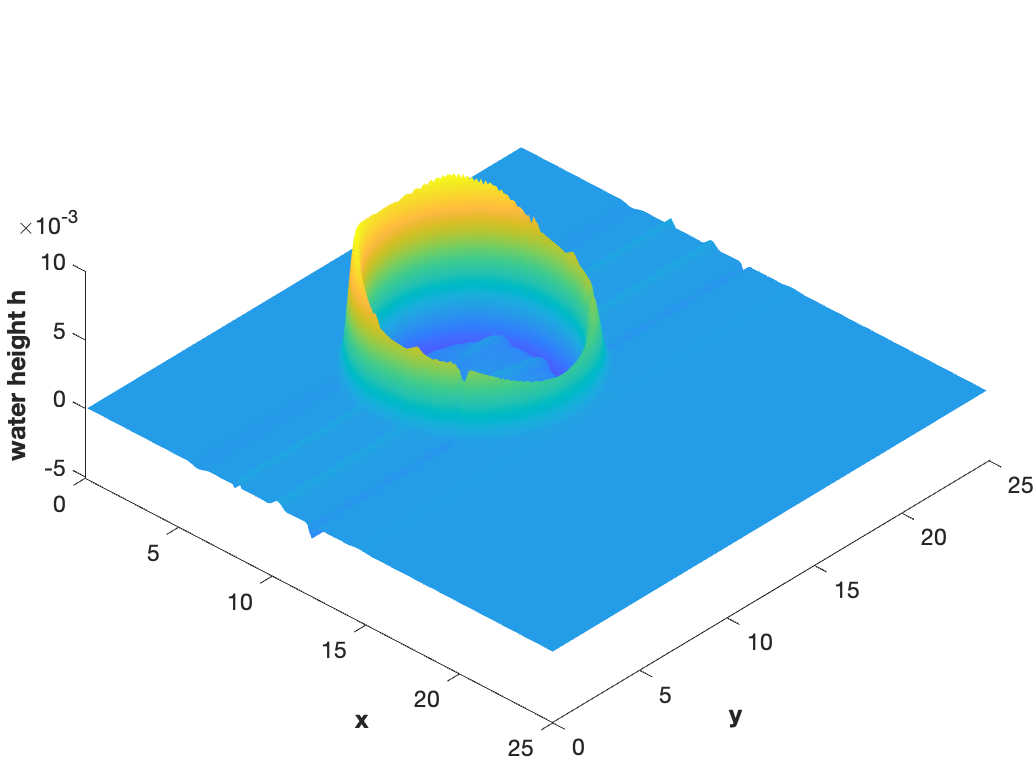}
  \caption{Example \ref{2dsper:moving}: The 3-D view of the difference between the numerical water height and the original moving equilibrium at time $t=1$. Left: well-balanced scheme, right: non well-balanced scheme.}\label{2dper:3dview}
\end{figure}

\begin{example}{\bf Small perturbation of pseudo two-dimensional moving-water equilibrium}\label{2dsper:moving}
\end{example} 
In the final example, we demonstrate the ability to capture the propagation of small perturbations in pseudo two-dimensional moving-water equilibrium, a capability that cannot be realized by a non-well-balanced scheme. Using almost the same setup as in Example \ref{2dswe:moving}, we introduce a small disturbance $A = 0.05$ to the initial water height $h$ in the box [6.5,7.5]$\times$[12,13]. We present the disturbance as it interacts with the hump using 200$\times$200 uniform cells at time $t=1$, and show the difference between the numerical water height and the original moving equilibrium in Fig. \ref{2dper:contour}. For comparison, we also perform the same numerical test using the well-balanced fifth-order finite volume WENO scheme for still-water equilibrium \cite{xing2006new}. The scheme designed for still-water steady-state \cite{xing2006new} is non-well-balanced for moving-water equilibrium. 
Fig. \ref{2dper:3dview} displays the 3D view of the difference for both the well-balanced and non-well-balanced schemes. The results reveal that our well-balanced scheme effectively captures the intricate and minute characteristics of the flow, while the non-well-balanced scheme produces spurious features, compromising the solution's accuracy.

\section{Conclusion}\label{se:co}
We have designed and tested a novel high-order well-balanced nodal finite volume WENO scheme for various systems of hyperbolic balance laws. The key difference from classical finite volume methods lies in evolving the nodal point values rather than the cell averages, which facilitates the transformation between the conservative and equilibrium variables. High-order accuracy is naturally achieved through the WENO reconstruction procedure applied to the equilibrium variables.
Another important technique for achieving the well-balanced property is the flux localization of the source terms. Moreover, the robustness of this method, particularly in the presence of non-conservative products when the bottom topography is discontinuous, is ensured through the application of the linear segment path function to the equilibrium variables. Both the theoretical analysis and numerical examples confirm the well-balanced nature of our proposed scheme. When extended to higher dimensions, the current scheme is limited to a family of quasi one-dimensional steady states. In future work, we aim to explore finite volume WENO schemes or discontinuous Galerkin methods capable of preserving genuinely two-dimensional equilibrium.

\begin{appendix}
\section{WENO reconstruction from nodal values}\label{a1}
For simplicity of presentation, we only describe the fifth-order WENO  reconstruction procedure for the proposed nodal FV-WENO scheme. Specifically, we detail the computation of the left-sided value $U_{j+\frac{1}{2}}^-$ from the ``nodal" values $\{U_j\}$, noting that the right-sided value $U_{j+\frac{1}{2}}^+$ can be obtained using the mirror-symmetric stencil.

Based on the small stencils 
$$s_0(j) = \{I_{j-2}, I_{j-1}, I_{j}\},s_1(j) = \{I_{j-1}, I_{j}, I_{j+1}\},s_2(j) = \{I_{j}, I_{j+1}, I_{j+2}\}, $$
we construct third-order accurate polynomials whose values at the cell interface $x_{j+\frac 1 2}$ are given by
\begin{align*}
  &u_0(x_{j+\frac{1}{2}}) = \frac{1}{8}(3 U_{j-2}-10U_{j-1}+15U_{j}),\\
  &u_1(x_{j+\frac{1}{2}}) = \frac{1}{8}(- U_{j-1}+6U_{j}+3U_{j+1}),\\
  &u_2(x_{j+\frac{1}{2}}) = \frac{1}{8}(3 U_{j}+6U_{j+1}-U_{j+2}),
\end{align*}
To obtain the linear weights, we use the large stencil 
$$S_j=\{I_{j-2}, I_{j-1}, I_{j},I_{j+1}, I_{j+2}\},$$
which yields the linear weights 
$$d_0 = \frac{1}{16}, \; d_1 = \frac{5}{8}, \; d_2 = \frac{5}{16}.$$
Next, we compute the smoothness indicators
\begin{align*}
\beta_r=\sum \limits _{l=1}^2 (\Delta x)^{2l-1} \displaystyle{\int_{x_{j-{1}/{2}}}^{x_{j+{1}/{2}}} \  \left( \dfrac{\partial^l u_r(x)}{\partial x^l} \right)^2 \ dx}, \ r=0,1,2,
\end{align*}
which measure the relative smoothness of the solution within each sub-stencil $s_r(j)$. For fifth-order accuracy, the smoothness indicators have the following explicit forms:
\begin{equation*}
  \begin{aligned}
     &\beta_0=\dfrac{13}{12}( U_{j-2}-2 U_{j-1}+ U_{j})^2+\dfrac{1}{4}(U_{j-2}-4 U_{j-1}+3 U_{j})^2.\\
     &\beta_1=\dfrac{13}{12}( U_{j-1}-2 U_{j}+ U_{j+1})^2+\dfrac{1}{4}( U_{j-1}-U_{j+1})^2.\\
     &\beta_2=\dfrac{13}{12}( U_{j}-2 U_{j+1}+ U_{j+2})^2+\dfrac{1}{4}(3U_{j}-4U_{j+1}+ U_{j+2})^2.
  \end{aligned}
 \end{equation*}
The nonlinear weights $\omega_r$ are then defined by
\begin{equation*}
  \omega_r =\dfrac{\alpha_r}{ \alpha_0 +\alpha_1+\alpha_2},\quad \alpha_r=\dfrac{d_r}{(\beta_r+\varepsilon)^2}, \quad r = 0,1,2.
\end{equation*}
Here $\varepsilon>0$ is a small parameter introduced to avoid division by zero. In our numerical experiments, we set $\varepsilon = 10^{-12}$ unless otherwise specified.
Finally, the fifth-order WENO reconstruction from nodal values is given by the convex combination
\begin{equation*}
U_{j+\frac{1}{2}}^- =\sum \limits _{r=0}^2\omega_r u_r(x_{j+\frac{1}{2}}).
\end{equation*}

\section{Equilibrium variables for two-dimensional Euler systems and Ripa model}\label{b1}
Two-dimensional Euler equations with gravity and Ripa model can be dealt with in the similar fashion as in the Section \ref{se:e2d}, with the modification of the relevant equilibrium variables and the matrix in \eqref{sse_2d}. Here, we present the specific expressions. 

For the Euler equations with gravity, the matrix and the equilibrium variable vectors in the $x$- and $y$-direction are defined as 
\begin{equation*}
  \mathcal{L}^x = \left[\begin{array}{ccccc}
    0 & 1 & 0 & 0 & 0 \\
    -\frac{\rho^{\gamma}}{\gamma-1} & u & 0 & \rho & 0\\ 
    0 & v & \rho u & 0 & 0\\
    0 & \frac{\gamma K\rho^{\gamma-1}}{\gamma-1}+\frac{u^2}{2} & 0 & \rho u & 0\\
    0 & 0 & 0 & 0 & 1 \\
  \end{array}\right], \ 
  \mathcal{L}^y = \left[\begin{array}{ccccc}
    0 & 0 & 1 & 0 & 0 \\
    0 & \rho v & u & 0 & 0 \\ 
    -\frac{\rho^{\gamma}}{\gamma-1} & 0 & v & \rho & 0 \\
    0 & 0 & \frac{\gamma K\rho^{\gamma-1}}{\gamma-1}+\frac{v^2}{2}  & \rho v & 0\\
    0 & 0 & 0 & 0 & 1 \\
  \end{array}\right],  
\end{equation*}
and
\begin{equation*}
  \widetilde{\boldsymbol{v}^x} = (\boldsymbol{v}^x,0)^T = (K,\rho u,v,\mathcal{E},0)^T,\quad \widetilde{\boldsymbol{v}^y} = (\boldsymbol{v}^y,0)^T = (K,u,\rho v,E,0)^T.
\end{equation*}
Here, $\mathcal{E} = \frac{1}{2}u^2+h+\phi$ and $E = \frac{1}{2}v^2+h+\phi$.

For the Ripa model, 
\begin{equation*}
  \mathcal{L}^x = \left[\begin{array}{ccccc}
    0 & 1 & 0 & 0 & 0 \\
    h & u & 0 & -\frac{1}{2}gh^2+ghb & 0\\ 
    0 & v & hu & 0 & 0\\
    0 & \theta & 0 & hu & 0\\
    0 & 0 & 0 & 0 & 1 \\
  \end{array}\right], \ 
  \mathcal{L}^y = \left[\begin{array}{ccccc}
    0 & 0 & 1 & 0 & 0 \\
    0 & hv & u & 0 & 0 \\ 
    h & 0 & v & -\frac{1}{2}gh^2+ghb & 0 \\
    0 & 0 & \theta  & hv & 0\\
    0 & 0 & 0 & 0 & 1 \\
  \end{array}\right],  
\end{equation*}
and
\begin{equation*}
  \widetilde{\boldsymbol{v}^x} = (\boldsymbol{v}^x,0)^T = (\mathcal{E},hu,v,\theta,0)^T,\quad \widetilde{\boldsymbol{v}^y} = (\boldsymbol{v}^y,0)^T = (E,u,hv,\theta,0)^T.
\end{equation*} 
Here, $\mathcal{E} = \frac{1}{2}u^2+g\theta(h+b)$ and $E = \frac{1}{2}v^2+g\theta(h+b)$.

\section{Numerical examples for scalar conservation law}\label{c1}
\begin{figure}[htb]
  \centering 
  \includegraphics[width=5.cm,scale=1]{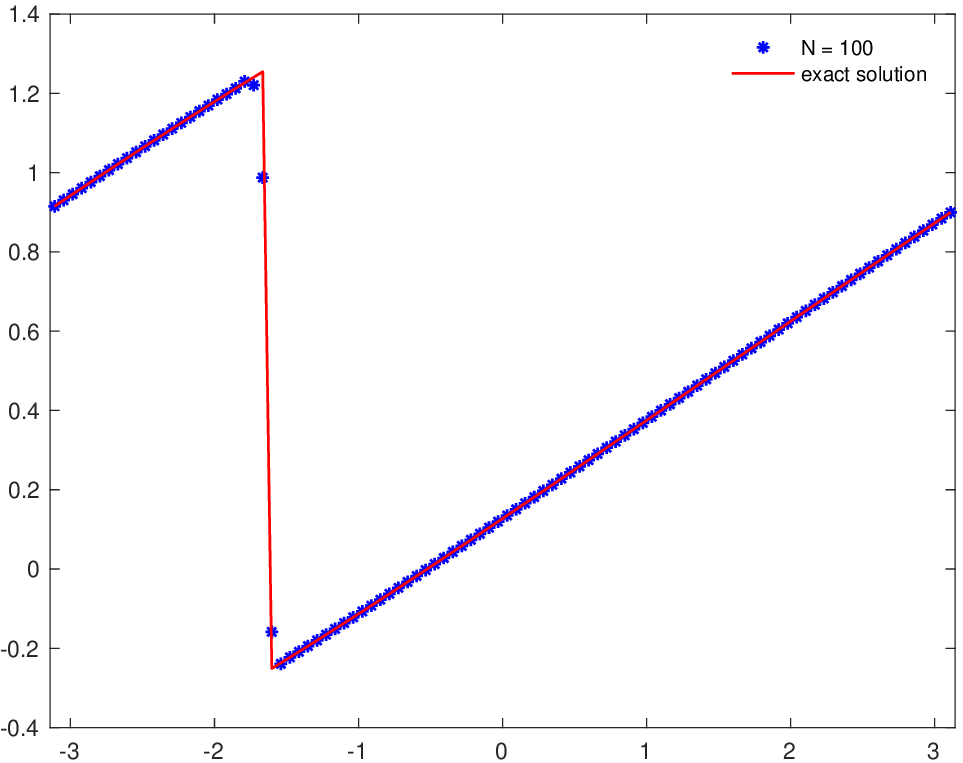}
  \caption{Example \ref{example:burgers}: Burgers equations with a shock, the numerical and analytical solutions at time $t=3$. }\label{fig:burgers}
\end{figure}

To assess the performance of our proposed FV-WENO scheme in the presence of discontinuities, we further consider several scalar conservation laws.
\begin{example}{\bf Burgers equations}\label{example:burgers}
\end{example}
We solve the scalar Burgers equations 
\begin{equation*}
  U_t + \left(\frac{U^2}{2}\right)_x = 0,
\end{equation*}
on the computational domain $[-\pi,\pi]$ together with the periodic boundary conditions. The initial condition is given by $U(x,0) = 0.5+\sin(x)$. We perform the simulation with 100 uniform cells up to the final time $t=3$ when the shock wave has formed in the solution. Comparisons of numerical solution computed by our proposed fifth-order FV-WENO scheme and the analytical solution are shown in Fig. \ref{fig:burgers}. We can observe that the discontinuity location is well captured without numerical oscillations.

\begin{figure}[htb]
  \centering 
  \subfigure[$t=0.1$]{
  \centering
     \includegraphics[width= 5.cm,scale=1]{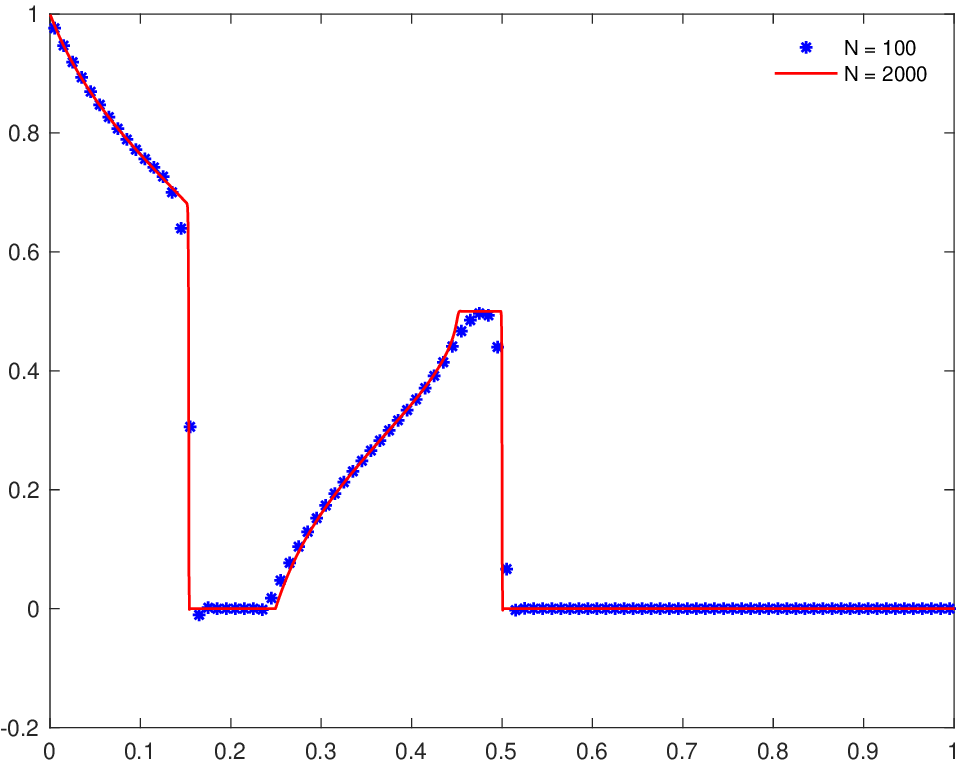}
  } 
  \subfigure[$t=0.2$]{
  \centering
     \includegraphics[width= 5.cm,scale=1]{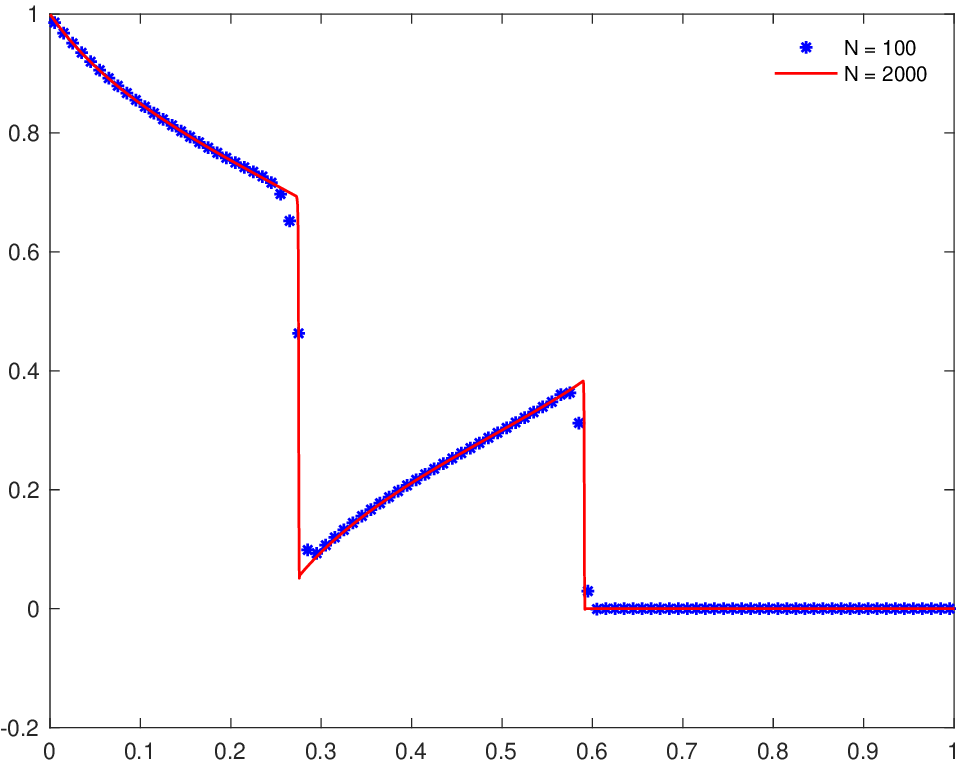}
  } 
  \subfigure[$t=0.4$]{
  \centering
     \includegraphics[width= 5.cm,scale=1]{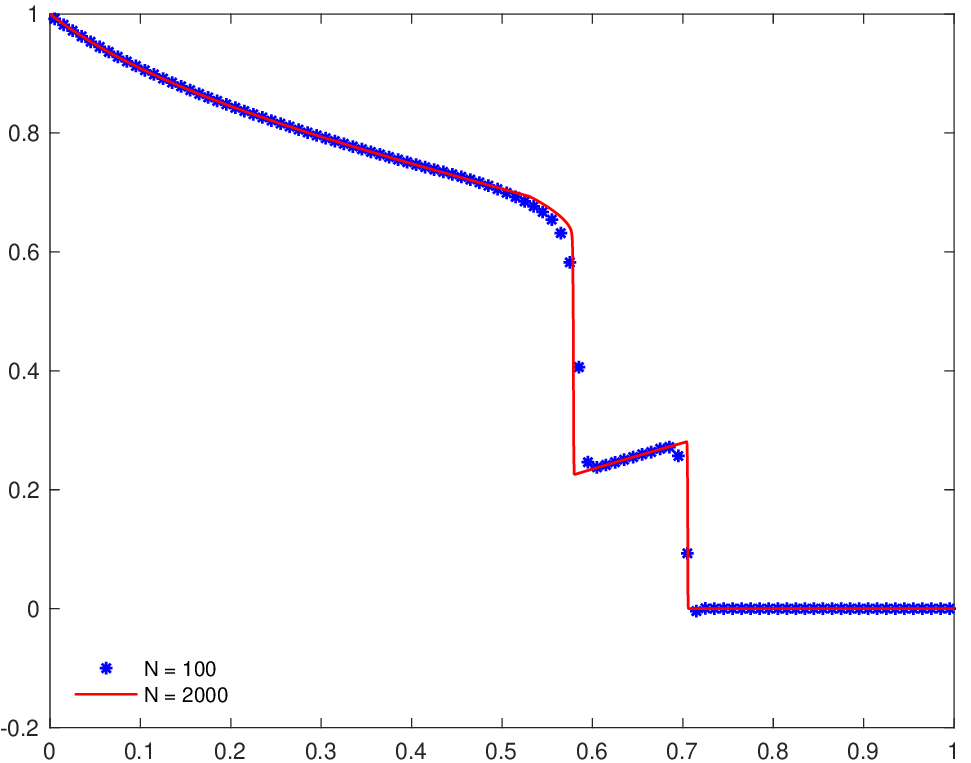}
  } 
  \subfigure[$t=0.5$]{
  \centering
     \includegraphics[width= 5.cm,scale=1]{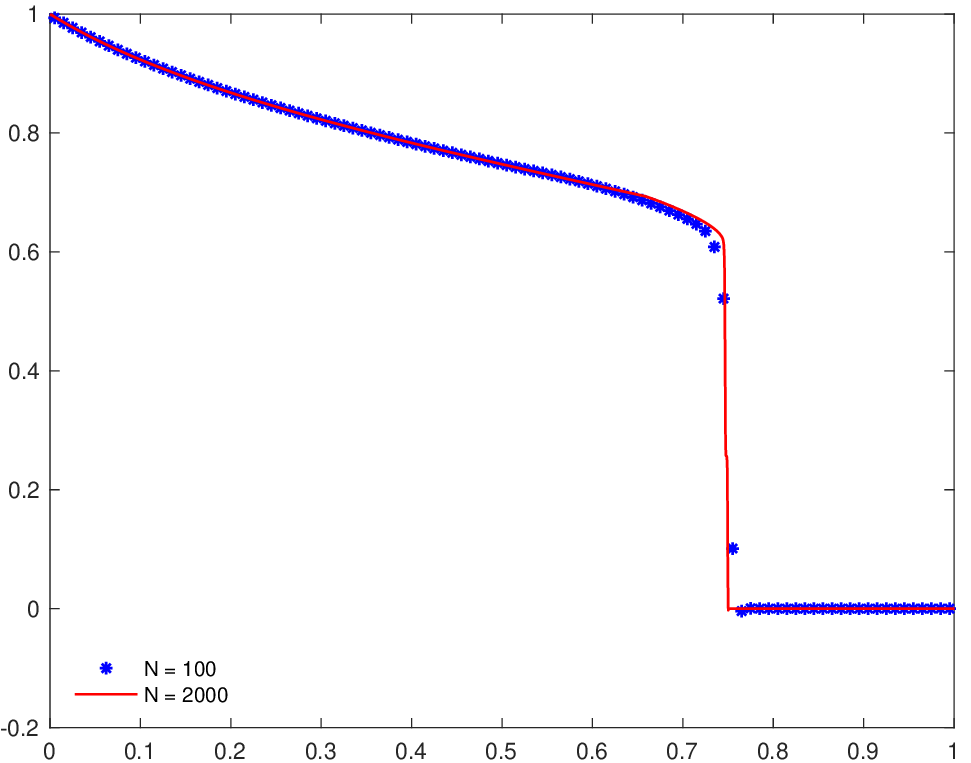}
  } 
  \caption{Example \ref{example:BL}: The numerical solution for the Buckley-Leverett equations at different times $t=0.1, 0.2, 0.4, 0.5$.}\label{fig:BL}
\end{figure}

\begin{example}{\bf Buckley-Leverett equations}\label{example:BL}
\end{example}
We consider the Buckley-Leverett equations with a non-convex flux,
\begin{equation*}
  U_t + \left(\frac{U^2}{U^2+(1-U)^2}\right)_x = 0.
\end{equation*}
The initial condition is given by
\begin{equation*}
  U(x,0) = \left\{\begin{array}{lll}
    1-20x,&   \text{if} \ \   0\leq x \leq 0.05, \\
    0.5,&   \text{if} \ \   0.25\leq x \leq 0.4, \\
    0,&\text{otherwise},\\
    \end{array}\right.
\end{equation*}
on the computational domain $[0,1]$. The evolution of the two initial pulses gives rise to shocks and rarefaction waves, which subsequently interact and merge. The reference solutions are calculated by the classical fifth-order WENO scheme \cite{shu1998essentially} on a refined 2000 uniform cells. Fig. \ref{fig:BL} shows our numerical solutions obtained on a uniform mesh of 100 cells at different times, which are in excellent agreement with the reference solution.

\end{appendix}


\addcontentsline{toc}{section}{References}
\bibliographystyle{abbrv}
\normalem

\bibliography{shortreference}
\end{document}